\documentclass[11pt]{article}
\usepackage{amsfonts}
\usepackage{amsfonts,latexsym,amsmath,amscd,geometry}
\usepackage{amssymb}
\usepackage{latexsym}
\renewcommand{\theequation}{\thesection.\arabic{equation}}
\newcommand \nc{\newcommand}
\newtheorem{theorem}{Theorem}[section]
\newtheorem{lemma}[theorem]{Lemma}

\newtheorem{corollary}[theorem]{Corollary}
\newtheorem{definition}[theorem]{Definition}

\newtheorem{remark}[theorem]{Remark}
\renewcommand{\thetheorem}{\thesubsection.\arabic{theorem}}
\nc{\ba}{\begin{array}}\nc{\ea}{\end{array}}
\nc{\be}{\begin{eqnarray}}\nc{\ee}{\end{eqnarray}}
\nc{\beq}{\begin{equation}}\nc{\eeq}{\end{equation}}
\nc{\bex}{\begin{eqnarray*}}\nc{\eex}{\end{eqnarray*}}
\nc{\btm}{\begin{theorem}} \nc{\etm}{\end{theorem}}
\nc{\blm}{\begin{lemma}} \nc{\elm}{\end{lemma}}
\nc{\R}{\mathbb{R}} \nc{\va}{\varepsilon} \nc{\ls}{\limits}

\def\de{\Delta}
\def\pf{\noindent{\bf Proof.\quad}}\def\endpf{\hfill$\Box$}
\def\les{\lesssim}\def\u{\dot{u}}\def\di{\mbox{div\,}}

\newcommand \qed {\hfill $\Box$}

\begin{document}
\title{{\bf Global symmetric classical and strong solutions of the full compressible Navier-Stokes equations with vacuum and large initial
data}}
\author{
\begin{tabular}{cc}
& Huanyao Wen, \ \ \ Changjiang Zhu\thanks{Corresponding author.\ \
Email:
cjzhu@mail.ccnu.edu.cn} \\[4mm]
 &  School of Mathematics and
Statistics \\
& Central China Normal University, Wuhan 430079, China\\
\end{tabular}
}
\date{}

\maketitle

\begin{abstract}
First of all, we get the global existence of classical and strong
solutions of the full compressible Navier-Stokes equations in three
space dimensions with initial data which is large and spherically or
cylindrically symmetric. The appearance of vacuum is allowed. In
particular, if the initial data is spherically symmetric, the space
dimension can be taken not less than two. The analysis is based on
some delicate {\it a priori} estimates globally in time which depend
on the assumption $\kappa=O(1+\theta^q)$ where $q>r$ ($r$ can be
zero), which relaxes the condition $q\ge2+2r$ in
\cite{Feireisl-book,Kawohl, Wen-Zhu}. This could be viewed as an
extensive work of \cite{Hoff} where the equations hold in the sense
of distributions in the set where the density is positive with
initial data which is large, discontinuous, and spherically or
cylindrically symmetric in three space dimension. Finally, with the
assumptions that vacuum may appear and that the solutions are not
necessarily symmetric, we establish a blow-up criterion in terms of
$\|\rho\|_{L^\infty_tL_x^\infty}$ and
$\|\rho\theta\|_{L^4_tL^\frac{12}{5}_x}$ for strong solutions.
\end{abstract}


\tableofcontents

\vspace{4mm}
\section {Introduction}
\setcounter{equation}{0} \setcounter{theorem}{0}
 The full compressible Navier-Stokes
equations can be written in the
sense of Eulerian coordinates in $\Omega\subset\mathbb{R}^N$ as
follows: \be\label{full N-S}
\begin{cases}
\rho_t+\nabla \cdot (\rho {\bf u})=0, \\
(\rho {\bf u})_t+\mathrm{div}(\rho{\bf u}\otimes{\bf u})+\nabla
P=\mathrm{div}(\mathcal {T})+\rho{\bf f},\\
(\rho E)_t+\mathrm{div}(\rho E{\bf u})+\mathrm{div}(P{\bf
u})=\mathrm{div}(\mathcal {T}{\bf
u})+\mathrm{div}(\kappa\nabla\theta)+\rho {\bf u}\cdot{\bf f}.
\end{cases}
\ee Here $\mathcal {T}$ is the stress tensor given by
$$\mathcal
{T}=\mu\left(\nabla {\bf u}+(\nabla {\bf u})^\prime\right)+\lambda
\di {\bf u} I_N,$$ where $I_N$ is a $N\times N$ unit matrix;
 $\rho=\rho({\bf x},t)$, ${\bf u}={\bf u}({\bf
x}, t)=(u_1,\cdots, u_N)({\bf x},t)$ and $\theta=\theta({\bf x},t)$
are unknown functions denoting  the density, velocity and absolute
temperature, respectively; $P=P(\rho,\theta)$, $E$, ${\bf f}={\bf
f}({\bf x}, t)=(f_1,\cdots, f_N)({\bf x},t)$ and $\kappa$ denote
respectively pressure, total energy, external forces and coefficient
of heat conduction, where $E=e+\frac{1}{2}u^2$ ($e$ is the internal
energy);
 $\mu$ and $\lambda$ are coefficients of viscosity, satisfying the following physical restrictions:
$$
\mu>0,\ 2\mu+N\lambda\ge0;
$$ $P$ and $e$ satisfy the second principle
of thermodynamics: \be\label{non-1.1} P=\rho^2\frac{\partial
e}{\partial\rho}+\theta\frac{\partial P}{\partial \theta}. \ee
 (\ref{full N-S}) is a well-known model which describes the
motion of compressible fluids. There were lots of works on the
existence,
 uniqueness, regularity and asymptotic behavior of the solutions
  during the last five decades. While, because of the stronger nonlinearity in
(\ref{full N-S}) compared with the Navier-Stokes equations for
isentropic flow (no temperature equation), many known mathematical
results focused on the absence of vacuum (vacuum means $\rho=0$),
refer for instance to \cite{Itaya, Jiang1, Kawohl,
Kazhikhov-Shelukhi, Matsumura-Nishida: Kyoto Un,
Matsumura-Nishida: CMP, Tani} for classical solutions. More
precisely, the local classical solutions to the Navier-Stokes
equations with heat-conducting fluid in H\"older spaces was
obtained respectively by Itaya in \cite{Itaya} for Cauchy problem
and by Tani in
 \cite{Tani} for IBVP with $\inf\rho_0>0$, where the space dimension $N=3$. Using
delicate energy methods in Sobolev spaces, Matsumura and Nishida in
\cite{Matsumura-Nishida: Kyoto Un, Matsumura-Nishida: CMP} showed
that the global classical solutions exist provided that the initial
data is small in some sense and away from vacuum in three space
dimension. For large initial data in one space dimension, Kazhikhov,
Shelukhi in \cite{Kazhikhov-Shelukhi} (for polytropic perfect gas
with $\mu,\lambda,\kappa=$const.) and Kawohl in \cite{Kawohl} (for
real gas with $\kappa=\kappa(\rho,\theta),\
\mu,\lambda=\mathrm{const.}$) respectively got the global classical
solutions to (\ref{full N-S}) in Lagrangian coordinates with
$\inf\rho_0>0$. The internal energy $e$ and the coefficient of heat
conduction $\kappa$ in \cite{ Kawohl} satisfy the following
assumptions for $\rho\le\overline{\varrho}$ and $\theta\ge0$ (we
translate these conditions in Eulerian coordinates)
\be\label{non-r1.5}\begin{cases} e(\rho,0)\ge0,\ \
\nu(1+\theta^r)\le\partial_\theta e(\rho,\theta)\le
N(\overline{\varrho})(1+\theta^r),\\
\kappa_0(1+\theta^q)\le\kappa(\rho,\theta)\le\kappa_1(1+\theta^q),\\
|\partial_\rho\kappa(\rho,
\theta)|+|\partial_{\rho\rho}\kappa(\rho,
\theta)|\le\kappa_1(1+\theta^q),
\end{cases}\ee
where $r\in[0,1]$, $q\ge2+2r$, and $\nu$, $N(\overline{\varrho})$,
$\kappa_0$ and $\kappa_1$ are positive constants. For the perfect
gas (i.e.,
$P=R\rho\theta$, $e=C_\nu \theta$ for some constants $R>0$ and
$C_\nu>0$) in the domain exterior to a ball in $\mathbb{R}^N$ ($N = 2$ or
$3$) with $\mu,\lambda,\kappa=$const., Jiang in \cite{Jiang1} got
the existence of global spherically symmetric classical large
solutions in H\"older spaces.

In fact, Kawohl in \cite{ Kawohl} also considered the case of
density-dependent viscosity for another boundary condition with
$\inf\rho_0>0$, where
$0<\underline{\mu}_0\le\mu(\rho)\le\overline{\mu}_0$ for any
$\rho\ge0$, and $\underline{\mu}_0$ and $\overline{\mu}_0$ are
positive constants. This was generalized to the case
$\mu(\rho)=\rho^\alpha$ by Jiang in \cite{Jiang:Math Nachr} for
$\alpha\in(0,\frac{1}{4})$, and by Qin, Yao in \cite{Qin-Yao} for
$\alpha\in(0,\frac{1}{2})$, respectively.

On the existence, asymptotic behavior of the weak solutions of
the full compressible Navier-Stokes equations with $\inf\rho_0>0$, please refer for instance to \cite{
Jiang2, Jiang3, Jiang-Zhang:weak solutions} for the existence of
weak solutions in 1D and for the existence of spherically
symmetric weak solutions in $\mathbb{R}^N$ ($N=2$, $3$), and refer
to \cite{Hoff} for the existence of spherically and cylindrically
symmetric weak solutions in $\mathbb{R}^3$, and refer to
\cite{Feireisl1} for the existence of {\em variational} solutions
in a bounded domain in $\mathbb{R}^N$ ($N=2$, $3$).

In the presence of vacuum (i.e. $\rho$ may vanish), to our best
knowledge, the mathematical results on global well-posedness of the
full compressible Navier-Stokes equations are usually limited to the
existence of weak solutions with special pressure, viscosity and
heat conductivity (see \cite{Bresch-Desjardins, Feireisl-book}).
More precisely, Feireisl in \cite{Feireisl-book} got the existence
of so-called {\em variational} solutions in dimension $N\ge 2$. The
temperature equation in \cite{Feireisl-book} is satisfied only as an
inequality in the sense of distributions. Anyhow, Feireisl's work is
the very first attempt towards the existence of weak solutions to
the full compressible Navier-Stokes equations in higher dimensions,
where the coefficients of viscosity are constants and
\be\label{non-1.6}\begin{cases} \kappa=\kappa(\theta)\in
C^2[0,\infty),\
\underline{\kappa}(1+\theta^q)\le\kappa(\theta)\le\overline{\kappa}(1+\theta^q)\
\ \mathrm{for}\ \mathrm{all}\ \theta\ge0,
\\ P=P(\rho,\theta)=\mathcal {P}_e(\rho)+\theta\mathcal{P}_\theta(\rho)\ \
\mathrm{for}\ \mathrm{all}\ \rho\ge0\ \mathrm{and}\ \theta\ge0,\\
\mathcal{P}_e, \mathcal{P}_\theta\in C[0,\infty)\cap
C^1(0,\infty);\ \mathcal{P}_e(0)=0,\ \mathcal{P}_\theta(0)=0,
\\ \mathcal{P}_e^\prime(\rho)\ge a_1\rho^{\overline{\gamma}-1}-b_1\ \ \mathrm{for}\
\mathrm{all}\ \rho>0;\ \mathcal{P}_e(\rho)\le
a_2\rho^{\overline{\gamma}}+b_1\ \
\mathrm{for}\ \mathrm{all}\ \rho\ge0,\\
\mathcal{P}_\theta \ \ \mathrm{is} \ \ \mbox{non-decreasing}\
\mathrm{in}\ [0,\infty);\ \mathcal{P}_\theta(\rho)\le
a_3(1+\rho^\Gamma)\ \ \mathrm{for}\ \mathrm{all}\ \rho\ge0,
\end{cases}\ee
where $\Gamma<\frac{\overline{\gamma}}{2}$ if $N=2$ and
$\Gamma=\frac{\overline{\gamma}}{N}$
 for $N\ge3$; $q\ge2$, $\overline{\gamma}>1$, and $a_1$, $a_2$, $a_3$, $b_1$, $\underline{\kappa}$
 and $\overline{\kappa}$ are positive constants. Note that the perfect gas
equation of state (i.e. $P=R\rho\theta$ for some constant $R>0$)
is not involved in (\ref{non-1.6}). In order that the equations
are satisfied as equalities in the sense of distribution, Bresch
and Desjardins in \cite{Bresch-Desjardins} proposed some different
assumptions from \cite{Feireisl-book}, and obtained the existence
of global weak solutions to the full compressible Navier-Stokes equations with large initial data in
$\mathbb{T}^3$ or $\mathbb{R}^3$. In \cite{Bresch-Desjardins}, the
viscosity $\mu=\mu(\rho)$ and $\lambda=\lambda(\rho)$ may vanish
when vacuum appears, and $\kappa$, $P$ and $e$ are assumed to be
satisfied\be\label{non-1.7}\begin{cases}
\kappa(\rho,\theta)=\kappa_0(\rho,\theta)(\rho+1)(\theta^q+1),\\
P=R\rho\theta+p_c(\rho), \\
e=C_\upsilon\theta+e_c(\rho),
\end{cases}\ee
where $q\ge2$, $R$ and $C_\upsilon$ are two positive constants,
$p_c(\rho)=O(\rho^{-\ell})$ and $e_c(\rho)=O(\rho^{-\ell-1})$ (for
some $\ell>1$) when $\rho$ is small enough, and
$\kappa_0(\rho,\theta)$ is assumed to satisfy
$$\underline{c}_0\le\kappa_0(\rho,\theta)\le\frac{1}{\underline{c}_0},$$
for $\underline{c}_0>0$. On the local existence and uniqueness of
strong solutions for $N=3$, please see \cite{cho-Kim: perfect gas}
for the perfect gas with $\mu,\lambda,\kappa=$const. While, there
are no global smooth solutions to (\ref{full N-S}) for Cauchy
problem when the initial density is of nontrivially compact support
and $\kappa=0$ (see \cite{Xin}).

Except for the Cauchy problems with initial density compactly
supported, it is still unknown whether the global strong (or
classical) solutions exist or not when vacuum appears (i.e., the
density may vanish) until recently. In our previous paper
\cite{Wen-Zhu}, we got existence and uniqueness of global classical
solutions to the full compressible Navier-Stokes equations in one
dimension with large initial data and vacuum. In \cite{Wen-Zhu}, the
coefficient of conduction $k$ depends on the temperature, growing as
$1+\theta^q$ where $q\ge2+2r$ ($r$ can be zero).

As a first step to study the problems (\ref{full N-S}) in high
dimensions, we study the problems in high dimensions with some
symmetry which reduces the whole system to an one dimensional system
with singular source terms. The singularity may be due to $x=0$,
$x=\infty$ or appearance of vacuum, where $x$ is the radius. Our
main concern here is to show the existence and uniqueness of global
classical and strong solutions to (\ref{full N-S}) with vacuum and
initial data which is large, and spherically or cylindrically
symmetric in three space dimension. In particular, if the initial
data is spherically symmetric, the space dimension can be taken not
less than two. This extends the results in \cite{Hoff} where the
equations hold in the sense of distributions in the set where the
density is positive with initial data which is large, discontinuous,
and spherically or cylindrically symmetric in three space dimension.
Besides, when the solutions are not necessarily symmetric, we shall
establish a blow-up criterion for strong solutions with vacuum.

 For compressible isentropic Navier-Stokes equations (i.e. no
temperature equation), there are so many results about the
well-posedness and asymptotic behaviors of the solutions when
vacuum appears. Refer to \cite{Feireisl2, Jiang-Zhang1, Lions2,
Luo} and \cite{Guo-Zhu, Liu-Xin-Yang, Vong-Yang-Zhu, Yang-Yao-Zhu,
Yang-Zhu, Zhu} for global weak solutions with constant viscosity
and with density-dependent viscosity, respectively. Refer to
\cite{Choe-Kim: Radially, Fan-Jiang-Ni} and \cite{Cho-Choe-Kim,
Cho-Kim, Choe-Kim: Strong, salvi} for global strong solutions and
for local strong (classical) solutions with constant viscosity,
respectively. Recently, Huang, Li, Xin in \cite{Huang-Li-Xin} and
Ding, Wen, Yao, Zhu in \cite{Ding-Wen-Zhu, Ding-Wen-Yao-Zhu}
independently got existence and uniqueness of global classical
solutions, where the initial energy in \cite{Huang-Li-Xin} is
assumed to be small in $\mathbb{R}^3$ and
$\rho-\widetilde{\rho}\in C\left([0,T]; H^3(\mathbb{R}^3)\right)$,
$u\in C\left([0,T];D^1(\mathbb{R}^3)\cap
D^3(\mathbb{R}^3)\right)\cap L^\infty\left([\tau, T];
D^4(\mathbb{R}^3)\right)$ (for $\tau>0$) which generalized the
results in \cite{Cho-Kim}, and the initial data in
\cite{Ding-Wen-Zhu, Ding-Wen-Yao-Zhu} could be large for dimension
$N=1$ and could be large but spherically symmetric for $N\ge2$,
and $(\rho,u)\in C([0,T]; H^4(I))$ ($I$ is bounded in
\cite{Ding-Wen-Zhu}, and is bounded or an exterior domain in
\cite{Ding-Wen-Yao-Zhu}).\\[1.5mm]

We would like to give some notations which will be used throughout
the paper.\\

 \noindent{\bf Notations:}\\

(i) $I=[a,b]$; $Q_T=I\times[0,T]$ for $T>0$.\\

(ii) $\displaystyle\int_{\sum} f =\int_{\sum} f \,dx$, for $\sum=I$ or $\mathbb{R}^3$.\\

(iii)\ For $1\le l\le \infty$, denote the $L^l$ spaces and the
standard Sobolev spaces as follows:
$$L^l=L^l(\Sigma),  \ D^{k,l}=\left\{ u\in L^1_{\rm{loc}}(\Sigma): \|\nabla^k u \|_{L^l}<\infty\right\},$$
$$W^{k,l}=L^l\cap D^{k,l},  \ H^k=W^{k,2}, \ D^k=D^{k,2},$$
$$D_0^1=\Big\{u\in L^6: \ \|\nabla u\|_{L^2}<\infty\},$$
$$\|u\|_{D^{k,l}}=\|\nabla^k u\|_{L^l}.$$

(iv)\ For two $3\times 3$ matrices $E=(E_{ij}), F=(F_{ij})$, denote
the scalar product between $E$ and $F$ by
$$E:F=\sum_{i,j=1}^3 E_{ij} F_{ij}.$$

(v)\ $G=(2\mu+\lambda)\mathrm{div}u-P$ is the effective viscous flux.\\

(vi)\ $\dot{h}=h_t+u\cdot\nabla h$ denotes the material
derivative.\\

The rest of the paper is organized as follows. In Section 2, we
present three main theorems of the paper. In Section 3, some useful
lemmas are stated, which will be used to prove the global existence
of classical and strong solutions. The main theorems will be proved
in Sections 4, 5 and 6, respectively.
\section{ Main results}

\setcounter{equation}{0} \setcounter{theorem}{0}
\subsection{Global symmetric classical and strong solutions with
vacuum}\label{sec 2.1} \setcounter{equation}{0}
\setcounter{theorem}{0}
\renewcommand{\theequation}{\thesubsection.\arabic{equation}}
Throughout Section \ref{sec 2.1} and the proofs of the main theorems
in the section, we take $\Sigma=I$ in the {\bf Notations}. For
simplicity, we assume that the external force ${\bf f}=0$. Assume
$\Omega=\{{\bf x}|a<|{\bf x}|<b\}$, for $0<a<b<\infty$. Then for the
symmetric cases, (\ref{full N-S}) takes the form: \beq
\label{non-1.2}
\begin{cases}\displaystyle
\rho_t+(\rho u)_x+\frac{m\rho u}{x}=0,\ \rho\ge0,\ a<x<b,\
t>0,\\[2mm]
\displaystyle \rho u_t+\rho uu_x-\frac{\rho
v^2}{x}+P_x=\beta(u_{xx}+\frac{mu_x}{x}-\frac{mu}{x^2}),\\[2mm]
\displaystyle\rho v_t+\rho uv_x+\frac{\rho
uv}{x}=\mu(v_{xx}+\frac{mv_x}{x}-\frac{mv}{x^2}),\\ \displaystyle
\rho w_t+\rho uw_x=\mu(w_{xx}+\frac{mw_x}{x}),\\ \displaystyle \rho
e_t+\rho ue_x+P(u_x+\frac{mu}{x})=(\kappa
\theta_x)_x+\frac{m\kappa\theta_x}{x}+\wp,
\end{cases} \eeq where $\beta=2\mu+\lambda$, $\wp=\lambda(u_x+\frac{mu}{x})^2+\mu\left(w_x^2+2u_x^2+(v_x-\frac{mv}{x})^2+\frac{2mu^2}{x^2}\right)$.
In the spherically symmetric case, $m = N-1,$ $x = |{\bf x}|$, ${\bf
u}({\bf x}, t)=u(x,t)\frac{{\bf x}}{x}$ and $v=w=0$. In the
cylindrically symmetric case, $m = 1,$ $x =\sqrt{x_1^2+x_2^2}$, and
$${\bf u}({\bf x}, t)=u(x, t)\frac{(x_1, x_2, 0)}{x}+v(x, t)
\frac{(-x_2, x_1, 0)}{x}+w(x, t)(0, 0, 1).$$
 We consider the initial and boundary
conditions: \be\label{non-1.3}(\rho, \ u,\ v,\ w,\
\theta)\big|_{t=0}=(\rho_0,\ u_0,\ v_0,\ w_0,\ \theta_0)(x) \ \
{\rm{in}} \ \ I, \ee and\be\label{non-1.4}
 (u,v,w,\ \theta_x)\big|_{x=a,b}=0, \ t\geq0.
\ee

 \subsubsection{Assumptions}

($A_1$): $\rho_0\ge0$, $\int_I\rho_0>0$.\\

\noindent($A_2$): $\mu$ and $\lambda$ are constants. $\mu>0$,
$2\mu+(m+1)\lambda>0$. $e=C_0Q(\theta)+e_c(\rho)$, $P=\rho
Q(\theta)+P_c(\rho)$, $\kappa=\kappa(\theta)$, for some constant
$C_0>0$. The constant $C_0$ plays no role in the analysis, we assume
$C_0=1$.\\

\noindent($A_3$): $P_c(\rho)\ge0$, $e_c(\rho)\ge0$, for $\rho\ge0$;
$P_c\in C^2[0,\infty)$; $\rho| \frac{\partial e_c}{\partial\rho}|\le
C_1 e_c(\rho)$, for some constant $C_1>0$.\\

\noindent$(A_4):$ $Q(\cdot)\in C^2[0,\infty)$ satisfies
$$\begin{cases}C_2\left(\beta+(1-\beta)\theta+\theta^{1+r}\right)\le
Q(\theta)\le
C_3\left(\beta+(1-\beta)\theta+\theta^{1+r}\right),\\
C_4(1+\theta^r)\le Q'(\theta)\le C_5(1+\theta^r) ,\end{cases}\ \ \
\ \ \  \ \ \ \ \ \ \
$$
for some constants $C_i>0$ ($i= 2$, $3$, $4$, $5$) and $r\ge 0$,
$\beta=0$
or $1$.\\

\noindent($A_5$): $\kappa\in C^2[0,\infty)$ satisfies
 $$C_6(1+\theta^q)\le \kappa(\theta)\le
C_7(1+\theta^q),$$ for
$q>r$, and some constants $C_i>0$ ($i=6$, $7$).\\

\subsubsection{Global strong solutions}

\begin{theorem}\label{th:1.1}( Strong solutions)
In addition to $(A_1)$-$(A_5)$, we assume $\rho_0\geq0$,
 $\rho_0\in H^2$, $u_0\in H^2\cap H_0^1$,
$\theta_0\in H^2$, $\partial_x\theta_0|_{x=0,1}=0$, and that the
following compatibility conditions are valid:
 \beq\label{compatibility conditions}
\begin{cases}\displaystyle
\beta(u_{0xx}+\frac{mu_{0x}}{x}-\frac{mu_0}{x^2})-P_x(\rho_0,\theta_0)=\sqrt{\rho}_0g_1,\\[3mm]
\displaystyle
\mu(v_{0xx}+\frac{mv_{0x}}{x}-\frac{mv_0}{x^2})=\sqrt{\rho}_0g_2,\\[3mm]
\displaystyle \mu(w_{0xx}+\frac{mw_{0x}}{x})=\sqrt{\rho}_0g_3,
\\[3mm] \displaystyle
 \left(\kappa(\theta_0)
\theta_{0x}\right)_x+\frac{m\kappa(\theta_0)\theta_{0x}}{x}+\wp(x,0)=\sqrt{\rho_0}\
g_4,\ x\in I,
\end{cases}\eeq for some $g_i\in L^2$, $i=1,2,3,4$. Then there exists
a unique global solution $(\rho,\ u,\ v,\ w,\ \theta)$ to
(\ref{non-1.2})-(\ref{non-1.4}) such that for any $T>0$  \bex
&\rho\in C([0,T];H^2),\ (u, v, w, \theta)\in C([0,T];H^2)\cap
L^2([0,T];H^3),&\\& (\sqrt{\rho}u_t,\sqrt{\rho}v_t,\sqrt{\rho}w_t,
\sqrt{\rho}e_t)\in L^\infty([0,T]; L^2),\ (u_t, v_t, w_t)\in
L^2([0,T]; H_0^1),\ \theta_t\in L^2([0,T]; H^1).& \eex
\end{theorem}
\begin{remark}
 From the assumptions ($A_2$)-($A_4$), we know that the polytropic
perfect gas (i.e., $P=R\rho\theta$, $e=C_\nu \theta$ for some
constants $R>0$ and $C_\nu>0$) is included if we take $r=\beta=0$
and $e_c(\rho)=P_c(\rho)\equiv0$ and $Q=\frac{C_\nu\theta}{C_0}$.
\end{remark}

\begin{remark} The global existence of strong solutions depends on $q>r$ in
our analysis. For the polytropic perfect gas, Theorem \ref{th:1.1}
works for any $q>0$, which relaxes the restriction $q\ge 2$ in
\cite{Feireisl-book,Kawohl, Wen-Zhu}.
\end{remark}

\begin{remark} Some similar compatibility conditions as (\ref{compatibility
conditions}) can be referred to \cite{cho-Kim: perfect gas} and
references therein. In \cite{cho-Kim: perfect gas}, the
$\mathrm{local}$ $H^2$-{\rm regularity} of $u$ and $\theta$ for
$\mathrm{the\ polytropic\ perfect\ gas}$ was obtained. The detailed
reasons why such conditions were needed can be found in
\cite{cho-Kim: perfect gas}. Roughly speaking, $g_1$, $g_2$, $g_3$
and $g_4$ are equivalent to $\sqrt{\rho}u_t$, $\sqrt{\rho}v_t$,
$\sqrt{\rho}w_t$ and $\sqrt{\rho}e_t$  at $t=0$,
respectively.
\end{remark}

\begin{remark} From the derivation of the Navier-Stokes equations from the
Boltzmann equation through the Chapman-Enskog expansion to the
second order (see \cite{Yang} and references therein), we know that
$\mu=\mu(\theta)$, $\lambda=\lambda(\theta)$ and
$\kappa=\kappa(\theta)$. As in \cite{Feireisl-book,Kawohl,Wen-Zhu},
$\mu$ and $\lambda$ are assumed to be constants here, because of the
restrictions of mathematical technique.
\end{remark}

\subsubsection{Global classical solutions}
\begin{theorem}\label{th:1.2}(Classical solutions)
In addition to $(A_1)$-$(A_5)$, we assume $\rho_0\geq0$,
 $\rho_0\in H^3$, $(\sqrt{\rho_0})_x\in L^\infty$, $u_0\in H^3\cap H_0^1$,
$\theta_0\in H^3$, $\partial_x\theta_0|_{x=0,1}=0$, $(Q,
P_c,\kappa)\in W^{3,\infty}$, and that the compatibility
conditions (\ref{compatibility conditions}) are satisfied for some
$g_i\in L^2$, $i=1,2,3,4$, and $\sqrt{\rho_0}g_j\in H_0^1$,
$j=1,2,3$, and $\sqrt{\rho_0}g_4\in H^1$. Then there exists a
unique global solution $(\rho,\ u,\ v,\ w,\ \theta)$ to
(\ref{non-1.2})-(\ref{non-1.4}) such that for any $T>0$  \bex
&\rho\in C([0,T];H^3),\ (u, v, w, \theta)\in
L^\infty([0,T];H^3)\cap L^2([0,T];H^4),&\\&
(\sqrt{\rho}u_t,\sqrt{\rho}v_t,\sqrt{\rho}w_t, \sqrt{\rho}e_t)\in
L^\infty([0,T]; L^2),\ (\rho u_t,\rho v_t,\rho w_t)\in
L^\infty([0,T]; H_0^1),&\\& (u_t, v_t, w_t)\in L^2([0,T]; H_0^1),\
\theta_t\in L^2([0,T]; H^1),\ , \rho e_t\in L^\infty([0,T]; H^1).&
\eex
\end{theorem}
\begin{remark}
  If $\rho_0\in H^4$, and ($u_0$, $v_0$, $w_0$) satisfies a
stronger compatibility conditions
 \beq\label{compatibility conditions1}
\begin{cases}\displaystyle\beta(u_{0xx}+\frac{mu_{0x}}{x}-\frac{mu_0}{x^2})-P_x(\rho_0,\theta_0)=\rho_0
\widetilde{g}_1,\\[3mm] \displaystyle
\mu(v_{0xx}+\frac{mv_{0x}}{x}-\frac{mv_0}{x^2})=\rho_0\widetilde{g}_2,\\[3mm]
\displaystyle
\mu(w_{0xx}+\frac{mw_{0x}}{x})=\rho_0\widetilde{g}_3,\ x\in I,
\end{cases}\eeq
for some $\widetilde{g}_i\in H_0^1$ and
$(\sqrt{\rho_0}\partial_x\widetilde{g}_{i})_x\in L^2$, we can obtain
by using the similar arguments as in \cite{Wen-Zhu} that $\rho\in
C([0,T];H^4)$ and $(u, v, w)\in C([0,T];H^4)\cap L^2([0,T];H^5)$.
\end{remark}

\subsection{A blow-up criterion in terms of
$\|\rho\|_{L^\infty_tL_x^\infty}$ and
$\|\rho\theta\|_{L^4_tL^{\frac{12}{5}}_x}$ for strong
solutions}\label{sec 2.2}
 \setcounter{equation}{0}
\setcounter{theorem}{0}
\renewcommand{\theequation}{\thesubsection.\arabic{equation}}

Throughout Section \ref{sec 2.2}  and the proofs of the main
theorems in the section, we take $\Sigma=\mathbb{R}^3$ in the {\bf
Notations}. In order to establish some sharp blow-up criterions, we
only consider that $\kappa=$constant, and that the state equations
of $P$ and $e$ is of ideal polytropic gas type: $ P=a\rho\theta,\
e=C_0\theta$, where $a$ and $C_0$ are two positive constants. The
constants $a$, $C_0$ and $\kappa$ in the equations play no roles in
the section, we assume $a=C_0=\kappa=1$. If the solutions are
regular enough (such as strong solutions), (\ref{full N-S}) is
equivalence to the following system:
 \be\label{full N-S+1}
\begin{cases}
\rho_t+\nabla \cdot (\rho u)=0, \\
\rho  u_t+\rho u\cdot\nabla u+\nabla
P(\rho,\theta)=\mu\Delta u+(\mu+\lambda)\nabla\mathrm{div}u,\\
\rho \theta_t+\rho
u\cdot\nabla\theta+\rho\theta\mathrm{div}u=\frac{\mu}{2}\left|\nabla
u+(\nabla u)^\prime\right|^2+\lambda(\mathrm{div}u)^2+\Delta\theta,\
\mathrm{in}\ \mathbb{R}^3.
\end{cases}
\ee System (\ref{full N-S+1}) is supplemented with initial
conditions \be\label{non-initial} (\rho, u, \theta)|_{t=0}=(\rho_0,
u_0, \theta_0),\ x\in\mathbb{R}^3, \ee with \be\label{non-boundary}
\rho(x,t)\rightarrow0,\ u(x,t)\rightarrow0,\
\theta(x,t)\rightarrow0,\ \mathrm{as}\ |x|\rightarrow\infty,\
\mathrm{for}\ t\ge0. \ee We give the definition of strong solutions
to (\ref{full N-S+1})-(\ref{non-boundary}) throughout Section
\ref{sec 2.2} and the proofs of main theorem in the section.
\begin{definition}(Strong solution) For $T>0$, $(\rho, u, \theta)$ is called a strong solution to the compressible Navier-Stokes equations
(\ref{full N-S+1})-(\ref{non-boundary}) in $\mathbb{R}^3\times [0,T]$, if for some $q\in (3, 6]$,
\bex &0\le \rho\in
C([0,T];W^{1,q}\cap H^1\cap L^1),\ \rho_t\in C([0,T];L^2\cap L^q),&\\
& (u, \theta)\in C([0,T];D^2\cap D^1_0)\cap L^2(0,T;D^{2,q}),\ (u_t,
\theta_t)\in L^2(0,T;D^1_0),&\\&
 (\sqrt{\rho} u_t, \sqrt{\rho} \theta_t)\in
L^\infty(0,T;L^2),& \eex and $(\rho,u,\theta)$ satisfies (\ref{full
N-S+1}) a.e. in $\mathbb{R}^3\times (0,T]$.
\end{definition}
We present our main theorem, which is on a blow-up criterion for
strong solutions to (\ref{full N-S+1})-(\ref{non-boundary}), as
follows:
\begin{theorem}\label{blowup-th:1.1}
Assume $\rho_0\geq0$,
 $\rho_0\in H^1\cap W^{1,q}\cap L^1$, for some $q\in(3,6]$, $(u_0, \theta_0)\in D^2\cap D_0^1$, and the
following compatibility conditions are satisfied:
\beq\label{non-compatibility}
\begin{cases}
\mu\Delta u_0+(\mu+\lambda)\nabla\mathrm{div} u_0-\nabla P(\rho_0, \theta_0)=\sqrt{\rho}_0g_1,\\
\kappa\Delta\theta_0+\frac{\mu}{2}\left|\nabla u_0+(\nabla
u_0)^\prime\right|^2+\lambda(\mathrm{div}u_0)^2=\sqrt{\rho_0}g_2,\
x\in\mathbb{R}^3,
\end{cases}
\eeq for some $g_i\in L^2$, $i=1,2$. Let $(\rho, u, \theta)$ be a
strong solution to (\ref{full N-S+1})-(\ref{non-boundary}) in
$\mathbb{R}^3\times[0, T]$. If $0<T^*<+\infty$ is the maximum time
of existence of the strong solution, then \be\label{non-result}
\lim\sup\limits_{T\nearrow T^*}\left(\|\rho\|_{L^\infty(0,T;
L^\infty)}+\|\rho\theta\|_{L^4(0,T; L^\frac{12}{5})}\right)=\infty,
\ee provided $3\mu>\lambda$.
\end{theorem}
\begin{remark}  Under the conditions of Theorem \ref{blowup-th:1.1},
the local existence of the strong solutions was obtained in
\cite{cho-Kim: perfect gas}. Thus, the assumption $T^*>0$ makes
sense.
\end{remark}

\begin{remark} For the ideal polytropic gas with $\kappa=$constant, we
noticed recently that Huang and Li in \cite{Huang-Li} got the global
existence of classical and weak solutions to Cauchy problem of
(\ref{full N-S}) in $\mathbb{R}^3$ with small initial energy and
non-vacuum state at infinity.
\end{remark}

\begin{remark}
Before Theorem \ref{blowup-th:1.1}, there have been several results
on the blow-up criterions for strong solutions to (\ref{full
N-S+1}),
 please refer for instance to \cite{Fan-Jiang-Ou: Blow-up criterion, Fang-Zi-Zhang, Sun-Wang-Zhang 1, Wen-Zhu 3} and references therein. More precisely,

$\bullet$ Fan-Jiang-Ou (\cite{Fan-Jiang-Ou: Blow-up criterion}, 3D)
\be\label{Fan-Jiang-Ou criterion} \lim\sup\limits_{t\nearrow
T^\star}\left(\|\theta\|_{L^\infty(0,t;L^\infty)}+\|\nabla
u\|_{L^1(0,t;L^\infty)}\right)=\infty, \ee provided $7\mu>\lambda$.
Here the appearance of vacuum is allowed.

It is well-known that the bound of $\|\nabla
u\|_{L^1(0,t;L^\infty)}$ yields that
$\|\rho\|_{L^\infty(0,t;L^\infty)}$ is bounded (see (2.2) in
\cite{Fan-Jiang-Ou: Blow-up criterion}), if the initial density is
bounded. When $\|\nabla u\|_{L^1(0,t;L^\infty)}$ in
(\ref{Fan-Jiang-Ou criterion}) is relaxed by the upper bound of the
density, the following blow-up criterions were obtained:

$\bullet$ Fang-Zi-Zhang (\cite{Fang-Zi-Zhang}, 2D)
\be\label{Fang-Zi-Zhang criterion} \lim\sup\limits_{t\nearrow
T^\star}\left(\|\theta\|_{L^\infty(0,t;L^\infty)}+\|\rho\|_{L^\infty(0,t;L^\infty)}\right)=\infty,
\ee where the appearance of vacuum is allowed;

$\bullet$ Sun-Wang-Zhang (\cite{Sun-Wang-Zhang 1}, 3D)
\be\label{Sun-Wang-Zhang criterion} \lim\sup\limits_{t\nearrow
T^\star}\left(\|\theta\|_{L^\infty(0,t;L^\infty)}+\|\rho\|_{L^\infty(0,t;L^\infty)}+\left\|\frac{1}{\rho}\right\|_{L^\infty(0,t;L^\infty)}\right)=\infty,
\ee provided $7\mu>\lambda$. Here  the appearance of vacuum is not
allowed.

$\bullet$ Wen-Zhu (\cite{Wen-Zhu 3}, 3D) \be\label{Wen-Zhu
criterion} \lim\sup\limits_{t\nearrow
T^\star}\left(\|\theta\|_{L^\infty(0,t;L^\infty)}+\|\rho\|_{L^\infty(0,t;L^\infty)}\right)=\infty,
\ee provided $3\mu>\lambda$. Here the appearance of vacuum is
allowed.
\end{remark}

\begin{remark}
Theorem \ref{blowup-th:1.1} is an extension of our former results in
\cite{Wen-Zhu 3} (see (\ref{Wen-Zhu criterion})). One of the main
ingredients is that the estimates of
$\|\sqrt{\rho}\theta\|_{L^\infty(0,T;L^2)}$ and $\|\nabla
u\|_{L^\infty(0,T;L^2)}$ are done together, i.e.,
\bex\sup\limits_{0\le t\le
T}\int_{\mathbb{R}^3}(\rho|\theta|^2+|\nabla u|^2)\,
dx+\int_0^T\int_{\mathbb{R}^3}(|\nabla\theta|^2+\rho
|u_t|^2)\,dxdt\le C.\eex In \cite{Wen-Zhu 3},
$\|\sqrt{\rho}\theta\|_{L^\infty(0,T;L^2)}$ and $\|\sqrt{\rho}
u\|_{L^\infty(0,T;L^2)}$ were done together, which needed the upper
bounds of the temperature and the density.
\end{remark}

\setcounter{section}{2} \setcounter{equation}{0}
\section{Preliminaries}
\setcounter{equation}{0} \setcounter{theorem}{0}
\renewcommand{\theequation}{\thesection.\arabic{equation}}
\renewcommand{\thetheorem}{\thesection.\arabic{theorem}}
The lemmas in the section will be useful in the next two sections.

\begin{lemma}(\cite{Wen-Zhu})\label{non-le:2.1}
Let $\Omega=[\widetilde{a}, \widetilde{b}]$ be a bounded domain in
$\mathbb{R}$, and $\rho$ be a non-negative function such that
$$
0<M\le\int_\Omega\rho\le K,
$$ for constants $M>0$ and $K>0$. Then
$$ \|v\|_{L^\infty(\Omega)}\le \frac{K}{M}\|v_x\|_{L^1(\Omega)}+\frac{1}{M}\left|\int_\Omega\rho v\right|,
$$
 for any $v\in H^1(\Omega)$.
\end{lemma}
\begin{remark}
The version of higher dimensions for Lemma \ref{non-le:2.1} can be
found in \cite{Feireisl1} or \cite{Feireisl-book}.
\end{remark}

\begin{corollary}(\cite{Wen-Zhu})\label{non-cor:2.1}
Consider the same conditions in Lemma \ref{non-le:2.1}, and in
addition assume $\Omega=I$, and
$$
\|\rho v\|_{L^1(I)}\le \overline{c}.
$$
Then for any $l>0$, there exists a positive constant $C=C(M, K, l,
\overline{c})$ such that $$\|v^l\|_{L^\infty(I)}\le
C\|(v^l)_x\|_{L^1(I)}+C,$$ for any $v^l\in H^1(I)$.
\end{corollary}

\begin{lemma}(Poincar$\acute{e}$ inequality)\label{non-le:2.2}
For any $v\in H^1_0(I)$, we have
$$
\|v\|_{L^\infty(I)}\le \|v_x\|_{L^1}.
$$
\end{lemma}

\begin{lemma}\label{non-le:2.3} (\cite{Simon}).
Assume $X\subset E\subset Y$ are Banach spaces and
$X\hookrightarrow\hookrightarrow E$. Then the following imbedding
are compact: $$(i)\ \  \left\{\varphi:\varphi\in L^q(0,T; X),
\frac{\partial\varphi}{\partial t}\in
L^1(0,T;Y)\right\}\hookrightarrow\hookrightarrow L^q(0,T; E),\ \
{\rm if}\ \  1\leq q\leq\infty; $$ $$(ii)\ \
\left\{\varphi:\varphi\in L^\infty(0,T; X),
\frac{\partial\varphi}{\partial t}\in
L^r(0,T;Y)\right\}\hookrightarrow\hookrightarrow C([0,T]; E),\ \
{\rm if}\ \   1< r\leq\infty. $$
\end{lemma}

\section{Proof of Theorem \ref{th:1.1}}
\setcounter{equation}{0} \setcounter{theorem}{0} In the section, we
denote by $C$ a generic constant depending only on $\|(\rho_0, u_0,
v_0, w_0, \theta_0)\|_{H^2}$, $\|g_i\|_{L^2}$ ($i=1,2,3,4$), $T$,
$\lambda$, $\mu$, $a$, $b$, and some other known constants, but
independent of the solutions and the lower bounds of the density. We
denote by
$$A\lesssim B$$ if there exists a generic constant $C$ such that $A\leq C
B$.

The strategies on proving Theorem \ref{th:1.1} are very classical.
More precisely, we derive various {\it a priori} estimates for
strong solutions of the Navier-Stokes equations
(\ref{non-1.2})-(\ref{non-1.4}), which are independent of positive
lower bounds of the initial density. Then we shall construct a
sequence of approximate initial data where the initial density has a
lower bound $\varepsilon>0$. With these {\it a priori} estimates
uniform for $\varepsilon$, we take the limits
$\varepsilon\rightarrow0^+$.

From now on, for any $T>0$, we shall derive some delicate {\it a
priori} estimates for the strong solutions $(\rho, u, v,
w,\theta)$ as in Theorem \ref{th:1.1} with $\inf\limits_{(x,t)\in
Q_T}\rho>0$. These energy estimates will be finished by five
steps.
\\

{\noindent\bf Step 1: Basic energy inequality}

\begin{lemma}\label{le:energy identity} Under the conditions of
Theorem \ref{th:1.1}, we have for any $t\in[0,T]$ \bex
\int_Ix^m\rho(1+e+u^2+v^2+w^2)\le C. \eex
\end{lemma}
\pf This bound is standard and follow directly from the equations
(\ref{full N-S})$_1$, (\ref{full N-S})$_3$ and the boundary
conditions.
\endpf\\

{\noindent\bf Step 2: Upper bound of density}
\begin{lemma}\label{le:upper bound of density}Under the conditions of
Theorem \ref{th:1.1}, we have
 \bex \|\rho\|_{L^\infty(Q_T)}\le C.
\eex
\end{lemma}
\pf The idea of the proof is essentially that of a similar result
of Frid, Shelukhin \cite{Frid} where $m=1$, but with a slightly
modification. We omit it for brevity.
\endpf\\

{\noindent\bf Step 3: $H^1$-estimates of $(\rho,u,v,w)$}

The next lemma plays an important role in the paper, whose proofs
are improved in contrast with \cite{Wen-Zhu}. The condition $q>r$
instead of $q\ge2+2r$ is enough here.
\begin{lemma}\label{le: theta}
Under the conditions of Theorem \ref{th:1.1},
 for $q>r$, and for any $0<\alpha<\min\{1, q-r\}$,
assume $2\mu+(m+1)\lambda>0$, we have \bex
\int_{Q_T}x^m\frac{(1+\theta^q)\theta_x^2}{\theta^{1+\alpha}}\le
C, \eex where the generic constant C depends on $\alpha$.
\end{lemma}
\begin{remark}
The proofs of this lemma depend on the boundary condition
${\theta_x}|_{x=a,b}=0$.
\end{remark}
\pf From (\ref{non-1.1}) and (\ref{non-1.2}), we get
\be\label{non-3.5} \rho e_\theta\theta_t+\rho u e_\theta
\theta_x+\theta P_\theta
(u_x+\frac{mu}{x})=\left(\kappa(\theta)\theta_x\right)_x+\frac{m\kappa\theta_x}{x}+\wp.
\ee Substituting $e=Q(\theta)+e_c(\rho)$ and $P=\rho
Q(\theta)+P_c(\rho)$ into (\ref{non-3.5}), we get
\be\label{non-3.6} \rho Q'(\theta)\theta_t+\rho u
Q'(\theta)\theta_x+\rho\theta
Q'(\theta)(u_x+\frac{mu}{x})=\left(\kappa(\theta)\theta_x\right)_x+\frac{m\kappa\theta_x}{x}+\wp,
\ee or \be\label{non-3.7} (\rho Q)_t+(\rho u Q)_x+\frac{m\rho
uQ}{x}+\rho\theta
Q'(\theta)(u_x+\frac{mu}{x})=\left(\kappa(\theta)\theta_x\right)_x+\frac{m\kappa\theta_x}{x}+\wp.
\ee Multiplying (\ref{non-3.6}) by $x^m\theta^{-\alpha}$, and
integrating by parts over $Q_T$, we have
\beq\label{non-3.8}\begin{split}&\int_{Q_T}x^m\left(\frac{\alpha\kappa(\theta)\theta_x^2}{\theta^{1+\alpha}}+\frac{\lambda(u_x+\frac{mu}{x})^2}{\theta^\alpha}+
\frac{\mu[w_x^2+2u_x^2+(v_x-\frac{mv}{x})^2+\frac{2mu^2}{x^2}]}{\theta^\alpha}\right)\\=&
\int_Ix^m\rho\int_0^\theta\frac{Q'(\xi)}{\xi^\alpha}-\int_Ix^m\rho_0\int_0^{\theta_0}\frac{Q'(\xi)}{\xi^\alpha}+\int_{Q_T}x^m\rho\theta^{1-\alpha}Q'(\theta)
(u_x+\frac{mu}{x})\\
\les&
\int_Ix^m\rho\left|\int_0^\theta\frac{1+\xi^r}{\xi^\alpha}\right|+\int_Ix^m\rho_0\left|\int_0^{\theta_0}\frac{1+\xi^r}{\xi^\alpha}\right|
+\int_{Q_T}x^m\rho\theta^{1-\alpha}(1+\theta^r)|u_x+\frac{mu}{x}|\\
\les&
\int_Ix^m\rho(1+\theta^{1+r})+\int_Ix^m\rho_0(1+\theta_0^{1+r})+\int_{Q_T}x^m\rho\theta^{1-\alpha}(1+\theta^r)|u_x+\frac{mu}{x}|,\end{split}
\eeq where we have used ($A_4$) and Young inequality.  Since
$\mu>0$, we have from (\ref{non-3.8}), ($A_2$), ($A_4$) and Lemma
\ref{le:energy identity}
\beq\label{non-3.8+1}\begin{split}\alpha\int_{Q_T}x^m\frac{\kappa\theta_x^2}{\theta^{1+\alpha}}
\le&C-\int_{Q_T}x^m\frac{\lambda(u_x+\frac{mu}{x})^2+2\mu(u_x^2+\frac{mu^2}{x^2})}{\theta^\alpha}\\&+C\int_{Q_T}x^m\rho\theta^{1-\alpha}(1+\theta^r)|u_x+\frac{mu}{x}|=
I_1+I_2+I_3. \end{split}\eeq Without loss of generality, we assume
$\lambda<0$. In fact, if $\lambda\ge0$, $I_2$ is obviously a good
term: \bex\begin{split}
I_2\le-2\mu\int_{Q_T}x^m\frac{u_x^2+\frac{mu^2}{x^2}}{\theta^\alpha}.
\end{split}\eex
For $\lambda<0$, we use Cauchy inequality to get
 \beq\label{non-3.8+2}\begin{split}I_2
=&-\int_{Q_T}x^m\frac{(2\mu+\lambda)u_x^2+(2\mu+m\lambda)\frac{mu^2}{x^2}+\frac{2m\lambda
uu_x}{x}}{\theta^\alpha}\\
\le&-[2\mu+(m+1)\lambda]\int_{Q_T}x^m\frac{u_x^2+\frac{mu^2}{x^2}}{\theta^\alpha}.
\end{split}\eeq
For $I_3$, using Cauchy inequality again, we get
\beq\begin{split}\label{non-3.8+3} I_3\le&
[2\mu+(m+1)\lambda]\int_{Q_T}x^m\frac{u_x^2+\frac{mu^2}{x^2}}{\theta^\alpha}+C\int_{Q_T}x^m\rho^2(1+\theta^{2+2r-\alpha})\\
\le&[2\mu+(m+1)\lambda]\int_{Q_T}x^m\frac{u_x^2+\frac{mu^2}{x^2}}{\theta^\alpha}+C\int_0^T\|\theta\|_{L^\infty}^{1+r-\alpha}\int_Ix^m\rho\theta^{1+r}
+C\\
\le&[2\mu+(m+1)\lambda]\int_{Q_T}x^m\frac{u_x^2+\frac{mu^2}{x^2}}{\theta^\alpha}+C\int_0^T\|\theta\|_{L^\infty}^{1+r-\alpha}
+C.
\end{split}
\eeq Substituting (\ref{non-3.8+2}) and (\ref{non-3.8+3}) into
(\ref{non-3.8+1}), we have
\beq\label{non-3.8+4}\begin{split}\alpha\int_{Q_T}x^m\frac{\kappa\theta_x^2}{\theta^{1+\alpha}}
\le&C\int_0^T\|\theta\|_{L^\infty}^{1+r-\alpha}+C.
\end{split}\eeq
 Now we estimate
the first term of the right hand side of (\ref{non-3.8+4}).

\noindent{\bf Case 1: $r<q<1+2r-\alpha$} \bex\begin{split}
C\int_0^T\|\theta\|_{L^\infty}^{1+r-\alpha}\les&
1+\int_0^T\|\theta^{r-\alpha}\theta_x\|_{L^2}\\ \le&
C+C\int_0^T\left(\int_I\frac{\theta_x^2\theta^q}{\theta^{1+\alpha}}\theta^{2r-\alpha+1-q}\right)^\frac{1}{2}\\
\le&\frac{1}{4}\alpha\int_{Q_T}x^m\frac{\kappa\theta_x^2}{\theta^{1+\alpha}}+C\int_0^T\|\theta\|_{L^\infty}^{2r-\alpha+1-q}+C\\
\le&\frac{1}{4}\alpha\int_{Q_T}x^m\frac{\kappa\theta_x^2}{\theta^{1+\alpha}}+\frac{1}{2}C\int_0^T\|\theta\|_{L^\infty}^{1+r-\alpha}+C,\end{split}
\eex where we have used Corollary \ref{non-cor:2.1}, ($A_5$) and
Young inequality.

 This gives \beq\label{non-3.9}\begin{split}
C\int_0^T\|\theta\|_{L^\infty}^{1+r-\alpha}\le&\frac{1}{2}\alpha\int_{Q_T}x^m\frac{\kappa\theta_x^2}{\theta^{1+\alpha}}+C.
\end{split}
\eeq

\noindent{\bf Case 2: $q\ge1+2r-\alpha$}

Using Young inequality, we have \beq\label{non-3.10}\begin{split}
C\int_0^T\|\theta\|_{L^\infty}^{1+r-\alpha}\les&
1+\int_0^T\|\theta^{r-\alpha}\theta_x\|_{L^2}\\ \les&
1+\int_0^T\left(\int_I\frac{\theta_x^2\theta^{1+2r-\alpha}}{\theta^{1+\alpha}}\right)^\frac{1}{2}\\
\le&C+C\int_0^T\left(\int_I\frac{\kappa\theta_x^2}{\theta^{1+\alpha}}\right)^\frac{1}{2}\\
\le&
\frac{1}{2}\alpha\int_{Q_T}x^m\frac{\kappa\theta_x^2}{\theta^{1+\alpha}}+C.\end{split}
\eeq Substituting (\ref{non-3.9}) or (\ref{non-3.10}) into
(\ref{non-3.8+4}), we complete the proof of Lemma \ref{le: theta}.
\endpf

The next estimate is a corollary of Lemma \ref{le: theta}, whose
proof can be found in \cite{Wen-Zhu} (Corollary 3.1). For
completeness, we present the proof.
\begin{corollary}\label{cor:int0T theta 1} Under the conditions of
Theorem \ref{th:1.1}, we have \bex
\int_0^T\|\theta\|_{L^\infty}^{q-\alpha+1}\le C. \eex
\end{corollary}
\pf By Corollary \ref{non-cor:2.1}, ($A_5$) and Lemma \ref{le:
theta}, we have \bex\begin{split}
\int_0^T\|\theta\|_{L^\infty}^{q-\alpha+1}=&\int_0^T\|\theta^\frac{q-\alpha+1}{2}\|_{L^\infty}^2\\
\le&
C\int_0^T\int_I\left(\theta^{\frac{q-\alpha-1}{2}}\theta_x\right)^2+C\\=&C\int_0^T\int_I\frac{\theta^{q}\theta_x^2}{\theta^{\alpha+1}}+C\\
\le& C. \end{split}\eex
\endpf

\begin{lemma}\label{le:int0T H^1 of uvw}Under the conditions of
Theorem \ref{th:1.1}, we have \bex
\int_{Q_T}x^m(u_x^2+v_x^2+x^{-2}u^2+x^{-2}v^2)\le C. \eex
\end{lemma}
\pf Multiplying (\ref{non-1.2})$_2$ and (\ref{non-1.2})$_3$ by
$x^mu$ and $x^m v$ respectively, and integrating by parts over
$I$, we have \beq\begin{split}\label{non-3.11}
\frac{1}{2}\frac{d}{dt}\int_Ix^m\rho u^2-\int_Ix^{m-1}\rho
uv^2+\int_Ix^mu P_x+\beta\int_Ix^m(u_x^2+mx^{-2}u^2)=0.
\end{split}
\eeq

\beq\begin{split}\label{non-3.12}
\frac{1}{2}\frac{d}{dt}\int_Ix^m\rho v^2+\int_Ix^{m-1}\rho
uv^2+\mu\int_Ix^m(v_x^2+mx^{-2}v^2)=0.
\end{split}
\eeq Adding (\ref{non-3.12}) into (\ref{non-3.11}), we have
\beq\begin{split}\label{non-3.13}
\frac{1}{2}\frac{d}{dt}\int_Ix^m\rho (u^2+v^2)+\int_Ix^mu
P_x+\beta\int_Ix^m(u_x^2+mx^{-2}u^2)+\mu\int_Ix^m(v_x^2+mx^{-2}v^2)=0.
\end{split}
\eeq Integrating (\ref{non-3.13}) over $(0,T)$, and using
integration by parts and Cauchy inequality, we have
\bex\begin{split}
&\beta\int_{Q_T}x^m(u_x^2+mx^{-2}u^2)+\mu\int_{Q_T}x^m(v_x^2+mx^{-2}v^2)\\
\le&\frac{1}{2}\int_Ix^m\rho_0 (u_0^2+v_0^2)+\int_0^T\int_Ix^mu_x
P+m\int_0^T\int_Ix^{m-1}u P\\ \le&
C+\int_0^T\int_Ix^m(u_x+mx^{-1}u) P\\
\le&\frac{1}{2}\beta\int_{Q_T}x^m(u_x^2+mx^{-2}u^2)+C\int_{Q_T}x^m
P^2.
\end{split}
\eex This gives \bex\begin{split}
&\frac{1}{2}\beta\int_{Q_T}x^m(u_x^2+mx^{-2}u^2)+\mu\int_{Q_T}x^m(v_x^2+mx^{-2}v^2)\\
\les&\int_{Q_T}x^m \rho^2(1+\theta^{2+2r})+1\\ \les&
\int_0^T\|\theta\|_{L^\infty}^{1+r}\int_Ix^m\rho\theta^{1+r}+1\\
\les&\int_0^T\|\theta\|_{L^\infty}^{1+r}+1,
\end{split}
\eex where we have used ($A_2$), ($A_3$), ($A_4$) and Lemmas
\ref{le:energy identity}-\ref{le:upper bound of density}.

Since $q>r$ and $0<\alpha<q-r$, we have $1+r<q-\alpha+1$. Thus,
using Young inequality and Corollary \ref{cor:int0T theta 1}, we
have \bex\begin{split}
&\frac{1}{2}\beta\int_{Q_T}x^m(u_x^2+mx^{-2}u^2)+\mu\int_{Q_T}x^m(v_x^2+mx^{-2}v^2)\\
\les&\int_0^T\|\theta\|_{L^\infty}^{q-\alpha+1}+1\le C.
\end{split}
\eex

\endpf

\begin{lemma}\label{le:H^1 of v}Under the conditions of
Theorem \ref{th:1.1}, we have for any $t\in[0,T]$ \bex
\int_Ix^m(v_x^2+x^{-2}v^2)+\int_{Q_T}x^m\rho v_t^2\le C. \eex
\end{lemma}
\pf Multiplying (\ref{non-1.2})$_3$ by $x^mv_t$, integrating by
parts over $I$, and using Cauchy inequality, we have
\bex\begin{split}
&\frac{\mu}{2}\frac{d}{dt}\int_Ix^m(v_x^2+mx^{-2}v^2)+\int_Ix^m\rho
v_t^2\\=&-\int_Ix^m\rho u v_xv_t-\int_Ix^{m-1}\rho uvv_t\\
\le&\frac{1}{2}\int_Ix^m\rho v_t^2+C\int_Ix^m\rho
u^2v_x^2+C\int_Ix^{m-2}\rho u^2v^2.
\end{split}
\eex Thus, we have \bex\begin{split}
&\mu\frac{d}{dt}\int_Ix^m(v_x^2+mx^{-2}v^2)+\int_Ix^m\rho
v_t^2\\
\les&\|u\|_{L^\infty}^2\|\rho\|_{L^\infty}\int_Ix^mv_x^2+\|v\|_{L^\infty}^2\int_Ix^{m}\rho
u^2\\
\les&\int_Ix^m(u_x^2+mx^{-2}u^2)\int_Ix^mv_x^2+\int_Ix^m(v_x^2+mx^{-2}v^2),
\end{split}
\eex where we have used Lemmas \ref{le:energy
identity}-\ref{le:upper bound of density} and
Poincar$\mathrm{\acute{e}}$ inequality.

By Gronwall inequality and Lemma \ref{le:int0T H^1 of uvw}, we
complete the proof of Lemma \ref{le:H^1 of v}.
\endpf

\begin{corollary}\label{cor:tx H^2 of v} Under the conditions of
Theorem \ref{th:1.1}, we have  \bex
\|v\|_{L^\infty(Q_T)}+\int_{Q_T}x^mv_{xx}^2\le C. \eex
\end{corollary}
\pf This is an immediately result from Lemmas \ref{le:upper bound
of density}, \ref{le:int0T H^1 of uvw}, \ref{le:H^1 of v},
Poincar$\mathrm{\acute{e}}$ inequality and (\ref{non-1.2})$_3$.
\endpf

\begin{lemma}\label{le:H^1 of w}Under the conditions of Theorem \ref{th:1.1}, we
have for any $t\in[0,T]$ \bex \int_Ix^mw_x^2+\int_{Q_T}x^m\rho
w_t^2\le C. \eex
\end{lemma}
\pf Multiplying (\ref{non-1.2})$_4$ by $x^mw_t$, integrating by
parts over $I$, and using Cauchy inequality, we have
\bex\begin{split} \int_Ix^m\rho
w_t^2+\frac{\mu}{2}\frac{d}{dt}\int_Ix^mw_x^2=&-\int_Ix^m\rho
uw_xw_t\\ \le&\frac{1}{2}\int_Ix^m\rho
w_t^2+\frac{1}{2}\int_Ix^m\rho u^2w_x^2.
\end{split}
\eex Thus, we apply Lemma \ref{le:upper bound of density} and
Poincar$\mathrm{\acute{e}}$ inequality to get \bex\begin{split}
\int_Ix^m\rho w_t^2+\mu\frac{d}{dt}\int_Ix^mw_x^2
\le&\|\rho\|_{L^\infty}
\|u\|_{L^\infty}^2\int_Ix^mw_x^2\\
\les&\int_Ix^m(u_x^2+mx^{-2}u^2)\int_Ix^mw_x^2.
\end{split}
\eex Using Gronwall inequality and Lemma \ref{le:int0T H^1 of
uvw}, we complete the proof of Lemma \ref{le:H^1 of w}.
\endpf \\

Similar to Corollary \ref{cor:tx H^2 of v}, we get the next
corollary.
\begin{corollary}\label{cor:tx H^2 of w} Under the conditions of
Theorem \ref{th:1.1}, we have \bex \int_{Q_T}x^mw_{xx}^2\le C.
\eex
\end{corollary}
\begin{lemma}\label{le:H^1 of u} Under the conditions of
Theorem \ref{th:1.1}, we have for any $t\in[0,T]$ \bex
\int_Ix^m(\rho
\theta^{q+r+2}+u_x^2+x^{-2}u^2)+\int_{Q_T}x^m\left(\rho
u_t^2+(1+\theta^q)^2\theta_x^2\right)\le C. \eex
\end{lemma}
\pf Multiplying (\ref{non-1.2})$_2$ by $x^mu_t$, integrating by
parts over $I$, and using Cauchy inequality, we have
\bex\begin{split} &\int_Ix^m\rho
u_t^2+\frac{\beta}{2}\frac{d}{dt}\int_Ix^m(u_x^2+mx^{-2}u^2)\\=&-\int_Ix^m\rho
uu_xu_t+\int_Ix^{m-1}\rho v^2u_t-\int_Ix^mP_xu_t\\
\le&\frac{1}{4}\int_Ix^m\rho
u_t^2+C\int_Ix^m\rho(u^2u_x^2+x^{-2}v^4)+\int_Ix^mPu_{xt}+m\int_I
x^{m-1}Pu_t.
\end{split}
\eex This, along with Lemma \ref{le:upper bound of density} and
Poincar$\mathrm{\acute{e}}$ inequality, deduces \be\label{d
ux-1}\begin{split} &\frac{3}{4}\int_Ix^m\rho
u_t^2+\frac{\beta}{2}\frac{d}{dt}\int_Ix^m(u_x^2+mx^{-2}u^2)\\
\le&C\|\rho\|_{L^\infty}\|u\|_{L^\infty}^2\int_Ix^mu_x^2+C\int_Ix^{m-2}\rho
v^4+\frac{d}{dt}\int_Ix^mPu_x-\int_Ix^mP_tu_x+m\int_I
x^{m-1}Pu_t\\
\le&C\left(\int_Ix^m(u_x^2+mx^{-2}u^2)\right)^2+\frac{d}{dt}\int_Ix^mPu_x+C\int_Ix^{m-2}\rho
v^4-\int_Ix^mP_tu_x+m\int_I
x^{m-1}Pu_t\\=&\sum\limits_{i=1}^5II_i.
\end{split}
\ee For $II_3$, using Lemma \ref{le:energy identity} and Corollary
\ref{cor:tx H^2 of v}, we have \be\label{II3}\begin{split}
II_3\les&\|v\|_{L^\infty}^2\int_Ix^{m}\rho v^2\le C.
\end{split}
\ee For $II_4$, we have \bex
\begin{split}
II_4=&-\beta^{-1}\int_Ix^mP_t(\beta
u_x-P)-\beta^{-1}\int_Ix^mP_tP\\
=&-\beta^{-1}\int_Ix^m(\rho Q)_t(\beta
u_x-P)-\beta^{-1}\int_Ix^mP_c(\rho)_t(\beta
u_x-P)-\frac{1}{2\beta}\frac{d}{dt}\int_Ix^mP^2\\=&\sum\limits_{i=1}^3II_{4,i}.
\end{split}
\eex For $II_{4,1}$, using (\ref{non-3.7}) and
(\ref{non-1.2})$_2$, we have \bex\begin{split}
II_{4,1}=&-\beta^{-1}\int_I(\beta
u_x-P)\left((x^m\kappa\theta_x)_x+x^m\wp-(x^m\rho u
Q)_x-x^m\rho\theta Q'(\theta)(u_x+\frac{mu}{x})\right)\\
=&\beta^{-1}\int_Ix^m(\beta u_x-P)_x(\kappa\theta_x-\rho u
Q)+\beta^{-1}\int_Ix^m(\beta u_x-P)\rho\theta
Q'(\theta)(u_x+\frac{mu}{x})\\&-\beta^{-1}\int_Ix^m(\beta
u_x-P)\wp\\ =&\beta^{-1}\int_Ix^m(\rho u_t+\rho uu_x-\frac{\rho
v^2}{x}-\frac{m\beta u_x}{x}+\frac{m\beta
u}{x^2})(\kappa\theta_x-\rho u Q)\\&+\beta^{-1}\int_Ix^m(\beta
u_x-P)\rho\theta
Q'(\theta)(u_x+\frac{mu}{x})-\beta^{-1}\int_Ix^m(\beta u_x-P)\wp.
\end{split}
\eex Recalling
$\wp=\lambda(u_x+\frac{mu}{x})^2+\mu\left(w_x^2+2u_x^2+(v_x-\frac{mv}{x})^2+\frac{2mu^2}{x^2}\right)$,
we have \bex\begin{split} II_{4,1} \le&\frac{1}{8}\int_Ix^m\rho
u_t^2+C\int_Ix^m\rho
u^2Q^2+C\int_Ix^m(\kappa\theta_x)^2+C\int_Ix^mu^2u_x^2+C\int_Ix^m\rho
v^4+C\int_Ix^mu_x^2\\&+C\int_Ix^{m-4}u^2+C\sup\limits_{x\in
I}(1+\theta^{1+r})\int_Ix^m\left(u_x^2+\rho^2Q^2+\rho+x^{-2}u^2\right)\\&+C\|\beta
u_x-P\|_{L^\infty}\int_Ix^m(u_x^2+x^{-2}u^2+w_x^2+v_x^2+x^{-2}v^2)\\
\le&\frac{1}{8}\int_Ix^m\rho
u_t^2+C\int_Ix^m(u_x^2+x^{-2}u^2)\int_Ix^m\rho(1+\theta^{2+2r})+C\int_Ix^m(\kappa\theta_x)^2\\&+C\left(\int_Ix^m(u_x^2+x^{-2}u^2)\right)^2
+C\sup\limits_{x\in
I}(1+\theta^{q-\alpha+1})\int_Ix^m\big(u_x^2+\rho(1+\theta^{2+2r})+x^{-2}u^2\big)\\&+C\|\beta
u_x-P\|_{L^\infty}\int_Ix^m(u_x^2+x^{-2}u^2+w_x^2+v_x^2+x^{-2}v^2)+C\\
\le&\frac{1}{8}\int_Ix^m\rho u_t^2+C\sup\limits_{x\in
I}\theta^{1+r}\int_Ix^m(u_x^2+x^{-2}u^2)+C\int_Ix^m(\kappa\theta_x)^2\\&+C\left(\int_Ix^m(u_x^2+x^{-2}u^2)\right)^2
+C\sup\limits_{x\in
I}(1+\theta^{q-\alpha+1})\int_Ix^m\left(u_x^2+\rho(1+\theta^{q+r+2})+x^{-2}u^2\right)\\&+C\|\beta
u_x-P\|_{L^\infty}\int_Ix^m(u_x^2+x^{-2}u^2+w_x^2+v_x^2+x^{-2}v^2)+C,
\end{split}
\eex where we have used Young inequality,
Poincar$\mathrm{\acute{e}}$ inequality, Lemmas \ref{le:energy
identity}-\ref{le:upper bound of density}, ($A_2$), ($A_4$),
(\ref{II3}), $q>r$ and $\alpha<q-r$.

This, along with Sobolev inequality, Cauchy inequality, Lemmas
\ref{le:upper bound of density}, \ref{le:H^1 of v}, \ref{le:H^1 of
w} and (\ref{non-1.2})$_2$, deduces

\bex\begin{split} II_{4,1} \le&\frac{1}{8}\int_Ix^m\rho
u_t^2+C\int_Ix^m(\kappa\theta_x)^2+C\left(\int_Ix^m(u_x^2+x^{-2}u^2)\right)^2
\\&+C\sup\limits_{x\in
I}(1+\theta^{q-\alpha+1})\left[1+\int_Ix^m\left(u_x^2+\rho\theta^{q+r+2}+x^{-2}u^2\right)\right]\\&+C\left(\|(\beta
u_x-P)_x\|_{L^2}+\|\beta
u_x-P\|_{L^2}\right)\left[1+\int_Ix^m(u_x^2+x^{-2}u^2)\right]+C\\
\le&\frac{1}{8}\int_Ix^m\rho
u_t^2+C\int_Ix^m(\kappa\theta_x)^2+C\left(\int_Ix^m(u_x^2+x^{-2}u^2)\right)^2
\\&+C\sup\limits_{x\in
I}(1+\theta^{q-\alpha+1})\left[1+\int_Ix^m\left(u_x^2+\rho\theta^{q+r+2}+x^{-2}u^2\right)\right]\\&+C\left(\|\rho
u_t+\rho uu_x-\frac{\rho v^2}{x}-\frac{m\beta u_x}{x}+\frac{m\beta
u}{x^2}\|_{L^2}+\|\beta
u_x-P\|_{L^2}\right)\left[1+\int_Ix^m(u_x^2+x^{-2}u^2)\right]+C\\
\le&\frac{1}{4}\int_Ix^m\rho
u_t^2+C\int_Ix^m(\kappa\theta_x)^2+C\left(\int_Ix^m(u_x^2+x^{-2}u^2)\right)^2
\\&+C\sup\limits_{x\in
I}(1+\theta^{q-\alpha+1})\left[1+\int_Ix^m\left(u_x^2+\rho\theta^{q+r+2}+x^{-2}u^2\right)\right]+C.
\end{split}
\eex For $II_{4,2}$, using (\ref{non-1.2})$_1$, integration by
parts, Lemma \ref{le:upper bound of density}, ($A_3$) and ($A_4$),
we have \bex\begin{split}
II_{4,2}=&\beta^{-1}\int_Ix^mP_c^\prime(\rho)(\rho_xu+\rho
u_x+mx^{-1}\rho u)(\beta
u_x-P)\\=&\beta^{-1}\int_Ix^m(P_c)_xu(\beta
u_x-P)+\beta^{-1}\int_Ix^m(P_c)^\prime \rho u_x(\beta
u_x-P)+\beta^{-1}\int_Imx^{m-1}(P_c)^\prime \rho u(\beta u_x-P)\\
\le&-\beta^{-1}\int_Ix^mP_cu_x(\beta
u_x-P)-\beta^{-1}\int_Ix^mP_cu(\beta
u_x-P)_x-\beta^{-1}\int_Imx^{m-1}P_cu(\beta u_x-P)\\&+C\int_Ix^m(
u_x^2+x^{-2}u^2)+C\int_Ix^m\rho(1+\theta^{2+2r})+C\\
\le&-\beta^{-1}\int_Ix^mP_cu(\beta u_x-P)_x+C\int_Ix^m(
u_x^2+x^{-2}u^2)+C\int_Ix^m\rho(1+\theta^{2+2r})+C.
\end{split}
\eex Using (\ref{non-1.2})$_2$, Lemmas \ref{le:energy
identity}-\ref{le:upper bound of density}, Corollary \ref{cor:tx
H^2 of v}, ($A_3$), Poincar$\mathrm{\acute{e}}$ inequality, Young
inequality and $\alpha<q-r$, we have \bex\begin{split} II_{4,2}
=&-\beta^{-1}\int_Ix^mP_cu(\rho u_t+\rho uu_x-\frac{\rho
v^2}{x}-\frac{m\beta u_x}{x}+\frac{m\beta u}{x^2})\\&+C\int_Ix^m(
u_x^2+x^{-2}u^2)+C\int_Ix^m\rho(1+\theta^{2+2r})+C\\
\le& \frac{1}{4}\int_Ix^m\rho u_t^2+C\int_Ix^m\rho
u^2+C\|u\|_{L^\infty}^2\|\rho\|_{L^\infty}\int_Ix^mu_x^2+C\|v\|_{L^\infty}^2\int_Ix^m\rho
v^2\\&+C\int_Ix^m(
u_x^2+x^{-2}u^2)+C\int_Ix^m\rho(1+\theta^{2+2r})+C\\
\le& \frac{1}{4}\int_Ix^m\rho u_t^2+C\left(\int_Ix^m(
u_x^2+x^{-2}u^2)\right)^2+C\sup\limits_{x\in
I}\theta^{q-\alpha+1}+C.
\end{split}
\eex Putting all the estimates about $II_{4,1}$ and $II_{4,2}$
into $II_4$, we have \be\label{II_4}
\begin{split}
II_4\le&-\frac{1}{2\beta}\frac{d}{dt}\int_Ix^mP^2+\frac{1}{2}\int_Ix^m\rho
u_t^2+C\int_Ix^m(\kappa\theta_x)^2+C\left(\int_Ix^m(u_x^2+x^{-2}u^2)\right)^2
\\&+C\sup\limits_{x\in
I}(1+\theta^{q-\alpha+1})\left[1+\int_Ix^m\left(u_x^2+\rho\theta^{q+r+2}+x^{-2}u^2\right)\right]+C.
\end{split}
\ee For $II_5$, recalling $P=\rho Q+P_c$, we have
\be\label{II_5}\begin{split}
II_5=&m\int_Ix^{m-1}\rho Qu_t+m\int_Ix^{m-1}P_c(\rho)u_t\\
\le&\frac{1}{8}\int_Ix^m\rho
u_t^2+C\int_Ix^m\rho(1+\theta^{2+2r})+C\int_Ix^m\rho\\
\le&\frac{1}{8}\int_Ix^m\rho
u_t^2+C(\|\theta\|_{L^\infty}^{q-\alpha+1}+1)\int_Ix^m\rho\theta^{1+r}+C\\
\le&\frac{1}{8}\int_Ix^m\rho
u_t^2+C\|\theta\|_{L^\infty}^{q-\alpha+1}+C,
\end{split}
\ee where we have used Lemma \ref{le:energy identity}, Young
inequality and ($A_4$). Note that $P_c(0)=0$ by (\ref{non-1.1})
and ($A_2$).

Putting (\ref{II3}), (\ref{II_4}) and (\ref{II_5}) into (\ref{d
ux-1}), we have \bex\begin{split} &\frac{3}{4}\int_Ix^m\rho
u_t^2+\frac{\beta}{2}\frac{d}{dt}\int_Ix^m(u_x^2+mx^{-2}u^2)\\
\le&C\left(\int_Ix^m(u_x^2+mx^{-2}u^2)\right)^2+\frac{d}{dt}\int_Ix^mPu_x-\frac{1}{2\beta}\frac{d}{dt}\int_Ix^mP^2+\frac{5}{8}\int_Ix^m\rho
u_t^2+C\int_Ix^m(\kappa\theta_x)^2\\&+C\sup\limits_{x\in
I}(1+\theta^{q-\alpha+1})\left[1+\int_Ix^m\left(u_x^2+\rho\theta^{q+r+2}+mx^{-2}u^2\right)\right]+C.
\end{split}
\eex Thus, \be\label{rho ut1}\begin{split}
&\frac{1}{8}\int_Ix^m\rho
u_t^2+\frac{\beta}{2}\frac{d}{dt}\int_Ix^m(u_x^2+mx^{-2}u^2)\\
\le&C\left(\int_Ix^m(u_x^2+mx^{-2}u^2)\right)^2+\frac{d}{dt}\int_Ix^mPu_x-\frac{1}{2\beta}\frac{d}{dt}\int_Ix^mP^2
+C\int_Ix^m(\kappa\theta_x)^2\\&+C\sup\limits_{x\in
I}(1+\theta^{q-\alpha+1})\left[1+\int_Ix^m\left(u_x^2+\rho\theta^{q+r+2}+mx^{-2}u^2\right)\right]+C.
\end{split}
\ee Integrating (\ref{rho ut1}) over $(0,t)$ for $t\in[0,T]$, and
using Corollary \ref{cor:int0T theta 1} and Cauchy inequality, we
have \bex\begin{split} &\frac{1}{8}\int_0^t\int_Ix^m\rho
u_t^2+\frac{\beta}{2}\int_Ix^m(u_x^2+mx^{-2}u^2)\\
\le&C\int_0^t\left(\int_Ix^m(u_x^2+mx^{-2}u^2)\right)^2+\int_Ix^mPu_x
+C\int_0^t\int_Ix^m(\kappa\theta_x)^2\\&+C\int_0^t\sup\limits_{x\in
I}(1+\theta^{q-\alpha+1})\int_Ix^m\left(u_x^2+\rho\theta^{q+r+2}+mx^{-2}u^2\right)+C\\
\le&C\int_0^t\left(\int_Ix^m(u_x^2+mx^{-2}u^2)\right)^2+\frac{\beta}{4}\int_Ix^mu_x^2+C\int_Ix^m(\rho^2Q^2+P_c^2)
+C\int_0^t\int_Ix^m(\kappa\theta_x)^2\\&+C\int_0^t\sup\limits_{x\in
I}(1+\theta^{q-\alpha+1})\int_Ix^m\left(u_x^2+\rho\theta^{q+r+2}+mx^{-2}u^2\right)+C.
\end{split}
\eex Thus, \be\label{rho ut2}\begin{split} &\int_0^t\int_Ix^m\rho
u_t^2+\int_Ix^m(u_x^2+mx^{-2}u^2)\\
\le&C\int_0^t\left(\int_Ix^m(u_x^2+mx^{-2}u^2)\right)^2+C\int_Ix^m\rho\theta^{q+r+2}
+C\int_0^t\int_Ix^m(\kappa\theta_x)^2\\&+C\int_0^t\sup\limits_{x\in
I}(1+\theta^{q-\alpha+1})\int_Ix^m\left(u_x^2+\rho\theta^{q+r+2}+mx^{-2}u^2\right)+C.
\end{split}
\ee The next step is to handle the second term and the third one
on the right hand side of (\ref{rho ut2}).

Multiplying (\ref{non-3.6}) by
$x^m\int_0^\theta\kappa(\xi)\,d\xi$, and integrating by parts over
$I$, we have \be\label{dt rho theta-1}\begin{split}
&\frac{d}{dt}\int_Ix^m\rho\left(\int_0^\theta
Q^\prime(\eta)\int_0^\eta\kappa(\xi)d\xi
d\eta\right)+\int_Ix^m(\kappa\theta_x)^2\\=&-\int_Ix^m\rho\theta
Q^\prime(\theta)(u_x+mx^{-1}u)\int_0^\theta\kappa(\xi)\,d\xi+\lambda\int_Ix^m(u_x+mx^{-1}u)^2\int_0^\theta\kappa(\xi)\,d\xi
\\&+\mu\int_Ix^m\left(w_x^2+2u_x^2+(v_x-mx^{-1}v)^2+2mx^{-2}u^2\right)\int_0^\theta\kappa(\xi)\,d\xi\\
\les&\|\theta(1+\theta^q)\|_{L^\infty}\left(\int_Ix^m\rho(1+\theta^{1+r})|u_x+mx^{-1}u|+\int_Ix^m(w_x^2+u_x^2+v_x^2+x^{-2}v^2+x^{-2}u^2)\right)\\
\les&\|\theta(1+\theta^q)\|_{L^\infty}\left(1+(\int_Ix^m\rho\theta^{2r+2})^\frac{1}{2}\left(\int_Ix^m(u_x^2+x^{-2}u^2)\right)^\frac{1}{2}+
\int_Ix^m(u_x^2+x^{-2}u^2)\right),
\end{split}
\ee where we have used ($A_4$), ($A_5$), Lemmas \ref{le:energy
identity}, \ref{le:upper bound of density}, \ref{le:H^1 of v} and
\ref{le:H^1 of w}, H\"older inequality and Cauchy inequality.

 By Corollary \ref{non-cor:2.1}, we have
\be\label{theta(1+theta q)}\begin{split}
\|\theta(1+\theta^q)\|_{L^\infty}\les \|\kappa\theta_x\|_{L^2}+1.
\end{split}
\ee Putting (\ref{theta(1+theta q)}) into (\ref{dt rho theta-1}),
and using Lemma \ref{le:energy identity} and Cauchy inequality, we
have \be\label{dt rho theta-2}\begin{split}
&\frac{d}{dt}\int_Ix^m\rho\left(\int_0^\theta
Q^\prime(\eta)\int_0^\eta\kappa(\xi)d\xi
d\eta\right)+\int_Ix^m(\kappa\theta_x)^2\\
\le&C(\|\kappa\theta_x\|_{L^2}+1)\left(1+(\int_Ix^m\rho\theta^{2r+2})^\frac{1}{2}\left(\int_Ix^m(u_x^2+x^{-2}u^2)\right)^\frac{1}{2}+
\int_Ix^m(u_x^2+x^{-2}u^2)\right)\\
\le&\frac{1}{2}\int_Ix^m(\kappa\theta_x)^2+C\|\theta\|_{L^\infty}^{r+1}\int_Ix^m(u_x^2+x^{-2}u^2)+C\left(\int_Ix^m(u_x^2+x^{-2}u^2)\right)^2+C.
\end{split}
\ee Integrating (\ref{dt rho theta-2}) over $(0,t)$, and using
($A_4$), ($A_5$) and Young inequality, we have \be\label{rho theta
q+r+2}\begin{split}
&\int_Ix^m\rho\theta^{q+r+2}+\int_0^t\int_Ix^m(\kappa\theta_x)^2\\
\le&
C\int_0^t\|\theta\|_{L^\infty}^{q-\alpha+1}\int_Ix^m(u_x^2+x^{-2}u^2)+C\int_0^t\left(\int_Ix^m(u_x^2+x^{-2}u^2)\right)^2+C.
\end{split}
\ee

By (\ref{rho ut2}), (\ref{rho theta q+r+2}), Corollary
\ref{cor:int0T theta 1}, Lemma \ref{le:int0T H^1 of uvw} and
Gronwall inequality, we complete the proof of Lemma \ref{le:H^1 of
u}.
\endpf \\

By Poincar$\mathrm{\acute{e}}$ inequality and Lemma \ref{le:H^1 of
u}, we get the next corollary.
\begin{corollary}\label{cor:xt L^infty of u} Under the conditions of
Theorem \ref{th:1.1}, we have \bex \|u\|_{L^\infty(Q_T)}\le C.
\eex
\end{corollary}
\begin{corollary}\label{cor:int0T theta 2} Under the conditions of
Theorem \ref{th:1.1}, we have \bex
\int_0^T\|\theta\|_{L^\infty}^{2q+2}\le C. \eex
\end{corollary}
\pf By Corollary \ref{non-cor:2.1} and Cauchy inequality, we have
\bex\begin{split} \int_0^T\sup\limits_{x\in
I}\theta^{2q+2}\les&\int_0^T\int_I\theta^{2q+1}|\theta_x|+1\\
\le&\frac{1}{2}\int_0^T\sup\limits_{x\in
I}\theta^{2q+2}+C\int_{Q_T}\theta^{2q}\theta_x^2.
\end{split}
\eex This together with Lemma \ref{le:H^1 of u} completes the
proof of Corollary \ref{cor:int0T theta 2}.
\endpf

\begin{lemma}\label{le:H^1 of rho}Under the conditions of
Theorem \ref{th:1.1}, we have for any $t\in[0,T]$ \bex
\int_Ix^m(\rho_x^2+\rho_t^2)+\int_{Q_T}x^mu_{xx}^2\le C. \eex
\end{lemma}
\pf Differentiating (\ref{non-1.2})$_1$, we have
\be\label{rhoxt+.}\begin{split}
\rho_{xt}+\rho_{xx}u+2\rho_xu_x+\rho
u_{xx}+mx^{-1}\rho_xu+mx^{-1}\rho u_x-mx^{-2}\rho u=0.
\end{split}
\ee Multiplying (\ref{rhoxt+.}) by $2x^m\rho_x$, and integrating
by parts over $I$, we have \be\label{dt rhox-1}\begin{split}
\frac{d}{dt}\int_Ix^m\rho_x^2=&-4\int_Ix^m\rho_x^2u_x-\int_Ix^m(\rho_x^2)_xu-2\int_Ix^m\rho\rho_xu_{xx}\\&-2m\int_Ix^{m-1}\rho_x^2u-2m\int_Ix^{m-1}\rho\rho_xu_x
+2m\int_Ix^{m-2}\rho\rho_xu\\=&-3\int_Ix^m\rho_x^2u_x-m\int_Ix^{m-1}\rho_x^2u-2\int_Ix^m\rho\rho_xu_{xx}\\&-2m\int_Ix^{m-1}\rho\rho_xu_x+2m\int_Ix^{m-2}\rho\rho_xu\\=
&\sum\limits_{i=1}^5III_i.
\end{split}
\ee For $III_1$, using Sobolev inequality, Young inequality,
($A_2$), ($A_3$), ($A_4$), and Lemmas \ref{le:energy identity},
\ref{le:upper bound of density}, \ref{le:H^1 of u}, we have
\be\label{III 1-1}\begin{split}
III_1=&-3\beta^{-1}\int_Ix^m\rho_x^2(\beta
u_x-P)-3\beta^{-1}\int_Ix^m\rho_x^2P\\ \les&\left(\|\beta
u_x-P\|_{L^\infty}+\|\rho
Q+P_c\|_{L^\infty}\right)\int_Ix^m\rho_x^2\\ \les&\left(\|\beta
u_x-P\|_{L^2}+\|(\beta
u_x-P)_x\|_{L^2}+\|1+\theta^{r+1}\|_{L^\infty}\right)\int_Ix^m\rho_x^2\\
\les&\left(1+\sup\limits_{x\in I}\theta^{q-\alpha+1}+\|(\beta
u_x-P)_x\|_{L^2}\right)\int_Ix^m\rho_x^2.
\end{split}
\ee For $\|(\beta u_x-P)_x\|_{L^2}$, using (\ref{non-1.2})$_2$,
Lemmas \ref{le:upper bound of density}, \ref{le:H^1 of u},
Corollary \ref{cor:tx H^2 of v} and Corollary \ref{cor:xt L^infty
of u}, we have \be\label{beta ux-P}\begin{split}
 \|(\beta
u_x-P)_x\|_{L^2}\les&\|\rho u_t\|_{L^2}+\|\rho
uu_x\|_{L^2}+\|x^{-1}\rho
v^2\|_{L^2}+\|x^{-1}u_x\|_{L^2}+\|x^{-2}u\|_{L^2}\\
\les&\|\sqrt{\rho} u_t\|_{L^2}+1.
\end{split}
\ee Substituting (\ref{beta ux-P}) into (\ref{III 1-1}), we obtain
\be\label{III 1}\begin{split} III_1 \les&\left(1+\sup\limits_{x\in
I}\theta^{q-\alpha+1}+\|\sqrt{\rho}
u_t\|_{L^2}\right)\int_Ix^m\rho_x^2.
\end{split}
\ee For $III_i$, $i=$2, 3, 4, 5, we have \be\label{III
2}\begin{split} III_2\les
\|u\|_{L^\infty}\int_Ix^m\rho_x^2\les\int_Ix^m\rho_x^2.
\end{split}
\ee \be\label{III 3}
\begin{split}
III_3\les\int_Ix^m\rho_x^2+\int_Ix^mu_{xx}^2.
\end{split}
\ee \be\label{III 4}
\begin{split}
III_4\les\int_Ix^m\rho_x^2+\int_Ix^mu_x^2\les \int_Ix^m\rho_x^2+1.
\end{split}
\ee \be\label{III 5}
\begin{split}
III_5\les\int_Ix^m\rho_x^2+\int_Ix^{m-2}u^2\les
\int_Ix^m\rho_x^2+1.
\end{split}
\ee Putting (\ref{III 1})-(\ref{III 5}) into (\ref{dt rhox-1}), we
have \be\label{dt rhox}\begin{split}
\frac{d}{dt}\int_Ix^m\rho_x^2\les\left(1+\sup\limits_{x\in
I}\theta^{q-\alpha+1}+\|\sqrt{\rho}
u_t\|_{L^2}\right)\int_Ix^m\rho_x^2+\int_Ix^mu_{xx}^2+1.
\end{split}\ee
By (\ref{beta ux-P}), ($A_2$), ($A_3$), ($A_4$) and Lemma
\ref{le:upper bound of density}, we have \bex\begin{split}
\int_Ix^mu_{xx}^2\les& \int_Ix^m\rho u_t^2+
\int_Ix^m\rho_x^2Q^2+\int_Ix^m\rho^2(Q^\prime)^2\theta_x^2+\int_Ix^m\rho_x^2+1\\
\les&\int_Ix^m\rho u_t^2+ (1+\sup\limits_{x\in
I}\theta^{2+2r})\int_Ix^m\rho_x^2+\int_Ix^m(1+\theta^{2r})\theta_x^2+1.
\end{split}
\eex Since $q>r$, using Young inequality, we have
\be\label{uxx}\begin{split} \int_Ix^mu_{xx}^2\les& \int_Ix^m\rho
u_t^2+ (1+\sup\limits_{x\in
I}\theta^{2+2q})\int_Ix^m\rho_x^2+\int_Ix^m(1+\theta^q)^2\theta_x^2+1.
\end{split}
\ee Substituting (\ref{uxx}) into (\ref{dt rhox}), and using
Corollary \ref{cor:int0T theta 1}, Lemma \ref{le:H^1 of u},
Corollary \ref{cor:int0T theta 2}, and Gronwall inequality, we get
\be\label{rhox} \int_Ix^m\rho_x^2\le C. \ee Substituting
(\ref{rhox}) into (\ref{uxx}), and using Lemma \ref{le:H^1 of u}
and Corollary \ref{cor:int0T theta 2} again, we have
$$\int_{Q_T}x^mu_{xx}^2\le C.$$ The estimate for $\rho_t$ can be
obtained easily by (\ref{non-1.2})$_1$.
\endpf\\

{\noindent\bf Step 4: $H^2$ estimates of $(\rho,u,v,w)$ and $H^1$
estimates of $\theta$}

\begin{lemma}\label{le:rho vt2}Under the conditions of Theorem
\ref{th:1.1}, we have for any $t\in[0,T]$ \bex \int_Ix^m\rho
v_t^2+\int_{Q_T}x^m\left( v_{xt}^2+x^{-2}v_t^2\right)\le C. \eex
\end{lemma}
\pf Differentiating (\ref{non-1.2})$_3$ with respect to t, we have
\be\label{rho vtt}\begin{split} &\rho
v_{tt}+\rho_tv_t+\rho_tuv_x+\rho u_t v_x+\rho
uv_{xt}+x^{-1}\rho_tuv+x^{-1}\rho u_tv+x^{-1}\rho
uv_t\\=&\mu(v_{xxt}+mx^{-1}v_{xt}-mx^{-2}v_t). \end{split}\ee
Multiplying (\ref{rho vtt}) by $x^mv_t$, and integrating by parts
over $I$, we have \be\label{dtrho vt}\begin{split}
&\frac{1}{2}\frac{d}{dt}\int_Ix^m\rho
v_t^2+\mu\int_Ix^m(v_{xt}^2+mx^{-2}v_t^2)\\=&-\frac{1}{2}\int_Ix^m\rho_tv_t^2-\int_Ix^m(\rho_tuv_x+\frac{\rho_tuv}{x})v_t-\int_Ix^m\left(\rho
u_tv_x+\rho uv_{xt}+\frac{\rho u_tv+\rho
uv_t}{x}\right)v_t\\=&\sum\limits_{i=1}^3IV_i.
\end{split}
\ee For $IV_1$, using (\ref{non-1.2})$_1$, integration by parts,
Cauchy inequality, Lemma \ref{le:upper bound of density} and
Corollary \ref{cor:xt L^infty of u}, we have
\be\label{IV1}\begin{split} IV_1=&\frac{1}{2}\int_I(x^m\rho
u)_xv_t^2\\=&-\int_Ix^m\rho uv_tv_{xt}\\
\le&\frac{\mu}{6}\int_Ix^mv_{xt}^2+C\int_Ix^m\rho v_t^2.
\end{split}
\ee Similarly, for $IV_2$, we have \be\label{IV2}
\begin{split}IV_2=&\int_I(x^m\rho
u)_x(uv_x+\frac{uv}{x})v_t\\=&-\int_Ix^m\rho
u(u_xv_x+uv_{xx}+\frac{u_xv+uv_x}{x}-\frac{uv}{x^2})v_t-\int_Ix^m\rho
u(uv_x+\frac{uv}{x})v_{xt}\\ \le&C\int_Ix^m\rho
v_t^2+C\int_Ix^m\rho\left(u^2u_x^2v_x^2+u^4v_{xx}^2+\frac{u^2u_x^2v^2+u^4v_x^2}{x^2}+\frac{u^4v^2}{x^4}\right)\\&+\frac{\mu}{6}\int_Ix^mv_{xt}^2
+C\int_Ix^m\rho^2u^2(u^2v_x^2+\frac{u^2v^2}{x^2})\\
\le&C\int_Ix^m\rho
v_t^2+C\|v_x\|_{L^\infty}^2\int_Ix^mu_x^2+C\int_Ix^mv_{xx}^2+C\|u_x\|_{L^\infty}^2\int_Ix^{m-2}v^2\\&+C\int_Ix^mv_x^2
+C\int_Ix^{m-2}v^2+\frac{\mu}{6}\int_Ix^mv_{xt}^2\\
\le&C\int_Ix^m\rho
v_t^2+C\int_Ix^mv_{xx}^2+C\int_Ix^mu_{xx}^2+\frac{\mu}{6}\int_Ix^mv_{xt}^2+C,
\end{split}
\ee where we have used (\ref{non-1.2})$_1$, integration by parts,
Cauchy inequality, Lemmas \ref{le:upper bound of density},
\ref{le:H^1 of v}, \ref{le:H^1 of u},  Corollary \ref{cor:xt
L^infty of u} and Sobolev inequality.

For $IV_3$, using Cauchy inequality, Lemmas \ref{le:upper bound of
density}, \ref{le:H^1 of v}, Corollary \ref{cor:tx H^2 of v},
Corollary \ref{cor:xt L^infty of u} and Sobolev inequality, we get
 \be\label{IV3}
\begin{split}IV_3=&-\int_Ix^m\rho u_tv_xv_t-\int_Ix^m\rho
uv_{xt}v_t-\int_Ix^{m-1}\rho u_tvv_t-\int_Ix^{m-1}\rho uv_t^2\\
\le&C\int_Ix^m\rho u_t^2+C\int_Ix^m\rho
v_x^2v_t^2+\frac{\mu}{6}\int_Ix^mv_{xt}^2+C\int_Ix^m\rho^2u^2v_t^2+C\int_Ix^{m-2}\rho
v^2v_t^2\\&+C\int_Ix^m\rho v_t^2\\ \le&C\int_Ix^m\rho
u_t^2+\frac{\mu}{6}\int_Ix^mv_{xt}^2+C(1+\int_Ix^mv_{xx}^2)\int_Ix^m\rho
v_t^2.
\end{split}
\ee Putting (\ref{IV1}), (\ref{IV2}) and (\ref{IV3}) into
(\ref{dtrho vt}), we have \be\label{dt rho vt-1}
\begin{split}
&\frac{d}{dt}\int_Ix^m\rho
v_t^2+\mu\int_Ix^m(v_{xt}^2+mx^{-2}v_t^2)\\
\le&C(1+\int_Ix^mv_{xx}^2)\int_Ix^m\rho
v_t^2+C\int_Ix^mv_{xx}^2+C\int_Ix^mu_{xx}^2+C\int_Ix^m\rho
u_t^2+C.
\end{split}
\ee It follows from (\ref{dt rho vt-1}), Corollary \ref{cor:tx H^2
of v}, Lemmas \ref{le:H^1 of u}, \ref{le:H^1 of rho}, the
compatibility conditions and Gronwall inequality, we complete the
proof of Lemma \ref{le:rho vt2}.
\endpf \\

By (\ref{non-1.2})$_3$, Lemmas \ref{le:upper bound of density},
\ref{le:H^1 of v}, \ref{le:H^1 of u}, \ref{le:H^1 of rho},
\ref{le:rho vt2}, Corollary \ref{cor:tx H^2 of v}, Corollary
\ref{cor:xt L^infty of u}, and Sobolev inequality, we get the next
estimate.
\begin{corollary}\label{cor:vxx}Under the conditions of Theorem
\ref{th:1.1}, we have for any $t\in[0,T]$ \bex
\|v\|_{W^{1,\infty}(Q_T)}+\int_Ix^m
v_{xx}^2+\int_{Q_T}x^mv_{xxx}^2\le C. \eex
\end{corollary}
Similar to Lemma \ref{le:rho vt2} and Corollary \ref{cor:vxx}, we
obtain the next lemma and the next corollary.
\begin{lemma}\label{le:rho wt}Under the conditions of Theorem
\ref{th:1.1}, we have for any $t\in[0,T]$ \bex \int_Ix^m\rho
w_t^2+\int_{Q_T}x^m w_{xt}^2\le C. \eex
\end{lemma}
\begin{corollary}\label{cor:wxx}Under the conditions of Theorem
\ref{th:1.1}, we have for any $t\in[0,T]$ \bex
\|w\|_{W^{1,\infty}(Q_T)}+\int_Ix^m
w_{xx}^2+\int_{Q_T}x^mw_{xxx}^2\le C. \eex
\end{corollary}
\begin{lemma}\label{le:rho ut}Under the conditions of Theorem
\ref{th:1.1}, we have for any $t\in[0,T]$ \bex \int_Ix^m\left(\rho
u_t^2+(1+\theta^q)^2\theta_x^2\right)+\int_{Q_T}x^m
\left(u_{xt}^2+x^{-2}u_t^2+\rho(1+\theta^{q+r})\theta_t^2\right)\le
C. \eex
\end{lemma}
\pf Differentiating (\ref{non-1.2})$_2$ with respect to $t$, we
have \be\label{rho utt+.=}\rho u_{tt}+\rho_tu_t+\rho_tuu_x+\rho
u_tu_x+\rho uu_{xt}-\frac{\rho_tv^2}{x}-\frac{2\rho
vv_t}{x}+P_{xt}=\beta
\left(u_{xxt}+\frac{mu_{xt}}{x}-\frac{mu_t}{x^2}\right). \ee

Multiplying (\ref{rho utt+.=}) by $x^m u_t$, and integrating by
parts over $I$, we have \be\label{dt rho ut-1}\begin{split}
&\frac{1}{2}\frac{d}{dt}\int_Ix^m\rho
u_t^2+\beta\int_Ix^m(u_{xt}^2+\frac{mu_t^2}{x^2})\\=&-\frac{1}{2}\int_Ix^m\rho_tu_t^2-\int_Ix^m\rho_tuu_xu_t+\int_Ix^{m-1}\rho_tv^2u_t-\int_Ix^m\rho(u_tu_x+uu_{xt}-
\frac{2vv_t}{x})u_t\\&+\int_IP_t(x^mu_{xt}+mx^{m-1}u_t)\\=&\sum\limits_{i=1}^5V_i.
\end{split}
\ee  For $V_1$, using (\ref{non-1.2})$_1$, integration by parts,
Cauchy inequality, Lemma \ref{le:upper bound of density},
Corollary \ref{cor:xt L^infty of u}, we have
\be\label{V1}\begin{split} V_1=&\frac{1}{2}\int_I(x^m\rho
u)_xu_t^2\\=&-\int_Ix^m\rho
uu_tu_{xt}\\
\le&\frac{\beta}{8}\int_Ix^mu_{xt}^2+C\int_Ix^m\rho u_t^2.
\end{split}
\ee For $V_2$ and $V_3$, using (\ref{non-1.2})$_1$, integration by
parts and Cauchy inequality again, along with Sobolev inequality,
Lemmas \ref{le:upper bound of density}, \ref{le:H^1 of v},
\ref{le:H^1 of u}, Corollary \ref{cor:tx H^2 of v} and Corollary
\ref{cor:xt L^infty of u}, we get \be\label{V2 and
V3}\begin{split} V_2+V_3=&\int_I(x^m \rho
u)_x(uu_xu_t-\frac{v^2u_t}{x})\\=&-\int_Ix^m\rho
u(u_x^2u_t+uu_{xx}u_t+uu_xu_{xt}-\frac{2vv_xu_t+v^2u_{xt}}{x}+\frac{v^2u_t}{x^2})\\
\le&C\int_Ix^m\rho u_t^2u_x^2+C\int_Ix^m\rho
u^2u_x^2+C\int_Ix^m\rho u_t^2+C\int_Ix^m\rho
u^4u_{xx}^2+\frac{\beta}{8}\int_Ix^mu_{xt}^2\\&+C\int_Ix^m\rho^2u^4u_x^2+C\int_Ix^m\rho
u^2v^2v_x^2+C\int_Ix^m\rho^2u^2v^4+C\int_Ix^m\rho u^2v^4\\
\le&C(1+\int_Iu_{xx}^2)\int_Ix^m\rho
u_t^2+C\int_Ix^mu_{xx}^2+\frac{\beta}{8}\int_Ix^mu_{xt}^2+C.
\end{split}
\ee For $V_4$, we have \be\label{V4}\begin{split}
V_4\le&\|u_x\|_{L^\infty}\int_Ix^m\rho
u_t^2+\frac{\beta}{8}\int_Ix^mu_{xt}^2+C\int_Ix^m\rho^2u^2u_t^2+C\int_Ix^m\rho
u_t^2+C\int_Ix^{m-2}\rho v^2v_t^2\\
\le&\frac{\beta}{8}\int_Ix^mu_{xt}^2+C(1+\|u_x\|_{H^1})\int_Ix^m\rho
u_t^2+C,
\end{split}
\ee where we have used Sobolev inequality, Lemmas \ref{le:upper
bound of density}, \ref{le:rho vt2}, Corollary \ref{cor:tx H^2 of
v} and Corollary \ref{cor:xt L^infty of u}.

For $V_5$, using Young inequality, ($A_2$), ($A_3$), ($A_4$), and
Lemmas \ref{le:upper bound of density}, \ref{le:H^1 of rho}, we
have

\be\label{V5}\begin{split} V_5\le
&\frac{\beta}{8}\int_Ix^m(u_{xt}^2+\frac{mu_t^2}{x^2})+C\int_Ix^mP_t^2\\
\le&\frac{\beta}{8}\int_Ix^m(u_{xt}^2+\frac{mu_t^2}{x^2})+C\int_Ix^m\rho_t^2(\partial_\rho
P)^2+C\int_Ix^m\theta_t^2(\partial_\theta P)^2\\
\le&\frac{\beta}{8}\int_Ix^m(u_{xt}^2+\frac{mu_t^2}{x^2})+C(1+\sup\limits_{x\in
I}\theta^{2+2r})\int_Ix^m\rho_t^2+C\int_Ix^m\rho^2(1+\theta^{2r})\theta_t^2\\
\le&\frac{\beta}{8}\int_Ix^m(u_{xt}^2+\frac{mu_t^2}{x^2})+C\sup\limits_{x\in
I}\theta^{2+2q}+C\int_Ix^m\rho(1+\theta^{q+r})\theta_t^2+C.
\end{split}
\ee Putting (\ref{V1}), (\ref{V2 and V3}), (\ref{V4}) and
(\ref{V5}) into (\ref{dt rho ut-1}), we have \be\label{dt rho
ut-2}\begin{split} &\frac{d}{dt}\int_Ix^m\rho
u_t^2+\beta\int_Ix^m(u_{xt}^2+\frac{mu_t^2}{x^2})\\
\le&C(1+\int_Iu_{xx}^2)\int_Ix^m\rho
u_t^2+C\int_Ix^mu_{xx}^2+C\sup\limits_{x\in
I}\theta^{2+2q}+C\int_Ix^m\rho(1+\theta^{q+r})\theta_t^2+C.
\end{split}
\ee Integrating (\ref{dt rho ut-2}) over $(0,t)$, and using the
compatibility conditions, Lemma \ref{le:H^1 of rho} and Corollary
\ref{cor:int0T theta 2}, we obtain \be\label{dt rho
ut-3}\begin{split} &\int_Ix^m\rho
u_t^2+\int_0^t\int_Ix^m(u_{xt}^2+\frac{mu_t^2}{x^2})\\
\le&C\int_0^t(1+\int_Iu_{xx}^2)\int_Ix^m\rho
u_t^2+C\int_0^t\int_Ix^m\rho(1+\theta^{q+r})\theta_t^2+C.
\end{split}
\ee The next step is to estimate the second integral of the right
hand side of (\ref{dt rho ut-3}).

Multiplying (\ref{non-3.6}) by
$x^m\left(\int_0^\theta\kappa(\xi)\,d\xi\right)_t$ (i.e.,
$x^m\kappa(\theta)\theta_t$), and integrating by parts over $I$,
we have \be\label{dt thetax-1}\begin{split} &\int_Ix^m\rho
Q^\prime\kappa\theta_t^2+\frac{1}{2}\frac{d}{dt}\int_Ix^m\kappa^2\theta_x^2\\=&-\int_Ix^m\rho
uQ^\prime\theta_x\kappa\theta_t-\int_Ix^m\rho\theta
Q^\prime(u_x+\frac{mu}{x})\kappa\theta_t+\lambda\int_Ix^m(u_x+\frac{mu}{x})^2\left(\int_0^\theta\kappa(\xi)\,d\xi\right)_t
\\&+\mu\int_Ix^m\left(w_x^2+2u_x^2+(v_x-\frac{mv}{x})^2+\frac{2mu^2}{x^2}\right)\left(\int_0^\theta\kappa(\xi)\,d\xi\right)_t\\
=&\sum\limits_{i=1}^4VI_i.
\end{split}
\ee For $VI_1$, using Cauchy inequality, ($A_4$), ($A_5$), Lemma
\ref{le:upper bound of density} and Corollary \ref{cor:xt L^infty
of u}, we have \be\label{VI1}\begin{split}
VI_1\le&\frac{1}{4}\int_Ix^m\rho
Q^\prime\kappa\theta_t^2+C\int_Ix^m\rho Q^\prime\kappa
u^2\theta_x^2\\ \le&\frac{1}{4}\int_Ix^m\rho
Q^\prime\kappa\theta_t^2+C\int_Ix^m(1+\theta^q)^2\theta_x^2.
\end{split}
\ee For $VI_2$, using Cauchy inequality, ($A_4$) and ($A_5$)
again, along with Lemmas \ref{le:energy identity}, \ref{le:H^1 of
u} and Corollary \ref{cor:xt L^infty of u}, we have
\be\label{VI2}\begin{split} VI_2\le&\frac{1}{4}\int_Ix^m\rho
Q^\prime\kappa\theta_t^2+C\int_Ix^m\rho
Q^\prime\kappa\theta^2(u_x^2+\frac{u^2}{x^2})\\
\le&\frac{1}{4}\int_Ix^m\rho
Q^\prime\kappa\theta_t^2+C(\|u_x\|_{L^\infty}^2+\|u\|_{L^\infty}^2)\int_Ix^m\rho(1+\theta^{q+r+2})\\
\le&\frac{1}{4}\int_Ix^m\rho
Q^\prime\kappa\theta_t^2+C\int_Ix^mu_{xx}^2+C.
\end{split}
\ee For $VI_3$, we have \be\label{VI3}\begin{split}
VI_3\le&\frac{d}{dt}\left(\lambda\int_Ix^m(u_x+\frac{mu}{x})^2\int_0^\theta\kappa(\xi)\,d\xi\right)-2\lambda\int_Ix^m(u_x+\frac{mu}{x})(u_{xt}+\frac{mu_t}{x})\int_0^\theta\kappa(\xi)\,d\xi
\\
\le&\lambda\frac{d}{dt}\left(\int_Ix^m(u_x+\frac{mu}{x})^2\int_0^\theta\kappa(\xi)\,d\xi\right)
+C\|\theta(1+\theta^q)\|_{L^\infty}\left(\int_Ix^m(u_{xt}^2+\frac{mu_t^2}{x^2})\right)^\frac{1}{2}\\
\le&\lambda\frac{d}{dt}\left(\int_Ix^m(u_x+\frac{mu}{x})^2\int_0^\theta\kappa(\xi)\,d\xi\right)
+C(\|\kappa\theta_x\|_{L^2}+1)\left(\int_Ix^m(u_{xt}^2+\frac{mu_t^2}{x^2})\right)^\frac{1}{2},
\end{split}
\ee where we have used H\"older inequality, Lemma \ref{le:H^1 of
u}, (\ref{theta(1+theta q)}) and ($A_5$).

For $VI_4$, using H\"older inequality, (\ref{theta(1+theta q)})
and Lemma \ref{le:H^1 of u} again, along with Lemmas \ref{le:H^1
of v}, \ref{le:H^1 of w}, we have \be\label{VI4}\begin{split}
&VI_4=\mu\frac{d}{dt}\int_Ix^m\left(w_x^2+2u_x^2+(v_x-\frac{mv}{x})^2+\frac{2mu^2}{x^2}\right)\int_0^\theta\kappa(\xi)\,d\xi
\\&-2\mu\int_Ix^m\left(w_xw_{xt}+2u_xu_{xt}+(v_x-\frac{mv}{x})(v_{xt}-\frac{mv_t}{x})+\frac{2muu_t}{x^2}\right)\int_0^\theta\kappa(\xi)\,d\xi\\
\le&\mu\frac{d}{dt}\int_Ix^m\left(w_x^2+2u_x^2+(v_x-\frac{mv}{x})^2+\frac{2mu^2}{x^2}\right)\int_0^\theta\kappa(\xi)\,d\xi
+C(\|\kappa\theta_x\|_{L^2}+1)\left(\int_Ix^mw_{xt}^2\right)^\frac{1}{2}\\&+C(\|\kappa\theta_x\|_{L^2}+1)\left[\left(\int_Ix^mu_{xt}^2\right)^\frac{1}{2}
+\left(\int_Ix^m(v_{xt}^2+\frac{v_t^2}{x^2})\right)^\frac{1}{2}+\left(\int_Ix^{m-2}u_t^2\right)^\frac{1}{2}\right].
\end{split}
\ee Substituting (\ref{VI1}), (\ref{VI2}), (\ref{VI3}) and
(\ref{VI4}) into (\ref{dt thetax-1}), we have \be\label{dt
thetax-2}\begin{split} &\int_Ix^m\rho
Q^\prime\kappa\theta_t^2+\frac{d}{dt}\int_Ix^m\kappa^2\theta_x^2\\
 \le&2\mu\frac{d}{dt}\int_Ix^m\left(w_x^2+2u_x^2+(v_x-\frac{mv}{x})^2+\frac{2mu^2}{x^2}\right)\int_0^\theta\kappa(\xi)\,d\xi\\&
 +2\lambda\frac{d}{dt}\left(\int_Ix^m(u_x+\frac{mu}{x})^2\int_0^\theta\kappa(\xi)\,d\xi\right)+
 C\int_Ix^m(1+\theta^q)^2\theta_x^2\\&+C\int_Ix^mu_{xx}^2
+C(\|\kappa\theta_x\|_{L^2}+1)\left(\int_Ix^mw_{xt}^2\right)^\frac{1}{2}\\&+C(\|\kappa\theta_x\|_{L^2}+1)\left[\left(\int_Ix^mu_{xt}^2\right)^\frac{1}{2}
+\left(\int_Ix^m(v_{xt}^2+\frac{v_t^2}{x^2})\right)^\frac{1}{2}+\left(\int_Ix^{m-2}u_t^2\right)^\frac{1}{2}\right]+C.
\end{split}
\ee Integrating (\ref{dt thetax-2}) over $(0,t)$, and using
($A_4$), ($A_5$) and Lemma \ref{le:H^1 of rho}, we have
\be\label{dt thetax-3}\begin{split} &\int_0^t\int_Ix^m\rho
(1+\theta^{q+r})\theta_t^2+\int_Ix^m(1+\theta^q)^2\theta_x^2\\
 \les&\int_Ix^m\left(w_x^2+u_x^2+v_x^2+\frac{v^2}{x^2}+\frac{u^2}{x^2}\right)\int_0^\theta\kappa(\xi)\,d\xi+\int_0^t\int_Ix^m(1+\theta^q)^2\theta_x^2\\&
+\int_0^t(\|\kappa\theta_x\|_{L^2}+1)\left(\int_Ix^mw_{xt}^2\right)^\frac{1}{2}+\int_0^t(\|\kappa\theta_x\|_{L^2}+1)\left(\int_Ix^mu_{xt}^2\right)^\frac{1}{2}\\&
+\int_0^t(\|\kappa\theta_x\|_{L^2}+1)\left[\left(\int_Ix^m(v_{xt}^2+\frac{v_t^2}{x^2})\right)^\frac{1}{2}+\left(\int_Ix^{m-2}u_t^2\right)^\frac{1}{2}\right]+1.
\end{split}
\ee Using (\ref{theta(1+theta q)}), (\ref{dt thetax-3}), Cauchy
inequality and Lemmas \ref{le:rho vt2}, \ref{le:rho wt}, we have
\be\label{dt thetax-4}\begin{split} &\int_0^t\int_Ix^m\rho
(1+\theta^{q+r})\theta_t^2+\int_Ix^m(1+\theta^q)^2\theta_x^2\\
\le&C\int_0^t\int_Ix^m(1+\theta^q)^2\theta_x^2
+C\int_0^t(\|\kappa\theta_x\|_{L^2}+1)\left[\left(\int_Ix^mu_{xt}^2\right)^\frac{1}{2}+\left(\int_Ix^{m-2}u_t^2\right)^\frac{1}{2}\right]+C.
\end{split}
\ee By (\ref{dt rho ut-3}), (\ref{dt thetax-4}), Cauchy inequality
and Gronwall inequality, we complete the proof of Lemma
\ref{le:rho ut}.
\endpf

By Lemmas \ref{non-le:2.1} and \ref{le:rho ut}, we get the next
corollary.
\begin{corollary}\label{cor:theta} Under the conditions of Theorem
\ref{th:1.1}, we have \bex \|\theta\|_{L^\infty(Q_T)}\le C. \eex
\end{corollary}
\begin{corollary}\label{cor:xt theta xx}Under the conditions of Theorem
\ref{th:1.1}, we have for any $t\in[0,T]$ \bex
\|u\|_{W^{1,\infty}(Q_T)}+\int_Ix^m
u_{xx}^2+\int_{Q_T}x^m\theta_{xx}^2\le C. \eex \pf By (\ref{uxx}),
Lemma \ref{le:H^1 of rho}, Lemma \ref{le:rho ut} and Corollary
\ref{cor:theta}, we get \be\label{uxx+1} \int_Ix^mu_{xx}^2\le C.
\ee It follows from (\ref{uxx+1}), Lemma \ref{le:H^1 of u},
Corollary \ref{cor:xt L^infty of u} and Sobolev inequality, we
obtain \bex \|u\|_{W^{1,\infty}}\le C. \eex By (\ref{non-3.6}),
(\ref{uxx+1}), ($A_4$), ($A_5$), Lemmas \ref{le:upper bound of
density}, \ref{non-le:2.2}, \ref{le:H^1 of u}, \ref{le:rho ut},
Corollaries \ref{cor:xt L^infty of u}, \ref{cor:wxx},
\ref{cor:theta}, and Cauchy inequality, we have
 \bex\begin{split}
\int_Ix^m\theta_{xx}^2\les&\int_Ix^m\rho\theta_t^2+\int_Ix^m\theta_x^2+\int_Ix^m(u_x^2+\frac{u^2}{x^2})+
\int_Ix^m\theta_x^4+\int_Ix^m(u_x^4+\frac{u^4}{x^4})\\&+\int_Ix^m(w_x^4+u_x^4+v_x^4+\frac{v^4}{x^4}+\frac{u^4}{x^4})\\
\les&\int_Ix^m\rho\theta_t^2+\int_Ix^m\theta_x^4+1\\
\les&\int_Ix^m\rho\theta_t^2+\|\theta_x^2\|_{L^\infty}\int_Ix^m\theta_x^2+1\\
\les&\int_Ix^m\rho\theta_t^2+\|\theta_x\theta_{xx}\|_{L^1}+1\\
\le&C\int_Ix^m\rho\theta_t^2+\frac{1}{2}\int_Ix^m\theta_{xx}^2+C.
\end{split}
\eex This deduces \be\label{int I theta xx1}\begin{split} \int_Ix^m\theta_{xx}^2\le C\int_Ix^m\rho\theta_t^2+C.
\end{split}
\ee
Integrating (\ref{int I theta xx1}) over $(0,T)$, and using Lemma \ref{le:rho ut}, we get
$$
\int_{Q_T}x^m\theta_{xx}^2\le  C.
$$

\endpf

\end{corollary}
\begin{lemma}\label{le: rho xx}Under the conditions of Theorem
\ref{th:1.1}, we have for any $t\in[0,T]$ \bex
\int_Ix^m\rho_{xx}^2+\int_{Q_T}x^mu_{xxx}^2\le C. \eex
\end{lemma}
\pf Differentiating (\ref{rhoxt+.}) with respect to $x$, we have
\be\label{rho xxt+.}\begin{split}
\rho_{xxt}=&-\rho_{xxx}u-3\rho_{xx}u_x-3\rho_xu_{xx}-\rho
u_{xxx}-\frac{m\rho_{xx}u}{x}-\frac{2m\rho_xu_x}{x}+\frac{2m\rho_xu}{x^2}+\frac{2m\rho
u_x}{x^2}\\&-\frac{2m\rho u}{x^3}-\frac{m\rho
u_{xx}}{x}.\end{split} \ee Multiplying (\ref{rho xxt+.}) by
$2x^m\rho_{xx}$, integrating by parts it over $I$, and using
H\"older inequality, we have \be\label{dt rho xx-1}\begin{split}
\frac{d}{dt}\int_Ix^m\rho_{xx}^2=&
-5\int_Ix^m\rho_{xx}^2u_x-m\int_Ix^{m-1}\rho_{xx}^2u-6\int_Ix^m\rho_x\rho_{xx}u_{xx}
-2\int_Ix^m\rho\rho_{xx}u_{xxx}\\&-4m\int_Ix^{m-1}\rho_x\rho_{xx}u_x+4m\int_Ix^{m-2}\rho_x\rho_{xx}u\\&
+4m\int_Ix^{m-2}\rho\rho_{xx}u_x-
4m\int_Ix^{m-3}\rho\rho_{xx}u-2m\int_Ix^{m-1}\rho\rho_{xx}u_{xx}\\=&\sum\limits_{i=1}^9VII_i.
\end{split}
\ee For $VII_i$, $i=$1, 2, 3, using Cauchy inequality, Lemma
\ref{le:H^1 of rho} and Corollary \ref{cor:xt theta xx}, we have
\be\label{VII1-3}\begin{split}
VII_1+VII_2+VII_3\les&(\|u_x\|_{L^\infty}+\|u\|_{L^\infty})\int_Ix^m\rho_{xx}^2+\int_Ix^m\rho_x^2u_{xx}^2+\int_Ix^m\rho_{xx}^2\\
\les&\int_Ix^m\rho_{xx}^2+\|\rho_x\|_{L^\infty}^2\int_Ix^mu_{xx}^2\\
\les&\int_Ix^m\rho_{xx}^2+1.
\end{split}
\ee For $VII_i$, $i=$4, 5, 6, using Cauchy inequality, Lemma
\ref{le:H^1 of rho} and Corollary \ref{cor:xt theta xx} again,
along with Lemma \ref{le:upper bound of density}, we have
\be\label{VII4-6}\begin{split}
VII_4+VII_5+VII_6\les&\int_Ix^m\rho_{xx}^2+\int_Ix^mu_{xxx}^2+\int_Ix^m\rho_{x}^2u_x^2+\int_Ix^m\rho_x^2u^2\\
\les&\int_Ix^m\rho_{xx}^2+\int_Ix^mu_{xxx}^2+(\|u_x\|_{L^\infty}^2+\|u\|_{L^\infty}^2)\int_Ix^m\rho_{x}^2\\
\les&\int_Ix^m\rho_{xx}^2+\int_Ix^mu_{xxx}^2+1.
\end{split}
\ee For $VII_i$, $i=$7, 8, 9, using Lemma \ref{le:upper bound of
density} and Corollary  \ref{cor:xt theta xx}, we have
\be\label{VII7-9}\begin{split}
VII_7+VII_8+VII_9\les&\int_Ix^m\rho_{xx}^2+\int_Ix^m\rho^2u_x^2+\int_Ix^m\rho^2u^2+\int_Ix^mu_{xx}^2\rho^2\\
\les&\int_Ix^m\rho_{xx}^2+1.
\end{split}
\ee Substituting (\ref{VII1-3}), (\ref{VII4-6}) and (\ref{VII7-9})
into (\ref{dt rho xx-1}), we have \be\label{dt rho
xx-2}\begin{split}
\frac{d}{dt}\int_Ix^m\rho_{xx}^2\les&\int_Ix^m\rho_{xx}^2+\int_Ix^mu_{xxx}^2+1.
\end{split}
\ee Differentiating (\ref{non-1.2})$_2$ with respect to $x$, we
have \be\label{u xxx+.}\begin{split} \beta u_{xxx}=&\rho_xu_t+\rho
u_{xt}+\rho_x uu_x+\rho u_x^2+\rho uu_{xx}-\frac{2\rho
vv_x}{x}-\frac{\rho_xv^2}{x}+\frac{\rho
v^2}{x^2}+P_{xx}\\&-\frac{m\beta(xu_{xx}-u_x)}{x^2}+\frac{m\beta(xu_x-2u)}{x^3}.
\end{split}
\ee By (\ref{u xxx+.}), we have \be\label{u xxx}\begin{split}
\int_Ix^mu_{xxx}^2\les&\int_Ix^m\rho_x^2
u_{t}^2+\int_Ix^m\rho^2u_{xt}^2+\int_Ix^m\rho_x^2u^2u_x^2+\int_Ix^m\rho^2u_x^4+\int_Ix^m\rho^2u^2u_{xx}^2\\&+\int_Ix^m\rho^2v^2v_x^2+\int_Ix^m\rho_x^2v^4
+\int_Ix^m\rho^2v^4+\int_Ix^mP_{xx}^2+\int_Ix^m(u_{xx}^2+u_x^2+\frac{u^2}{x^2})\\
\les&\int_Ix^m\rho_x^2u_t^2+\int_Ix^m\rho^2u_{xt}^2+\|u\|_{L^\infty}^2\|u_x\|_{L^\infty}^2\int_Ix^m\rho_{x}^2+\|\rho\|_{L^\infty}^2\|u_x\|_{L^\infty}^2\int_Ix^mu_{x}^2
\\&+\|\rho\|_{L^\infty}^2\|u\|_{L^\infty}^2\int_Ix^mu_{xx}^2+\|\rho\|_{L^\infty}^2\|v\|_{L^\infty}^2\int_Ix^mv_{x}^2+\|v\|_{L^\infty}^4\int_Ix^m\rho_x^2
\\&+\|\rho\|_{L^\infty}\|v\|_{L^\infty}^2\int_Ix^m\rho
v^2+\int_Ix^mP_{xx}^2+1\\
\les&\int_Ix^m\rho_x^2u_t^2+\int_Ix^m\rho^2u_{xt}^2+\int_Ix^mP_{xx}^2+1,
\end{split}
\ee where we have used Lemmas \ref{le:upper bound of density},
\ref{le:H^1 of rho}, and Corollaries \ref{cor:vxx}, \ref{cor:xt
theta xx}. We need to control the right hand side of (\ref{u
xxx}).

For the first one, we have \be\label{1}\begin{split}
\int_Ix^m\rho_x^2u_t^2\les
\|u_t\|_{L^\infty}^2\int_Ix^m\rho_x^2\les\int_Ix^m(u_{xt}^2+\frac{u_t^2}{x^2}).
\end{split}
\ee For the second one, we have \be\label{2}\begin{split}
\int_Ix^m\rho^2u_{xt}^2\les
\|\rho\|_{L^\infty}^2\int_Ix^mu_{xt}^2\les \int_Ix^mu_{xt}^2.
\end{split}
\ee For the third one, recalling $P=\rho Q+P_c(\rho)$, we have
\be\label{3}\begin{split}
\int_Ix^mP_{xx}^2=&\int_Ix^m|\rho_{xx}Q+2\rho_xQ^\prime\theta_x+\rho
Q^{\prime\prime}\theta_x^2+\rho
Q^\prime\theta_{xx}+P_c^{\prime\prime}\rho_x^2+P_c^\prime\rho_{xx}|^2\\
\les&\int_Ix^m\rho_{xx}^2+\int_Ix^m\theta_{xx}^2+1,
\end{split}
\ee where we have used ($A_3$), ($A_4$), Lemmas \ref{le:upper
bound of density}, \ref{le:H^1 of rho}, \ref{le:rho ut}, Corollary
\ref{cor:theta} and Sobolev inequality.

Substituting (\ref{1}), (\ref{2}) and (\ref{3}) into (\ref{u
xxx}), we have \be\label{u xxx-1}\begin{split}
\int_Ix^mu_{xxx}^2\les&\int_Ix^m(u_{xt}^2+\frac{u_t^2}{x^2})+\int_Ix^m\rho_{xx}^2+\int_Ix^m\theta_{xx}^2+1.
\end{split}
\ee Putting (\ref{u xxx-1}) into (\ref{dt rho xx-2}), we have
\be\label{dt rho xx-3}\begin{split}
\frac{d}{dt}\int_Ix^m\rho_{xx}^2\les&\int_Ix^m(u_{xt}^2+\frac{u_t^2}{x^2})+\int_Ix^m\rho_{xx}^2+\int_Ix^m\theta_{xx}^2+1.
\end{split}
\ee By (\ref{dt rho xx-3}), Lemma \ref{le:rho ut}, Corollary
\ref{cor:xt theta xx} and Gronwall inequality, we have
\be\label{rho xx} \int_Ix^m\rho_{xx}^2\le C. \ee  Using Lemma
\ref{le:rho ut} and Corollary \ref{cor:xt theta xx} again, along
with (\ref{u xxx-1}) and (\ref{rho xx}), we get \bex
\int_{Q_T}x^mu_{xxx}^2\le C. \eex
\endpf

\begin{corollary}\label{cor: rho w1,infty}Under the conditions of Theorem
\ref{th:1.1}, we have for any $t\in[0,T]$ \bex
\int_Ix^m\rho_{xt}^2+\|\rho\|_{W^{1,\infty}(Q_T)}+\|\rho_t\|_{L^\infty(Q_T)}+\int_{Q_T}x^m\rho_{tt}^2\le
C. \eex
\end{corollary}
\pf
From (\ref{rhoxt+.}), we have
\be\label{dt rho xt}\begin{split}
\int_Ix^m\rho_{xt}^2\les&\int_Ix^m\rho_{xx}^2u^2+\int_Ix^m\rho_x^2u_x^2+\int_Ix^m\rho^2u_{xx}^2+\int_Ix^m\rho_x^2u^2+\int_Ix^m\rho^2u_x^2+\int_Ix^m\rho^2u^2\\ \les&\|u\|_{L^\infty}^2\int_Ix^m\rho_{xx}^2+\|\rho_x\|_{L^\infty}^2\int_Ix^mu_{x}^2+\|\rho\|_{L^\infty}^2\int_Ix^mu_{xx}^2+\|u\|_{L^\infty}^2\int_Ix^m\rho_{x}^2
\\&+\|\rho\|_{L^\infty}^2\int_Ix^mu_{x}^2+\|\rho\|_{L^\infty}\int_Ix^m\rho u^2\le
C,
\end{split}
\ee where we have used Lemmas \ref{le:upper bound of density},
\ref{le:H^1 of rho}, Corollary \ref{cor:xt theta xx}, and
(\ref{rho xx}).

By Sobolev inequality, Lemmas \ref{le:upper bound of density},
\ref{le:H^1 of rho}, (\ref{rho xx}), (\ref{dt rho xt}), we have
\be\label{rho w1+rho t L}
\|\rho\|_{W^{1,\infty}(Q_T)}+\|\rho_t\|_{L^\infty(Q_T)}\le C. \ee
Differentiating (\ref{non-1.2})$_1$ with respect to $t$, we have
\be\label{rho tt+.}
\rho_{tt}=-(\rho_{xt}u+\rho_xu_t+\rho_tu_x+\rho
u_{xt})-\frac{m(\rho u_t+\rho_tu)}{x}. \ee This deduces
\bex\begin{split}
\int_Ix^m\rho_{tt}^2\les&\int_Ix^m\rho_{xt}^2u^2+\int_Ix^m\rho_{x}^2u_t^2+\int_Ix^m\rho_{t}^2u_x^2+\int_Ix^m\rho^2u_{xt}^2+\int_Ix^m(\rho^2u_t^2+\rho_t^2u^2)\\
\les&\|u\|_{L^\infty}^2\int_Ix^m\rho_{xt}^2+\|u_t\|_{L^\infty}^2\int_Ix^m\rho_{x}^2+\|u_x\|_{L^\infty}^2\int_Ix^m\rho_{t}^2
+\|\rho\|_{L^\infty}^2\int_Ix^mu_{xt}^2\\&+\|\rho\|_{L^\infty}\int_Ix^m\rho
u_t^2+\|u\|_{L^\infty}^2\int_Ix^m\rho_t^2\\
\les&\int_Ix^m(u_{xt}^2+\frac{u_t^2}{x^2})+1,
\end{split}
\eex where we have used (\ref{dt rho xt}), Lemmas \ref{le:upper
bound of density}, \ref{le:H^1 of rho}, \ref{le:rho ut}, and
Corollary \ref{cor:xt theta xx}.

This combining Lemma \ref{le:rho ut} gives \bex
\int_{Q_T}x^m\rho_{tt}^2\le C. \eex
\endpf\\

{\noindent\bf Step 5: $H^2$ estimates of $\theta$}
\begin{lemma}\label{le:dt rho thtea t}Under the conditions of
Theorem \ref{th:1.1}, we have for any $t\in[0,T]$ \bex
\int_Ix^m\rho\theta_t^2+\int_{Q_T}x^m |(\kappa\theta_x)_t|^2\le C.
\eex
\end{lemma}
\pf Differentiating (\ref{non-3.6}) with respect to $t$, we have
\be\label{rho theta tt+.} \begin{split}\rho
Q^\prime\theta_{tt}+\rho
Q^{\prime\prime}\theta_t^2+\rho_tQ^\prime\theta_t+(\rho
uQ^\prime\theta_x)_t+\left(\rho\theta Q^\prime
(u_x+\frac{mu}{x})\right)_t=(\kappa\theta_x)_{xt}+\frac{m(\kappa\theta_x)_t}{x}+\wp_t.\end{split}
\ee Multiplying (\ref{rho theta tt+.}) by
$x^m\left(\int_0^\theta\kappa(\xi)\,d\xi\right)_t$ (i.e.,
$x^m\kappa(\theta)\theta_t$), and integrating by parts over $I$,
we have \be\label{dt rho theta t-1}\begin{split}
&\frac{1}{2}\frac{d}{dt}\int_Ix^m\rho
 Q^\prime\kappa\theta_t^2+\int_Ix^m\left|(\kappa\theta_x)_t\right|^2\\=&-\frac{1}{2}\int_Ix^m\rho_tQ^\prime\kappa\theta_t^2-\frac{1}{2}\int_Ix^m\rho
 Q^{\prime\prime}\theta_t^3\kappa+\frac{1}{2}\int_Ix^m\rho
 Q^\prime\kappa^\prime\theta_t^3-\int_Ix^m(\rho u
 Q^\prime\theta_x)_t\kappa\theta_t \\&-\int_Ix^m\left(\rho\theta Q^\prime
(u_x+\frac{mu}{x})\right)_t\kappa\theta_t+\int_Ix^m\wp_t\kappa\theta_t\\=&\sum\limits_{i=1}^6VIII_i.
\end{split}
\ee For $VIII_1$, using (\ref{non-1.2})$_1$, integration by parts,
($A_4$), ($A_5$), Lemma \ref{le:upper bound of density},
Corollaries \ref{cor:xt L^infty of u}, \ref{cor:theta}, Cauchy
inequality and Poincar$\mathrm{\acute{e}}$ inequality, we have
\be\label{VIII 1}
\begin{split}
VIII_1=&\frac{1}{2}\int_I(x^m\rho u)_xQ^\prime\kappa\theta_t^2\\=&-\frac{1}{2}\int_Ix^m\rho uQ^{\prime\prime}\kappa\theta_x\theta_t^2-\frac{1}{2}\int_Ix^m\rho uQ^\prime\kappa^\prime\theta_x\theta_t^2-\int_Ix^m\rho uQ^\prime\kappa\theta_t\theta_{xt}\\ \les&\|\theta_x\|_{L^\infty}\int_Ix^m\rho\theta_t^2-\int_Ix^m\rho uQ^\prime\theta_t(\kappa\theta_{x})_t+\int_Ix^m\rho uQ^\prime\kappa^\prime\theta_t^2\theta_{x}\\ \le&C\|\theta_x\|_{L^\infty}\int_Ix^m\rho\theta_t^2+\frac{1}{10}\int_Ix^m|(\kappa\theta_{x})_t|^2+C\int_Ix^m\rho\theta_t^2\\ \le&C(1+\|\theta_{xx}\|_{L^2})\int_Ix^m\rho\theta_t^2+\frac{1}{10}\int_Ix^m|(\kappa\theta_{x})_t|^2.
\end{split}
\ee For $VIII_i$, $i=$2, 3, we have
\be\label{VIII2-3}\begin{split}
VIII_2+VIII_3\les&\|\kappa\theta_t\|_{L^\infty}\int_Ix^m\rho
\theta_t^2\\
\les&(\|(\kappa\theta_t)_x\|_{L^2}+\int_Ix^m\rho\kappa|\theta_t|)\int_Ix^m\rho
\theta_t^2\\
\le&\frac{1}{10}\int_Ix^m|(\kappa\theta_{x})_t|^2+\left(\int_Ix^m\rho\theta_t^2\right)^2+C,
\end{split}
\ee where we have used ($A_4$), ($A_5$), Lemmas \ref{non-le:2.1},
\ref{le:energy identity}, Corollary \ref{cor:theta} and Cauchy
inequality.

For $VIII_4$, we have \be\label{VIII4}\begin{split}
VIII_4=&-\int_Ix^m\rho_t u
 Q^\prime\theta_x\kappa\theta_t-\int_Ix^m\rho u_t
 Q^\prime\theta_x\kappa\theta_t-\int_Ix^m\rho u
 Q^{\prime\prime}\theta_x\kappa\theta_t^2-\int_Ix^m\rho u
 Q^\prime\theta_{xt}\kappa\theta_t\\=&\sum\limits_{j=1}^4VIII_{4,j}.
\end{split}
\ee For $VIII_{4,1}$, using (\ref{non-1.2})$_1$, ($A_4$), ($A_5$),
integration by parts, Cauchy inequality, Lemmas \ref{le:upper
bound of density}, \ref{le:rho ut}, and Corollaries
\ref{cor:theta}, \ref{cor:xt theta xx}, we have \be\label{VIII
4,1}\begin{split} VIII_{4,1}=&\int_I(x^m\rho u)_x u
 Q^\prime\theta_x\kappa\theta_t\\=&-\int_Ix^m\rho u u_x
 Q^\prime\theta_x\kappa\theta_t-\int_Ix^m\rho u^2
 Q^{\prime\prime}\theta_x^2\kappa\theta_t-\int_Ix^m\rho u^2
 Q^\prime\theta_{xx}\kappa\theta_t-\int_Ix^m\rho u^2
 Q^\prime\theta_x(\kappa\theta_t)_x\\ \le&C\int_Ix^m\rho \theta_t^2+C\int_Ix^m\theta_x^4+C\int_Ix^m\theta_{xx}^2+\frac{1}{20}\int_Ix^m|(\kappa\theta_x)_t|^2+C
 \\ \le&C\int_Ix^m\rho \theta_t^2+C\|\theta_x\|_{L^\infty}^2\int_Ix^m\theta_x^2+C\int_Ix^m\theta_{xx}^2+\frac{1}{20}\int_Ix^m|(\kappa\theta_x)_t|^2+C\\ \le&C\left(\int_Ix^m\rho \theta_t^2\right)^2+C\int_Ix^m\theta_{xx}^2+\frac{1}{20}\int_Ix^m|(\kappa\theta_x)_t|^2+C.
\end{split}
\ee For $VIII_{4,2}$ and $VIII_{4,3}$, we have \be\label{VIII
4,2-4,3}\begin{split} VIII_{4,2}+VIII_{4,3}
\les&\|\theta_x\|_{L^\infty}^2\int_Ix^m\rho
u_t^2+\int_Ix^m\rho\theta_t^2+\|u\|_{L^\infty}
 \|\theta_x\|_{L^\infty}\int_Ix^m\rho\theta_t^2\\ \les&\|\theta_x\|_{L^2}^2+\|\theta_{xx}\|_{L^2}^2+(1+\|\theta_x\|_{L^2}+\|\theta_{xx}\|_{L^2})\int_Ix^m\rho\theta_t^2
 \\
 \les&\|\theta_{xx}\|_{L^2}^2+\left(\int_Ix^m\rho\theta_t^2\right)^2+1,
\end{split}
\ee where we have used ($A_4$), ($A_5$), Lemma \ref{le:rho ut},
Corollaries \ref{cor:xt L^infty of u}, \ref{cor:theta}, and
Sobolev inequality.

For $VIII_{4,4}$, using ($A_4$), ($A_5$), Lemma \ref{le:rho ut},
Corollaries \ref{cor:xt L^infty of u}, \ref{cor:theta}, and
Sobolev inequality again, along with  Lemma \ref{le:upper bound of
density}, we have \be\label{VIII 4,4}
\begin{split}
VIII_{4,4}=&-\int_Ix^m\rho u
 Q^\prime(\kappa\theta_{x})_t\theta_t+\int_Ix^m\rho u
 Q^\prime\kappa^\prime\theta_{x}\theta_t^2\\ \le&\frac{1}{20}\int_Ix^m|(\kappa\theta_x)_t|^2+C\left(\|\rho\|_{L^\infty}\|u\|_{L^\infty}^2
+\|u\|_{L^\infty}\|\theta_x\|_{L^\infty}\right)\int_Ix^m\rho\theta_t^2\\
\le&\frac{1}{20}\int_Ix^m|(\kappa\theta_x)_t|^2+C\int_Ix^m\theta_{xx}^2
+C\left(\int_Ix^m\rho\theta_t^2\right)^2+C.
\end{split}
\ee Putting (\ref{VIII 4,1}), (\ref{VIII 4,2-4,3}) and (\ref{VIII
4,4}) into (\ref{VIII4}), we have \be\label{VIII 4+1}\begin{split}
VIII_4\le
\frac{1}{10}\int_Ix^m|(\kappa\theta_x)_t|^2+C\left(\int_Ix^m\rho
\theta_t^2\right)^2 +C\int_Ix^m\theta_{xx}^2+C.
\end{split}
\ee
For $VIII_5$, we have
\be\label{VIII 5}\begin{split}
VIII_5=&-\int_Ix^m\rho_t\theta Q^\prime
(u_x+\frac{mu}{x})\kappa\theta_t-\int_Ix^m\rho Q^\prime
(u_x+\frac{mu}{x})\kappa\theta_t^2\\&-\int_Ix^m\rho\theta Q^{\prime\prime}
(u_x+\frac{mu}{x})\kappa\theta_t^2-\int_Ix^m\rho\theta Q^\prime
(u_{xt}+\frac{mu_t}{x})\kappa\theta_t\\=&\sum\limits_{j=1}^4VIII_{5,j}.
\end{split}
\ee For $VIII_{5,1}$, we have \be\label{VIII 5,1}\begin{split}
VIII_{5,1}=&\int_I(x^m\rho u)_x\theta Q^\prime
(u_x+\frac{mu}{x})\kappa\theta_t\\=&-\int_Ix^m\rho u\theta_x
Q^\prime (u_x+\frac{mu}{x})\kappa\theta_t-\int_Ix^m\rho
u\theta\theta_x Q^{\prime\prime}
(u_x+\frac{mu}{x})\kappa\theta_t\\&-\int_Ix^m\rho u\theta Q^\prime
(u_{xx}+\frac{mu_x}{x}-\frac{mu}{x^2})\kappa\theta_t-\int_Ix^m\rho
u\theta Q^\prime (u_x+\frac{mu}{x})(\kappa\theta_t)_x\\
\le&C\int_Ix^m\rho\theta_t^2+C\|u\|_{W^{1,\infty}}^2\int_Ix^m\theta_x^2+C\int_Ix^m(u_{xx}^2+u_x^2
+\frac{u^2}{x^2})+\frac{1}{10}\int_Ix^m|(\kappa\theta_x)_t|^2\\
\le&\frac{1}{10}\int_Ix^m|(\kappa\theta_x)_t|^2+C\int_Ix^m\rho\theta_t^2+C,
\end{split}
\ee where we have used (\ref{non-1.2})$_1$, integration by parts,
Cauchy inequality, ($A_4$), ($A_5$), Lemmas \ref{le:upper bound of
density}, \ref{le:rho ut}, and Corollaries \ref{cor:theta},
\ref{cor:xt theta xx}.

For $VIII_{5,2}$ and $VIII_{5,3}$, we have
\be\label{VIII 5,2-5,3}\begin{split}
VIII_{5,2}+VIII_{5,3}\le&C\int_Ix^m\rho\theta_t^2.
\end{split}
\ee
For $VIII_{5,4}$, we have
\be\label{VIII 5,4}
\begin{split}
VIII_{5,4}\le&C\int_Ix^m\rho\theta_t^2+C\int_Ix^m(u_{xt}^2+\frac{u_t^2}{x^2}).
\end{split}
\ee
Putting (\ref{VIII 5,1}), (\ref{VIII 5,2-5,3}) and (\ref{VIII 5,4}) into (\ref{VIII 5}), we have
\be\label{VIII 5+1}\begin{split}
VIII_5\le\frac{1}{10}\int_Ix^m|(\kappa\theta_x)_t|^2+C\int_Ix^m\rho\theta_t^2+C\int_Ix^m(u_{xt}^2+\frac{u_t^2}{x^2})+C.
\end{split}
\ee For $VIII_6$, recalling
$\wp=\lambda(u_x+\frac{mu}{x})^2+\mu\left(w_x^2+2u_x^2+(v_x-\frac{mv}{x})^2+\frac{2mu^2}{x^2}\right)$,
we have \be\label{VIII 6}\begin{split}
VIII_6\les&\|\kappa\theta_t\|_{L^\infty}\int_Ix^m|\wp_t|\\
\les&\|\kappa\theta_t\|_{L^\infty}\left((\int_Ix^mu_{xt}^2)^\frac{1}{2}
+(\int_Ix^m\frac{u_t^2}{x^2})^\frac{1}{2}+(\int_Ix^mw_{xt}^2)^\frac{1}{2}
+(\int_Ix^m(v_{xt}^2+\frac{v_t^2}{x^2}))^\frac{1}{2}\right)\\ \le
&\frac{1}{10}\int_Ix^m|(\kappa\theta_x)_t|^2+C\int_Ix^m\rho\theta_t^2+C\int_Ix^mu_{xt}^2
+C\int_Ix^m\frac{u_t^2}{x^2}+C\int_Ix^mw_{xt}^2
\\&+C\int_Ix^m(v_{xt}^2+\frac{v_t^2}{x^2}),
\end{split}
\ee where we have used Lemma \ref{non-le:2.1} and Cauchy
inequality.

Putting (\ref{VIII 1}), (\ref{VIII2-3}), (\ref{VIII 4+1}), (\ref{VIII 5+1}) and (\ref{VIII 6}) into (\ref{dt rho theta t-1}), we have
\be\label{dt rho theta t-2}\begin{split}
&\frac{d}{dt}\int_Ix^m\rho
 Q^\prime\kappa\theta_t^2+\int_Ix^m\left|(\kappa\theta_x)_t\right|^2\\ \le & C\left(\int_Ix^m\rho \theta_t^2\right)^2
+C\int_Ix^m\theta_{xx}^2+C\int_Ix^m(u_{xt}^2+\frac{u_t^2}{x^2}+w_{xt}^2+v_{xt}^2+\frac{v_t^2}{x^2})+C.
\end{split}
\ee By Gronwall inequality, we complete the proof of Lemma
\ref{le:dt rho thtea t}.
\endpf
\begin{corollary}\label{cor:t theta t}Under the conditions of
Theorem \ref{th:1.1}, we have  \bex
\int_0^T\|\theta_t\|_{L^\infty}^2\le C. \eex
\end{corollary}
\pf By Lemma \ref{non-le:2.1} and ($A_5$), we get
\bex\begin{split}
\int_0^T\|\theta_t\|_{L^\infty}^2\les&\int_0^T\|\kappa\theta_t\|_{L^\infty}^2\\
\les&\int_0^T\left(\int_I\rho\kappa|\theta_t|\right)^2+\int_{Q_T}|(\kappa\theta_t)_x|^2\\
\les&
\int_0^T\int_I\rho\theta_t^2+\int_{Q_T}|(\kappa\theta_x)_t|^2\le
C.
\end{split}
\eex
\endpf

\begin{corollary}\label{cor:xt theta xt}Under the conditions of
Theorem \ref{th:1.1}, we have  \bex \int_{Q_T}x^m\theta_{xt}^2\le
C. \eex
\end{corollary}
\pf
Since
$$
\kappa\theta_{xt}=(\kappa\theta_x)_t-\kappa^\prime\theta_t\theta_x,
$$
we obtain \bex\begin{split}
\int_{Q_T}x^m\theta_{xt}^2\les&\int_{Q_T}x^m\kappa^2\theta_{xt}^2\\ \les&\int_{Q_T}x^m\left|(\kappa\theta_x)_t\right|^2+\int_{Q_T}x^m(\kappa^\prime)^2\theta_t^2\theta_x^2
\\ \les&1+\int_0^T\sup\limits_{x\in
I}\theta_t^2\int_Ix^m\theta_x^2 \le C, \end{split}\eex where we
have used ($A_5$), Lemmas \ref{le:rho ut}, \ref{le:dt rho thtea
t}, and Corollaries \ref{cor:theta}, \ref{cor:t theta t}.

\endpf

\begin{corollary}\label{cor:theta w 1,infty} Under the conditions
of Theorem \ref{th:1.1}, we have for any $t\in[0,T]$ \bex
\|\theta\|_{W^{1,\infty}(Q_T)}+\int_Ix^m
\theta_{xx}^2+\int_{Q_T}x^m\theta_{xxx}^2\le C. \eex
\end{corollary}
\pf
By (\ref{int I theta xx1}) and Lemma \ref{le:dt rho thtea t}, we get
$$
\int_Ix^m\theta_{xx}^2\le C.
$$
This, along with Lemma \ref{le:rho ut}, Corollary \ref{cor:theta}
and Sobolev inequality, gives
$$
\|\theta\|_{W^{1,\infty}(Q_T)}\le C.
$$
Differentiating (\ref{non-3.6}) w.r.t. $x$, we have
\be\label{theta xxx+.}\begin{split}
&\kappa\theta_{xxx}=-(\frac{m\kappa\theta_x}{x})_x-3\kappa^\prime\theta_x\theta_{xx}-\kappa^{\prime\prime}\theta_x^3-\wp_x+\rho
Q^\prime\theta_{xt}+\rho_xQ^\prime\theta_t+\rho
Q^{\prime\prime}\theta_x\theta_t\\&\ \ \ \ \ \ \ \ \ \ \ \ +(\rho u
Q^\prime\theta_x)_x+\left(\rho\theta Q^\prime (u_x+\frac{mu}{x})\right)_x. \end{split}\ee
This deduces
\be\label{int I theta xxx 1}\begin{split}
\int_Ix^m\theta_{xxx}^2\les&\int_Ix^m\rho^2\theta_{xt}^2+\int_Ix^m\rho_x^2\theta_t^2+1\\ \les&\int_Ix^m\theta_{xt}^2+\sup\limits_{x\in I}\theta_t^2+1.
\end{split}
\ee
Integrating (\ref{int I theta xxx 1}) over $(0,T)$, and using Corollary \ref{cor:t theta t} and Corollary \ref{cor:xt theta xt}, we get
$$
\int_{Q_T}x^m\theta_{xxx}^2\le C.
$$
\endpf

{\noindent\bf Constructing a sequence of approximate initial data
and passing to the limits:}

Denote $\rho_0^\varepsilon=\rho_0+\varepsilon$ and
$\theta_0^\varepsilon=\theta_0+\varepsilon$. $u_0^\varepsilon$
 is the solution to the
 following elliptic problem:
 \beq\label{Auxiliary problem 1}
\begin{cases}\beta(u_{0xx}^\varepsilon+\frac{mu_{0x}^\varepsilon}{x}-\frac{mu_0^\varepsilon}{x^2})-P_x(\rho_0^\varepsilon,\theta_0^\varepsilon)=\sqrt{\rho}_0g_1,\ x\in I,\\
u_0^\varepsilon|_{x=a,b}=0.
\end{cases}\eeq
Since $\rho_0^\varepsilon=\rho_0+\varepsilon\in H^2$,
$\theta_0^\varepsilon=\theta_0+\varepsilon\in H^2$ and $g_1\in
L^2$, we obtain by using the standard elliptic estimates
\beq\label{u_0 H2 converges}
\begin{cases}
u_0^\varepsilon\in H^2\cap H_0^1,\\
\|u_0^\varepsilon-u_0\|_{H^2}\le C\varepsilon.
\end{cases}
\eeq Consider (\ref{non-1.2})-(\ref{non-1.4}) with initial data
replaced by ($\rho_0^\varepsilon$, $u_0^\varepsilon$, $v_0$,
$w_0$, $\theta_0^\varepsilon$), we get a global solution
($\rho^\varepsilon$, $u^\varepsilon$, $v^\varepsilon$,
$w^\varepsilon$, $\theta^\varepsilon$) to
(\ref{non-1.2})-(\ref{non-1.4}) for each $\varepsilon>0$. These
{\it a priori} estimates in the section are valid for
($\rho^\varepsilon$, $u^\varepsilon$, $v^\varepsilon$,
$w^\varepsilon$, $\theta^\varepsilon$). Then taking the limits
$\varepsilon\rightarrow0^+$ (taking subsequence if necessary), and
using Lemma \ref{non-le:2.3}, we get a solution denoted by
($\rho$, $u$, $v$, $w$, $\theta$) to
(\ref{non-1.2})-(\ref{non-1.4}) with regularities as in Theorem
\ref{th:1.1}. The uniqueness of the solutions can be done by the
standard methods, see for instance \cite{cho-Kim: perfect gas} and
references therein. We end the proof of Theorem \ref{th:1.1}. \qed
\begin{remark}
The global existence of ($\rho^\varepsilon$, $u^\varepsilon$,
$v^\varepsilon$, $w^\varepsilon$, $\theta^\varepsilon$) can be
done by {\it local existence} (using the similar arguments as in
\cite{cho-Kim: perfect gas})+ {\it some a priori estimates
globally in time throughout the section}.
\end{remark}

\section{Proof of Theorem \ref{th:1.2}}
\setcounter{equation}{0} \setcounter{theorem}{0}
In the section, we
denote by $C$ a generic constant depending only on $\|(\rho_0, u_0,
v_0, w_0, \theta_0)\|_{H^3}$, $\|g_i\|_{L^2}$,
$\|(\sqrt{\rho_0}g_i)_x\|_{L^2}$ ($i=1,2,3,4$), $T$, $\lambda$,
$\mu$, $a$, $b$, and some other known constants, but independent of
the solutions and the lower bounds of the density.

Under the conditions of Theorem \ref{th:1.2}, all those {\it a
priori} estimates in section 3 are also satisfied in this section.
Following the quite similar strategies with the proof of Theorem
\ref{th:1.1}, to prove Theorem \ref{th:1.2}, it suffices to get
$H^3$-{\it a-priori}-estimates of the classical solutions ($\rho$,
u, v, w, $\theta$) as in Theorem \ref{th:1.2} with
$\inf\limits_{(x,t)\in Q_T}\rho>0$.

To begin with, we get a $W^{1,\infty}(Q_T)$-estimate of
$\sqrt{\rho}$ which was obtained in \cite{Ding-Wen-Zhu,
Ding-Wen-Yao-Zhu} as a crucial estimate to get $H^4$-estimate of
$u$.
\begin{lemma}\label{le:rho x12infty} Under the conditions
of Theorem \ref{th:1.2}, we have
$$
\|(\sqrt{\rho})_x\|_{L^\infty(Q_T)}+\|(\sqrt{\rho})_t\|_{L^\infty(Q_T)}\le
C.
$$
\end{lemma}
\pf Multiplying $(\ref{non-1.2})_1$ by
$\frac{1}{2\sqrt{\rho}}$, we have \be\label{rho12t+.}
(\sqrt{\rho})_t+\frac{mx^{-1}}{2}\sqrt{\rho}u+\frac{1}{2}\sqrt{\rho}u_x+(\sqrt{\rho})_xu=0.
\ee Differentiating (\ref{rho12t+.}) with respect to $x$, we get
$$
(\sqrt{\rho})_{xt}+\frac{mx^{-1}}{2}(\sqrt{\rho})_xu+\frac{mx^{-1}}{2}\sqrt{\rho}u_x-\frac{m\sqrt{\rho}u}{2x^2}+\frac{3}{2}(\sqrt{\rho})_xu_x+
\frac{1}{2}\sqrt{\rho}u_{xx}+(\sqrt{\rho})_{xx}u=0.
$$
Denote $h=(\sqrt{\rho})_x$, we have
$$
h_t+h_xu+h(\frac{mx^{-1}u}{2}+\frac{3u_x}{2})+\frac{mx^{-1}}{2}\sqrt{\rho}u_x-\frac{m\sqrt{\rho}u}{2x^2}+
\frac{1}{2}\sqrt{\rho}u_{xx}=0,
$$
which implies \be\label{dt h}\begin{split}
&\frac{d}{dt}\left\{h\exp\left[\int_0^t(\frac{mx^{-1}u}{2}+
\frac{3u_x}{2})(x(\tau,y),\tau)d\tau\right]\right\}\\=&-\left(\frac{mx^{-1}\sqrt{\rho}u_x}{2}-\frac{m\sqrt{\rho}u}{2x^2}+
\frac{\sqrt{\rho}u_{xx}}{2}\right)\\
&\ \ \ \ \ \times \exp\left[\int_0^t\left(\frac{mx^{-1}u}{2}+
\frac{3u_x}{2}\right)(x(\tau,y),\tau)d\tau\right], \end{split}\ee where
$x(t,y)$ satisfies
 $$
\begin{cases}
\frac{d x(t,y)}{d t}=u(x(t,y),t),\ 0\le t<s,\\
x(s,y)=y.
\end{cases} $$
Integrating (\ref{dt h}) over $(0,s)$, we get
$$
\arraycolsep=1.5pt
\begin{array}{rl}
h(y,s)= & \displaystyle \exp\left(-\int_0^s\left(\frac{mx^{-1}u}{2}+
\frac{3u_x}{2}\right)(x(\tau,y),\tau)d\tau\right)h(x(0,y),0)
\\ [5mm]
& \displaystyle
-\int_0^s\Big[(\frac{mx^{-1}\sqrt{\rho}u_x}{2}-\frac{m\sqrt{\rho}u}{2x^2}+
\frac{\sqrt{\rho}u_{xx}}{2})\\
&\ \ \ \  \ \ \ \ \ \
\times\exp\left(\int_s^t\left(\frac{mx^{-1}u}{2}+
\frac{3u_x}{2}\right)(x(\tau,y),\tau)d\tau\right)\Big]dt.
\end{array}
$$
This implies \be\label{rho12 x}
\|(\sqrt{\rho})_x\|_{L^\infty(Q_T)}\le C. \ee From (\ref{rho12t+.}) and
(\ref{rho12 x}),
we get
$$
\|(\sqrt{\rho})_t\|_{L^\infty(Q_T)}\le C.
$$
The proof of Lemma \ref{le:rho x12infty} is complete. \endpf
\begin{lemma}\label{le:rho theta xt}Under the conditions
of Theorem \ref{th:1.2}, we have for any $t\in[0,T]$
$$
\int_Ix^m\rho^2\left|(\kappa\theta_x)_t\right|^2+\int_{Q_T}x^m\rho^3\theta_{tt}^2\le
C.
$$
\end{lemma}
\pf Multiplying $(\ref{rho theta tt+.})$ by
$x^m\rho^\alpha(\kappa\theta_t)_t$ (i.e.
$x^m\rho^{\alpha}\kappa\theta_{tt}+x^m\rho^{\alpha}\kappa^\prime\theta_t^2$,
where $\alpha>0$ is to be decided later), and integrating by parts
over $I$, we have
 \be\label{dt rho theta xt-1}\begin{split}
&\int_Ix^m\rho^{\alpha+1}\kappa
Q^\prime\theta_{tt}^2+\frac{1}{2}\frac{d}{dt}\int_Ix^m\rho^{\alpha}\left|(\kappa\theta_x)_t\right|^2\\=&
\frac{\alpha}{2}\int_Ix^m\rho^{\alpha-1}\rho_t\left|(\kappa\theta_x)_t\right|^2
-\alpha\int_Ix^m\rho^{\alpha-1}\rho_x\kappa\theta_{tt}(\kappa\theta_x)_t-\alpha\int_Ix^m\rho^{\alpha-1}\rho_x\kappa^\prime\theta_t^2(\kappa\theta_x)_t
\\&+\int_Ix^m\rho^{\alpha}\wp_t(\kappa\theta_{tt}+\kappa^\prime\theta_t^2)-\int_Ix^m\rho^{\alpha+1}Q^\prime\kappa^\prime\theta_t^2\theta_{tt}
-\int_Ix^m\rho^{\alpha}\left[\rho
Q^{\prime\prime}\theta_t^2+\rho_tQ^\prime\theta_t\right]\left(\kappa\theta_{tt}+\kappa^\prime\theta_t^2\right)\\&-\int_Ix^m\rho^{\alpha}(\rho
uQ^\prime\theta_x)_t\left(\kappa\theta_{tt}+\kappa^\prime\theta_t^2\right)-\int_Ix^m\rho^{\alpha}\left(\rho\theta
Q^\prime
(u_x+\frac{mu}{x})\right)_t\left(\kappa\theta_{tt}+\kappa^\prime\theta_t^2\right)=\sum\limits_{i=1}^8VIV_i.
\end{split}
\ee For $VIV_2$, using ($A_4$), ($A_5$), Corollary \ref{cor:theta}
and Cauchy inequality, we have \be\label{VIV 2}
\begin{split}VIV_2
\le&\frac{1}{10}\int_Ix^m\rho^{\alpha+1}\kappa
Q^\prime\theta_{tt}^2+C\int_Ix^m\rho^{\alpha-3}\rho_x^2|(\kappa\theta_x)_t|^2\\
\le&\frac{1}{10}\int_Ix^m\rho^{\alpha+1}\kappa
Q^\prime\theta_{tt}^2+C\|\rho^{\alpha-3}\rho_x^2\|_{L^\infty}\int_Ix^m|(\kappa\theta_x)_t|^2.
\end{split}\ee
For $VIV_1$ and $VIV_3$, we have \be\label{VIV 1 and
3}\begin{split}
VIV_1+VIV_3\le&C\|\rho^{\alpha-1}\rho_t\|_{L^\infty}\int_Ix^m\left|(\kappa\theta_x)_t\right|^2
+C\|\theta_t\|_{L^\infty}^2\int_Ix^m\rho^\alpha|(\kappa\theta_x)_t|^2\\&+C\|\theta_t\|_{L^\infty}^2\int_Ix^m\rho^{\alpha-2}\rho_x^2.
\end{split}
\ee For $VIV_4$, we have \be\label{VIV 4}\begin{split}
VIV_4\le&\frac{1}{10}\int_Ix^m\rho^{\alpha+1}\kappa
Q^\prime\theta_{tt}^2+C\|\rho^{\alpha-1}\|_{L^\infty}\int_Ix^m\wp_t^2+C\int_Ix^m\rho^{\alpha+1}\theta_t^4\\
\le&\frac{1}{10}\int_Ix^m\rho^{\alpha+1}\kappa
Q^\prime\theta_{tt}^2+C\|\rho^{\alpha-1}\|_{L^\infty}\int_Ix^m(u_{xt}^2+\frac{u_t^2}{x^2}+w_{xt}^2+v_{xt}^2+\frac{v_t^2}{x^2})+C\|\theta_t\|_{L^\infty}^2,
\end{split}
\ee where we have used Cauchy inequality, Corollaries
\ref{cor:vxx}, \ref{cor:wxx}, \ref{cor:xt theta xx} and Lemmas
\ref{le:upper bound of density}, \ref{le:dt rho thtea t}.

For $VIV_5$ and $VIV_6$, using ($A_4$), ($A_5$), Cauchy
inequality, Corollary \ref{cor:theta}, and Lemmas \ref{le:upper
bound of density}, \ref{le:dt rho thtea t}, we have \be\label{VIV
5 and 6}\begin{split}
VIV_5+VIV_6=&-\int_Ix^m\rho^{\alpha+1}Q^\prime\kappa^\prime\theta_t^2\theta_{tt}
-\int_Ix^m\rho^\alpha\left[\rho
Q^{\prime\prime}\theta_t^2+\rho_tQ^\prime\theta_t\right]\kappa\theta_{tt}\\&-\int_Ix^m\rho^\alpha\left[\rho
Q^{\prime\prime}\theta_t^2+\rho_tQ^\prime\theta_t\right]\kappa^\prime\theta_t^2\\
\le&\frac{1}{10}\int_Ix^m\rho^{\alpha+1}\kappa
Q^\prime\theta_{tt}^2+C\int_Ix^m\rho^{\alpha+1}\theta_t^4+C\int_Ix^m\rho^{\alpha-1}\rho_t^2\theta_t^2\\
\le&\frac{1}{10}\int_Ix^m\rho^{\alpha+1}\kappa
Q^\prime\theta_{tt}^2+C\|\theta_t\|_{L^\infty}^2(1+\int_Ix^m\rho^{\alpha-1}\rho_t^2)+C.
\end{split}
\ee For $VIV_7$, using ($A_4$), ($A_5$), Cauchy inequality, and
Lemmas \ref{le:upper bound of density}, \ref{le:dt rho thtea t}
again, along with Corollaries \ref{cor:xt L^infty of u},
\ref{cor:theta w 1,infty} and Lemma \ref{le:rho ut}, we have
\be\label{VIV 7}\begin{split} VIV_7=&-\int_Ix^m\rho^\alpha\rho_t
uQ^\prime\theta_x\left(\kappa\theta_{tt}+\kappa^\prime\theta_t^2\right)-\int_Ix^m\rho^{\alpha+1}
u_t
Q^\prime\theta_x\left(\kappa\theta_{tt}+\kappa^\prime\theta_t^2\right)\\&-\int_Ix^m\rho^{\alpha+1}
uQ^{\prime\prime}\theta_t\theta_x\left(\kappa\theta_{tt}+\kappa^\prime\theta_t^2\right)-\int_Ix^m\rho^{\alpha+1}
uQ^\prime\theta_{xt}\left(\kappa\theta_{tt}+\kappa^\prime\theta_t^2\right)\\
\le&\frac{1}{10}\int_Ix^m\rho^{\alpha+1}\kappa
Q^\prime\theta_{tt}^2+C\int_Ix^m\rho^{\alpha-1}\rho_t^2+C\int_Ix^m\rho^{\alpha+1}\theta_t^4+C\int_Ix^m\rho^{\alpha+1}
u_t^2\\&+C\int_Ix^m\rho^{\alpha+1}
\theta_t^2+C\int_Ix^m\rho^{\alpha+1}\theta_{xt}^2\\
\le&\frac{1}{10}\int_Ix^m\rho^{\alpha+1}\kappa
Q^\prime\theta_{tt}^2+C\int_Ix^m\rho^{\alpha-1}\rho_t^2+C\|\theta_t\|_{L^\infty}^2+C\int_Ix^m\theta_{xt}^2+C.
\end{split}
\ee Similarly, for $VIV_8$, we have \be\label{VIV 8}\begin{split}
VIV_8=&-\int_Ix^m\rho^\alpha\rho_t\theta Q^\prime
(u_x+\frac{mu}{x})\left(\kappa\theta_{tt}+\kappa^\prime\theta_t^2\right)-\int_Ix^m\rho^{\alpha+1}\theta_t
Q^\prime
(u_x+\frac{mu}{x})\left(\kappa\theta_{tt}+\kappa^\prime\theta_t^2\right)\\&-\int_Ix^m\rho^{\alpha+1}\theta
Q^{\prime\prime}
\theta_t(u_x+\frac{mu}{x})\left(\kappa\theta_{tt}+\kappa^\prime\theta_t^2\right)-\int_Ix^m\rho^{\alpha+1}\theta
Q^\prime
(u_{xt}+\frac{mu_t}{x})\left(\kappa\theta_{tt}+\kappa^\prime\theta_t^2\right)\\
\le&\frac{1}{10}\int_Ix^m\rho^{\alpha+1}\kappa
Q^\prime\theta_{tt}^2+C\int_Ix^m\rho^{\alpha-1}\rho_t^2+C\int_Ix^m\rho^{\alpha+1}\theta_t^4\\&+C\int_Ix^m\rho^{\alpha+1}\theta_t^2
+C\int_Ix^m\rho^{\alpha+1}(u_{xt}^2+\frac{u_t^2}{x^2})\\
\le&\frac{1}{10}\int_Ix^m\rho^{\alpha+1}\kappa
Q^\prime\theta_{tt}^2+C\int_Ix^m\rho^{\alpha-1}\rho_t^2+C\|\theta_t\|_{L^\infty}^2+C\int_Ix^m(u_{xt}^2+\frac{u_t^2}{x^2})+C.
\end{split}
\ee Putting (\ref{VIV 2}), (\ref{VIV 1 and 3}), (\ref{VIV 4}),
(\ref{VIV 5 and 6}), (\ref{VIV 7}) and (\ref{VIV 8}) into (\ref{dt
rho theta xt-1}), we have \be\label{dt rho theta
xt-2}\begin{split} &\int_Ix^m\rho^{\alpha+1}\kappa
Q^\prime\theta_{tt}^2+\frac{d}{dt}\int_Ix^m\rho^\alpha\left|(\kappa\theta_x)_t\right|^2\\
\le&C\|\theta_t\|_{L^\infty}^2\int_Ix^m\rho^\alpha|(\kappa\theta_x)_t|^2+C(\|\rho^{\alpha-3}\rho_x^2\|_{L^\infty}+\|\rho^{\alpha-1}\rho_t\|_{L^\infty})\int_Ix^m\left|(\kappa\theta_x)_t\right|^2
\\&+C(1+\int_Ix^m\rho^{\alpha-1}\rho_t^2+\int_Ix^m\rho^{\alpha-2}\rho_x^2)\|\theta_t\|_{L^\infty}^2
+C\int_Ix^m\rho^{\alpha-1}\rho_t^2+C\int_Ix^m\theta_{xt}^2
\\&+C\|\rho^{\alpha-1}\|_{L^\infty}\int_Ix^m(u_{xt}^2+\frac{u_t^2}{x^2}+w_{xt}^2+v_{xt}^2+\frac{v_t^2}{x^2})
+C\int_Ix^m(u_{xt}^2+\frac{u_t^2}{x^2})+C.
\end{split}
\ee {\noindent\bf Claim:}
\be\label{claim}\int_Ix^m\rho^\alpha\left|(\kappa\theta_x)_t\right|^2+\int_{Q_T}x^m\rho^{\alpha+1}\kappa
Q^\prime\theta_{tt}^2\le C, \ee if \be\label{claim 1}
\|\rho^{\alpha-3}\rho_x^2\|_{L^\infty}+\|\rho^{\alpha-1}\rho_t\|_{L^\infty}+\|\rho^{\alpha-1}\|_{L^\infty}\le
C. \ee In fact, if (\ref{claim 1}) is satisfied, (\ref{claim})
will be obtained by using (\ref{dt rho theta xt-2}), Corollaries
\ref{cor:t theta t}, \ref{cor:xt theta xt}, Lemmas \ref{le:rho
vt2}, \ref{le:rho wt}, \ref{le:rho ut}, \ref{le:dt rho thtea t},
(\ref{compatibility conditions}) and Gronwall inequality.

To ensure that (\ref{claim 1}) is valid, the restriction
$\alpha\ge3$ seems necessary. While, to get
${L^\infty_t}H^3_x$-estimate of $\theta$ (Corollary \ref{cor:theta
xxx}), by (\ref{non-3.6}), $\alpha$ in (\ref{claim}) should
satisfy $\alpha\le 2$ . Assume $\alpha\in[1,2]$, to get
(\ref{claim 1}), it suffices to prove
$$
\|\rho^{\alpha-3}\rho_x^2\|_{L^\infty}\le C.
$$
In fact, \bex\begin{split}
\|\rho^{\alpha-3}\rho_x^2\|_{L^\infty}=4\left\|\rho^{\alpha-2}|(\sqrt{\rho})_x|^2\right\|_{L^\infty}.
\end{split}
\eex Assume $\alpha=2$, and use Lemma \ref{le:rho x12infty},
(\ref{claim 1}) is to be obtained. This means that (\ref{claim})
is obtained for $\alpha=2$. The proof of Lemma \ref{le:rho theta
xt} is complete.
\endpf

\begin{corollary}\label{cor:theta xxx}Under the conditions
of Theorem \ref{th:1.2}, we have for any $t\in[0,T]$
$$
\int_Ix^m\left(\theta_{xxx}^2+\rho^2\theta_{xt}^2\right)\le C.
$$
\end{corollary}
\pf A direct calculation gives
$$
\rho\kappa\theta_{xt}=\rho(\kappa\theta_x)_t-\rho\kappa^\prime\theta_t\theta_x,
$$
which implies \be\begin{split}\label{rho thetaxt}
\int_Ix^m\rho^2\theta_{xt}^2\le&C\int_Ix^m\rho^2\left|(\kappa\theta_x)_t\right|^2+C\|\theta_x\|_{L^\infty}^2\int_Ix^m\rho\theta_t^2 \\
 \le & C, \end{split}\ee
where we have used $(A_5)$, Lemma \ref{le:upper bound of density},
Lemma \ref{le:dt rho thtea t}, Lemma \ref{le:rho theta xt} and
Corollary \ref{cor:theta w 1,infty}.

 From the first inequality of (\ref{int I theta xxx 1}), using Lemma \ref{le:dt rho thtea
t}, (\ref{rho thetaxt}) and Lemma \ref{le:rho x12infty}, we obtain
\bex\begin{split}\int_Ix^m\theta_{xxx}^2\le&C\int_Ix^m\rho^2\theta_{xt}^2+C\int_Ix^m\rho_x^2\theta_t^2+C\\
\le&C\int_Ix^m\rho^2\theta_{xt}^2+C\int_Ix^m\rho|(\sqrt{\rho})_x|^2\theta_t^2+C\le
C.\end{split}\eex \endpf

\begin{lemma}\label{le:rho uxt}Under the conditions
of Theorem \ref{th:1.2}, we have for any $t\in[0,T]$
$$
\int_Ix^m\rho^2u_{xt}^2+\int_{Q_T}x^m\rho^3u_{tt}^2\le C.
$$
\end{lemma}
\pf
 Similar to Lemma \ref{le:rho theta xt},
multiplying (\ref{rho utt+.=}) by $x^m\rho^2u_{tt}$, and integrating by parts
over $I$, we have \be\label{dt rho uxt-1}
\begin{split}&\int_Ix^m\rho^3u_{tt}^2+\frac{\beta}{2}\frac{d}{dt}\int_Ix^m\rho^2
u_{xt}^2\\=&\beta\int_Ix^m\rho\rho_tu_{xt}^2-2\beta\int_Ix^m\rho\rho_xu_{xt}u_{tt}-\int_Ix^m\rho^2u_{tt}(\rho_tu_t+\rho_tuu_x+\rho
u_tu_x+\rho uu_{xt}-\frac{\rho_tv^2}{x}-\frac{2\rho
vv_t}{x})\\&-\beta\int_Ix^m\rho^2u_{tt}\frac{mu_t}{x^2}-\int_Ix^m\rho^2u_{tt}P_{xt}\\
=&\sum\limits_{i=1}^5VV_i.\end{split} \ee For $VV_1$ and $VV_2$,
using Corollary \ref{cor: rho w1,infty}, Lemma \ref{le:rho
x12infty} and Cauchy inequality, we have \be\label{VV 1 and
2}\begin{split} VV_1+VV_2\le
&\|\rho\|_{L^\infty}\|\rho_t\|_{L^\infty}\int_Ix^mu_{xt}^2-4\beta\int_Ix^m\rho^\frac{3}{2}(\sqrt{\rho})_xu_{xt}u_{tt}\\
\le&C\int_Ix^mu_{xt}^2+\frac{1}{8}\int_Ix^m\rho^3u_{tt}^2.
\end{split}
\ee For $VV_3$, we have \be\label{VV 3}\begin{split} VV_3
\le&\frac{1}{8}\int_Ix^m\rho^3u_{tt}^2+C\int_Ix^m\rho
\rho_t^2u_t^2+C\int_Ix^m\rho\rho_t^2u^2u_x^2+C\int_Ix^m\rho^3u_t^2u_x^2\\&+C\int_Ix^m\rho^3u^2u_{xt}^2+C\int_Ix^m\frac{\rho\rho_t^2v^4}{x^2}
+C\int_Ix^m\frac{\rho^3v^2v_t^2}{x^2}\\
\le&\frac{1}{8}\int_Ix^m\rho^3u_{tt}^2+C\|\rho_t\|_{L^\infty}^2\int_Ix^m\rho
u_t^2+C\|\rho\|_{L^\infty}\|\rho_t\|_{L^\infty}^2\|u\|_{L^\infty}^2\int_Ix^mu_x^2\\&+C\|u_x\|_{L^\infty}^2\|\rho\|_{L^\infty}^2\int_Ix^m\rho
u_t^2+C\|\rho\|_{L^\infty}^3\|u\|_{L^\infty}^2\int_Ix^mu_{xt}^2\\&+C\|\rho_t\|_{L^\infty}^2\|v\|_{L^\infty}^2\int_Ix^m\rho
v^2 +C\|v\|_{L^\infty}^2\|\rho\|_{L^\infty}^2\int_Ix^m\rho v_t^2\\
\le&\frac{1}{8}\int_Ix^m\rho^3u_{tt}^2+C\int_Ix^mu_{xt}^2 +C,
\end{split}\ee where we have used Cauchy inequality, Corollaries \ref{cor:tx H^2 of v}, \ref{cor:xt theta xx}, \ref{cor: rho
w1,infty}, and Lemmas \ref{le:rho vt2}, \ref{le:rho ut}.

For $VV_4$, using Cauchy inequality and Lemma \ref{le:rho ut}, we
have \be\label{VV 4}\begin{split}
VV_4\le\frac{1}{8}\int_Ix^m\rho^3u_{tt}^2+C\int_Ix^m\rho
u_{t}^2\le\frac{1}{8}\int_Ix^m\rho^3u_{tt}^2+C.
\end{split}
\ee For $VV_5$, recalling $P=\rho Q+P_c$, we have \be\label{VV
5}\begin{split} VV_5=&-\int_Ix^m\rho^2u_{tt}\left(\rho_{xt}
Q+\rho_xQ^\prime\theta_t+\rho_t Q^\prime\theta_x+\rho
Q^{\prime\prime}\theta_t\theta_x+\rho
Q^\prime\theta_{xt}+P_c^\prime\rho_{xt}+P_c^{\prime\prime}\rho_{x}\rho_t\right)\\
\le&\frac{1}{8}\int_Ix^m\rho^3u_{tt}^2+C\int_Ix^m\rho\rho_{xt}^2
Q^2+\int_Ix^m\rho\rho_x^2(Q^\prime)^2\theta_t^2+\int_Ix^m\rho\rho_t^2
(Q^\prime)^2\theta_x^2\\&+\int_Ix^m\rho^3
(Q^{\prime\prime})^2\theta_t^2\theta_x^2+\int_Ix^m\rho^3
(Q^\prime)^2\theta_{xt}^2+\int_Ix^m\rho(P_c^\prime)^2\rho_{xt}^2+\int_Ix^m\rho(P_c^{\prime\prime})^2\rho_{x}^2\rho_t^2\\
\le&\frac{1}{8}\int_Ix^m\rho^3u_{tt}^2+C,
\end{split}
\ee  where we have used Cauchy inequality, ($A_3$), ($A_4$),
Corollaries \ref{cor: rho w1,infty}, \ref{cor:theta w 1,infty},
\ref{cor:theta xxx}, and Lemma \ref{le:dt rho thtea t}.

Putting (\ref{VV 1 and 2}), (\ref{VV 3}), (\ref{VV 4}) and
(\ref{VV 5}) into (\ref{dt rho uxt-1}), we have \be\label{dt rho
uxt-2}
\begin{split}&\int_Ix^m\rho^3u_{tt}^2+\beta\frac{d}{dt}\int_Ix^m\rho^2
u_{xt}^2 \le C\int_Ix^mu_{xt}^2 +C.\end{split} \ee By (\ref{dt rho
uxt-2}), Lemma \ref{le:rho ut} and Gronwall inequality, we
complete the proof of Lemma \ref{le:rho uxt}.\endpf

\begin{corollary}\label{cor: x uxxx} Under the conditions
of Theorem \ref{th:1.2}, we have for any $t\in[0,T]$
$$\int_Ix^mu_{xxx}^2\le C.$$
\end{corollary}
\pf By (\ref{u xxx}), (\ref{3}), Lemmas \ref{le:rho ut}, \ref{le:
rho xx}, \ref{le:rho x12infty}, \ref{le:rho uxt}, Corollaries
\ref{cor:theta w 1,infty}, we get \be\label{u xxx+1}\begin{split}
\int_Ix^mu_{xxx}^2
\le&C\int_Ix^m\rho_x^2u_t^2+C\int_Ix^m\rho^2u_{xt}^2+C\int_Ix^mP_{xx}^2+C\\
\le& C\int_Ix^m\rho|(\sqrt{\rho})_x|^2u_t^2+C\le C.
\end{split}
\ee
\endpf
\begin{lemma}\label{le:rho vxt}Under the conditions
of Theorem \ref{th:1.2}, we have for any $t\in[0,T]$
$$
\int_Ix^m\rho^2(v_{xt}^2+w_{xt}^2)+\int_{Q_T}x^m\rho^3(v_{tt}^2+w_{tt}^2)\le C.
$$
\end{lemma}
\pf Multiplying (\ref{rho vtt}) by $x^m\rho^2v_{tt}$, integrating
by parts over $I$, we have \be\label{dt rho vxt}\begin{split}
&\int_Ix^m\rho^3v_{tt}^2+\frac{\mu}{2}\frac{d}{dt}\int_Ix^m\rho^2
v_{xt}^2\\=&\mu\int_Ix^m\rho\rho_tv_{xt}^2-2\mu\int_Ix^m\rho\rho_xv_{xt}v_{tt}-\int_Ix^m\rho^2v_{tt}(\rho_tv_t+\rho_tuv_x+\rho
u_tv_x+\rho
uv_{xt})\\&-\int_Ix^m\rho^2v_{tt}(x^{-1}\rho_tuv+x^{-1}\rho
u_tv+x^{-1}\rho uv_t)-\mu\int_Ix^m\rho^2v_{tt}\frac{mv_t}{x^2}\\
=&\sum\limits_{i=1}^5VVI_i.
\end{split}
\ee For $VVI_1$ and $VVI_2$, using Cauchy inequality, Corollary
\ref{cor: rho w1,infty} and Lemma \ref{le:rho x12infty}, we have
\be\label{VVI 1 and 2}\begin{split}
VVI_1+VVI_2 \le&\|\rho\|_{L^\infty}\|\rho_t\|_{L^\infty}\int_Ix^mv_{xt}^2-4\mu\int_Ix^m\rho^\frac{3}{2}(\sqrt{\rho})_xv_{xt}v_{tt}\\
\le&\frac{1}{8}\int_Ix^m\rho^3v_{tt}^2+C\int_Ix^mv_{xt}^2.
\end{split}
\ee For $VVI_3$, using Cauchy inequality, Corollary \ref{cor: rho
w1,infty} again, together with Corollaries \ref{cor:xt L^infty of
u}, \ref{cor:vxx}, and Lemmas \ref{le:rho vt2}, \ref{le:rho ut},
we have \be\label{VVI 3}\begin{split} VVI_3
\le&\frac{1}{8}\int_Ix^m\rho^3v_{tt}^2+C\|\rho_t\|_{L^\infty}^2\int_Ix^m\rho
v_t^2+C\|\rho\|_{L^\infty}\|\rho_t\|_{L^\infty}^2\|u\|_{L^\infty}^2\int_Ix^mv_x^2\\&+C\|v_x\|_{L^\infty}^2\|\rho\|_{L^\infty}^2\int_Ix^m\rho
u_t^2+C\|\rho\|_{L^\infty}^3\|u\|_{L^\infty}^2\int_Ix^mv_{xt}^2\\
\le&\frac{1}{8}\int_Ix^m\rho^3v_{tt}^2+C\int_Ix^mv_{xt}^2+C.
\end{split}
\ee Similarly, for $VVI_4$ and $VVI_5$, we have \be\label{VVI
4}\begin{split} VVI_4
\le&\frac{1}{8}\int_Ix^m\rho^3v_{tt}^2+C\|\rho\|_{L^\infty}\|u\|_{L^\infty}^2\|v\|_{L^\infty}^2\int_Ix^m\rho_t^2+C\|\rho\|_{L^\infty}^2\|v\|_{L^\infty}^2\int_Ix^m\rho
u_t^2\\&+C\|\rho\|_{L^\infty}^2\|u\|_{L^\infty}^2\int_Ix^m\rho
v_t^2\\ \le&\frac{1}{8}\int_Ix^m\rho^3v_{tt}^2+C;
\end{split}
\ee  \be\label{VVI 5}\begin{split}
VVI_5\le\frac{1}{8}\int_Ix^m\rho^3v_{tt}^2+C\int_Ix^m\rho v_t^2
\le\frac{1}{8}\int_Ix^m\rho^3v_{tt}^2+C.
\end{split}
\ee Putting (\ref{VVI 1 and 2}), (\ref{VVI 3}), (\ref{VVI 4}) and
(\ref{VVI 5}) into (\ref{dt rho vxt}), and using Lemma \ref{le:rho
vt2}, Gronwall inequality, we have \be\label{dt rho
vxt+1}\begin{split} \int_Ix^m\rho^2
v_{xt}^2+\int_{Q_T}x^m\rho^3v_{tt}^2\le C.
\end{split}
\ee
Similarly, we get
\bex
\int_Ix^m\rho^2
w_{xt}^2+\int_{Q_T}x^m\rho^3w_{tt}^2\le C.
\eex
\endpf

\begin{corollary}\label{cor: x wxxx and v xxx}Under the conditions
of Theorem \ref{th:1.2}, we have for any $t\in[0,T]$
$$\int_Ix^m(v_{xxx}^2+w_{xxx}^2)+\int_{Q_T}(v_{xxxx}^2+w_{xxxx}^2)\le C.$$
\end{corollary}
\pf
 Similar to Corollary \ref{cor: x uxxx}, we get
$$
\int_Ix^m(v_{xxx}^2+w_{xxx}^2)\le C.
$$
By (\ref{rho vtt}), Corollaries \ref{cor:xt L^infty of u},
\ref{cor:vxx}, \ref{cor: rho w1,infty}, and Lemmas \ref{le:rho
vt2}, \ref{le:rho ut}, \ref{le:rho vxt}, we have
$$
\int_{Q_T}\rho v_{xxt}^2\le C.
$$
This together with (\ref{non-1.2})$_3$ gives
$$
\int_{Q_T}v_{xxxx}^2\le C.
$$
Similarly, we can get
$$
\int_{Q_T}w_{xxxx}^2\le C.
$$
\endpf
\begin{lemma}\label{le: rho xxx} Under the conditions
of Theorem \ref{th:1.2}, we have for any $t\in[0,T]$
$$
\int_Ix^m\rho_{xxx}^2+\int_{Q_T}x^m(u_{xxxx}^2+\theta_{xxxx}^2)\le
C.
$$
\end{lemma}
\pf
 Differentiating (\ref{rho xxt+.})
with respect to $x$, multiplying it by $x^m\rho_{xxx}$, and
integrating by parts over $I$, we have \bex
\begin{split}&\frac{1}{2}\frac{d}{dt}\int_Ix^m\rho_{xxx}^2
=\int_Ix^m\rho_{xxx}\Big[-mx^{-1}\rho
u_{xxx}-3mx^{-1}\rho_xu_{xx}-3mx^{-1}\rho_{xx}u_x+\frac{3m\rho
u_{xx}}{x^2}\\&-\frac{6m\rho
u_x}{x^3}-mx^{-1}\rho_{xxx}u+\frac{3m\rho_{xx}
u}{x^2}+\frac{6m\rho_x u_x}{x^2}-\frac{6m\rho_x
u}{x^3}+\frac{6m\rho u}{x^4}-4\rho_{xxx}u_x-6\rho_{xx}u_{xx}\\&
-4\rho_xu_{xxx}\Big]+\frac{1}{2}\int_Ix^m \rho_{xxx}^2
 u_x+\frac{1}{2}\int_Imx^{m-1} \rho_{xxx}^2 u-\int_Ix^m\rho\rho_{xxx}
u_{xxxx}.\end{split} \eex By Sobolev inequality, Cauchy
inequality, Corollaries \ref{cor:xt theta xx}, \ref{cor: rho
w1,infty}, \ref{cor: x uxxx}, and Lemma \ref{le: rho xx},  we get
\be\label{dt rho xxx} \frac{d}{dt}\int_Ix^m\rho_{xxx}^2\le
C\int_Ix^m\rho_{xxx}^2+C\int_Ix^mu_{xxxx}^2+C.\ee
 Differentiating (\ref{u xxx+.}) with respect to $x$, we have
\be\label{u xxxx+.} \begin{split} \beta
u_{xxxx}=&\rho_{xx}u_t+2\rho_xu_{xt} +\rho
u_{xxt}+(\rho_xuu_x+\rho u_x^2+\rho uu_{xx})_x-(\frac{2\rho
vv_x}{x}+\frac{\rho_xv^2}{x}-\frac{\rho v^2}{x^2})_x\\&+P_{xxx}
-\frac{m\beta u_{xxx}}{x}+\frac{3m\beta u_{xx}}{x^2}-\frac{6m\beta
u_x}{x^3}+\frac{6m\beta u}{x^4}. \end{split}\ee Using Lemma
\ref{le: rho xx}, Corollaries \ref{cor:vxx}, \ref{cor:xt theta
xx}, \ref{cor: rho w1,infty}, \ref{cor:theta w 1,infty},
\ref{cor:theta xxx}, \ref{cor: x uxxx}, and ($A_3$), ($A_4$), we
obtain \be\label{int u xxxx-1}\begin{split} \int_Ix^mu_{xxxx}^2
\le& C\int_Ix^m|P_{xxx}|^2+C\int_Ix^m (\rho u_{xxt}^2+u_{xt}^2)+C\\
\le&C\int_Ix^m\rho_{xxx}^2+C\int_Ix^m(\rho u_{xxt}^2+u_{xt}^2)+C.
\end{split}\ee Substituting (\ref{int u xxxx-1}) into (\ref{dt rho
xxx}), and using (\ref{rho utt+.=}), Lemmas \ref{le:rho ut},
\ref{le:rho uxt}, Gronwall inequality, we obtain \be\label{int rho
xxx} \int_Ix^m\rho_{xxx}^2\le C. \ee By (\ref{int u xxxx-1}) and
(\ref{int rho xxx}), we have
$$
\int_{Q_T}x^mu_{xxxx}^2\le C.
$$
Differentiating (\ref{theta xxx+.}) with respect to $x$, and using
($A_4$), ($A_5$), Corollaries \ref{cor:vxx}, \ref{cor:wxx},
\ref{cor:xt theta xx}, \ref{cor: rho w1,infty}, \ref{cor:t theta
t}, \ref{cor:xt theta xt}, \ref{cor:theta w 1,infty},
\ref{cor:theta xxx}, and Lemmas \ref{le: rho xx}, \ref{le:rho
theta xt}, we have
$$
\int_{Q_T}x^m\theta_{xxxx}^2\le C.
$$

The proof of Lemma \ref{le: rho xxx} is complete. \endpf

\section{Proof of Theorem \ref{blowup-th:1.1}}
\setcounter{equation}{0} \setcounter{theorem}{0}

Let $0<T^*<\infty$ be the maximum time of existence of strong
solution $(\rho, u, \theta)$ to (\ref{full
N-S+1})-(\ref{non-boundary}). Namely, $(\rho, u, \theta)$ is a
strong solution to (\ref{full N-S+1})-(\ref{non-boundary}) in
$\mathbb{R}^3\times [0, T]$ for any $0<T<T^*$, but not a strong
solution in $\mathbb{R}^3\times [0, T^*]$. We shall prove Theorem
\ref{blowup-th:1.1} by using a contradiction argument. Suppose that
(\ref{non-result}) were false, i.e. \beq\label{blowup-2.1}
M:=\lim\sup\limits_{t\nearrow
T^*}(\|\rho(t)\|_{L^\infty}+\int_0^t\|\rho\theta(s)\|_{L^\frac{12}{5}}^4ds)<\infty.
\eeq The goal is to show that under the assumption
(\ref{blowup-2.1}), there is a bound $C>0$ depending only on $M,
\rho_0, u_0, \theta_0, \mu,\lambda, \kappa$, and $T^*$ such that
\beq\label{non-uniform_est1} \sup_{0\le
t<T^*}\|\theta(t)\|_{L^\infty}\le C. \eeq With
(\ref{non-uniform_est1}) and (\ref{blowup-2.1}), we showed in our
previous paper \cite{Wen-Zhu 3} that $T^*$ is not the maximum time,
which is the desired contradiction.

Throughout the rest of the section, we denote by $C$ a generic
constant depending only on $\rho_0$, $u_0$, $\theta_0$, $T^*$, $M$,
$\lambda$, $\mu$, $\kappa$. We denote by
$$A\lesssim B$$ if there exists a generic constant $C$ such that $A\leq C B$.
\begin{lemma}\label{blowup-le:2.1}
Under the conditions of Theorem \ref{blowup-th:1.1} and
(\ref{blowup-2.1}), it holds that \be\label{non-energy inequality}
\displaystyle\sup\limits_{0\leq t\leq T}\int_{\mathbb{R}^3} \rho\le
C,\ \mathrm{for}\ \mathrm{any}\ T\in[0,T^*). \ee
\end{lemma}
\pf Integrating (\ref{full N-S+1})$_1$ over
$\mathbb{R}^3\times[0,t]$, for $t<T^*$, and using the assumption
$\rho_0\in L^1$, we get (\ref{non-energy inequality}).
\endpf
\begin{lemma}\label{blowup-le:2.2}
Under the conditions of Theorem \ref{blowup-th:1.1} and
(\ref{blowup-2.1}), if $3\mu>\lambda$, it holds that \be\label{u
2nabla u 2} \sup\limits_{0\leq t\leq
T}\int_{\mathbb{R}^3}\rho|u|^4+\int_0^T\int_{\mathbb{R}^3}|u|^2|\nabla
u|^2\,dx\leq C, \ee for any $T\in(0,T^*)$.
\end{lemma}
\pf  The detailed proof of Lemma \ref{blowup-le:2.2} could be found
in \cite{Wen-Zhu 3}, which might be slightly modified.
\endpf
\begin{lemma}\label{blowup-le: int nabla u}Under the conditions of Theorem \ref{blowup-th:1.1} and (\ref{blowup-2.1}), it holds that for any $T\in[0,T^*)$
\be\label{H 1 of u}\sup\limits_{0\le t\le
T}\int_{\mathbb{R}^3}(\rho|\theta|^2+|\nabla u|^2)\,
dx+\int_0^T\int_{\mathbb{R}^3}(|\nabla\theta|^2+\rho
|u_t|^2)\,dxdt\le C.\ee
\end{lemma}
\pf Multiplying (\ref{full N-S+1})$_2$ by $u_t$, and integrating by
parts over $\mathbb{R}^3$, we have \beq\label{dt nabla u
2-1}\begin{split} &\int_{\mathbb{R}^3}\rho
|u_t|^2+\frac{1}{2}\frac{d}{dt}\int_{\mathbb{R}^3}\left(\mu|\nabla
u|^2+(\mu+\lambda)|\mathrm{div}u|^2\right)\\=&\frac{d}{dt}\int_{\mathbb{R}^3}P\mathrm{div}u-\frac{1}{2(2\mu+\lambda)}
\frac{d}{dt}\int_{\mathbb{R}^3}P^2
-\frac{1}{2\mu+\lambda}\int_{\mathbb{R}^3}P_tG-\int_{\mathbb{R}^3}\rho
u\cdot\nabla u \cdot u_t\\=&\sum\limits_{i=1}^4VVII_i,
\end{split}
\eeq where $G=(2\mu+\lambda)\mathrm{div}u-P$.

Recalling $P=\rho\theta$, we obtain from (\ref{full N-S+1})$_1$ and
(\ref{full N-S+1})$_3$ \beq\label{equation of P t}\begin{split}
P_t=-\mathrm{div}(P u)-\rho\theta\mathrm{div}u+\mu\left(\nabla
u+(\nabla u)^\prime\right):\nabla
u+\lambda\mathrm{div}u\mathrm{div}u+\Delta\theta.
\end{split}
\eeq Substituting (\ref{equation of P t}) into $VVII_3$, and using
integration by parts and H\"older inequality, we have
 \beq\label{VVII 3}\begin{split}
VVII_3=&-\frac{1}{2\mu+\lambda}\int_{\mathbb{R}^3}P u\cdot\nabla
G+\frac{1}{2\mu+\lambda}\int_{\mathbb{R}^3}\rho\theta\mathrm{div}u
G\\&+\frac{\mu}{2\mu+\lambda}\int_{\mathbb{R}^3}\left(\nabla
u+(\nabla u)^\prime\right):\left(\nabla G\otimes
u\right)+\frac{\lambda}{2\mu+\lambda}\int_{\mathbb{R}^3}\mathrm{div}u
u\cdot\nabla
G\\&+\frac{1}{2\mu+\lambda}\int_{\mathbb{R}^3}\big(\mu\Delta
u+(\mu+\lambda)\nabla \mathrm{div}u\big)\cdot uG
+\frac{1}{2\mu+\lambda}\int_{\mathbb{R}^3}\nabla\theta\cdot\nabla
G\\ &\le C\|\rho u\theta\|_{L^2}\|\nabla G\|_{L^2}+\frac{1}{2\mu+\lambda}\int_{\mathbb{R}^3}\rho\theta\mathrm{div}u
G+C\|\nabla
G\|_{L^2}\big\|u|\nabla u|\big\|_{L^2}\\&+C\|\nabla
G\|_{L^2}\|\nabla\theta\|_{L^2}+\frac{1}{2\mu+\lambda}\int_{\mathbb{R}^3}\big(\mu\Delta
u+(\mu+\lambda)\nabla \mathrm{div}u\big)\cdot uG.
\end{split}
\eeq Substituting (\ref{full N-S+1})$_2$ into (\ref{VVII 3}), and
using Sobolev inequality, (\ref{blowup-2.1}) and integration by
parts, we have \beq\label{VVII 3+1}\begin{split} VVII_3\le&
C\|\nabla G\|_{L^2}\Big(\|\rho u\theta\|_{L^2}+\big\|u|\nabla
u|\big\|_{L^2}+\|\nabla\theta\|_{L^2}\Big)+\frac{1}{2\mu+\lambda}\int_{\mathbb{R}^3}\rho
u_t\cdot uG\\&+\frac{1}{2\mu+\lambda}\int_{\mathbb{R}^3}\rho
u\cdot\nabla u\cdot uG-\frac{1}{2\mu+\lambda}\int_{\mathbb{R}^3} Pu\cdot\nabla G\\
\le&C\|\nabla G\|_{L^2}\Big(\|\rho u\theta\|_{L^2}+\big\|u|\nabla
u|\big\|_{L^2}+\|\nabla\theta\|_{L^2}\Big)+\frac{1}{6}\int_{\mathbb{R}^3}\rho
|u_t|^2\\&+C\int_{\mathbb{R}^3}\rho|u|^2|G|^2+C\big\|u|\nabla u|\big\|_{L^2}^2\\
\le&C\|\nabla G\|_{L^2}\Big(\|\rho u\theta\|_{L^2}+\big\|u|\nabla
u|\big\|_{L^2}+\|\nabla\theta\|_{L^2}\Big)+\frac{1}{6}\int_{\mathbb{R}^3}\rho
|u_t|^2\\&+C\int_{\mathbb{R}^3}|u|^2|\nabla
u|^2+C\int_{\mathbb{R}^3}\rho|u|^2|\rho\theta|^2.
\end{split}
\eeq
 Taking $div$ on both side of
(\ref{full N-S+1})$_2$, we get \be\label{equation of G} \Delta
G=\mathrm{div}(\rho u_t+\rho u\cdot\nabla u). \ee
 By (\ref{equation of G}) and the standard
$L^2$-estimates together with (\ref{blowup-2.1}), we get
\beq\label{H 1 of G}
\begin{split}\|\nabla G\|_{L^{2}}\les \|\rho u_t\|_{L^{2}}+\|\rho u\cdot\nabla u\|_{L^2}\les
 \|\sqrt{\rho} u_t\|_{L^{2}}+\big\| |u||\nabla u|\big\|_{L^2}.
\end{split}
\eeq
 Substituting (\ref{H 1 of G}) into (\ref{VVII 3+1}), and using
Cauchy inequality, we have \beq\label{VVII 3+2-1}\begin{split}
VVII_3 \le&C\|\rho u\theta\|_{L^2}^2+C\big\|u|\nabla
u|\big\|_{L^2}^2+C\|\nabla\theta\|_{L^2}^2+\frac{1}{3}\int_{\mathbb{R}^3}\rho
|u_t|^2.
\end{split}
\eeq For the first term of the right hand side of (\ref{VVII
3+2-1}), using H\"older inequality, Sobolev inequality and Cauchy
inequality, we have \be\label{rho u theta}\begin{split} \|\rho
u\theta\|_{L^2}^2\le&\big\||u|^2\big\|_{L^6}\|\rho\theta\|_{L^\frac{12}{5}}^2\le
C\big\|u|\nabla u|\big\|_{L^2}\|\rho\theta\|_{L^\frac{12}{5}}^2\\
\le&C\big\|u|\nabla
u|\big\|_{L^2}^2+C\|\rho\theta\|_{L^\frac{12}{5}}^4.
\end{split}
\ee Substituting (\ref{rho u theta}) into (\ref{VVII 3+2-1}), we
have
 \beq\label{VVII 3+2}\begin{split} VVII_3
\le&C\|\rho\theta\|_{L^\frac{12}{5}}^4+C\big\|u|\nabla
u|\big\|_{L^2}^2+C\|\nabla\theta\|_{L^2}^2+\frac{1}{3}\int_{\mathbb{R}^3}\rho
|u_t|^2.
\end{split}
\eeq
 For $VVII_4$, using Cauchy inequality and (\ref{blowup-2.1}), we have \be\label{VVII
4}\begin{split} VVII_4\le&\frac{1}{6}\int_{\mathbb{R}^3}\rho
|u_t|^2+C\int_{\mathbb{R}^3} |u|^2 |\nabla u|^2.
\end{split}
\ee Putting (\ref{VVII 3+2}) and (\ref{VVII 4}) into (\ref{dt nabla
u 2-1}), and integrating it over $[0,t]$, for $t<T^*$, we have
\bex\begin{split} &\int_0^t\int_{\mathbb{R}^3}\rho
|u_t|^2+\int_{\mathbb{R}^3}\left(\mu|\nabla
u|^2+(\mu+\lambda)|\mathrm{div}u|^2\right)\\
\le&2\int_{\mathbb{R}^3}P\mathrm{div}u+C\int_0^t\|\nabla\theta\|_{L^2}^2+C\\
\le&(\mu+\lambda)\int_{\mathbb{R}^3}|\mathrm{div}u|^2+C\big(\int_{\mathbb{R}^3}\rho\theta^2+\int_0^t\int_{\mathbb{R}^3}|\nabla\theta|^2\big)+C,
\end{split}
\eex where we have used Cauchy inequality, (\ref{blowup-2.1}) and
(\ref{u 2nabla u 2}). Therefore, \be\label{zy}\begin{split}
\int_0^t\int_{\mathbb{R}^3}\rho|u_t|^2+\int_{\mathbb{R}^3}|\nabla
u|^2\leq
C\big(\int_{\mathbb{R}^3}\rho\theta^2+\int_0^t\int_{\mathbb{R}^3}|\nabla\theta|^2\big)+C.
\end{split}
\ee Multiplying (\ref{full N-S+1})$_3$ by $\theta$ and integrating
by parts over $\mathbb{R}^3$, we have \beq\label{zyy}\begin{split}
&\int_{\mathbb{R}^3}|\nabla\theta|^2+\frac{1}{2}\frac{d}{dt}\int_{\mathbb{R}^3}\rho|\theta|^2\\
=&-\int_{\mathbb{R}^3}\rho\theta^2\mathrm{div}u+\int_{\mathbb{R}^3}\frac{\mu}{2}|\nabla u+(\nabla u)'|^2\theta+\int_{\mathbb{R}^3}\lambda|\mathrm{div}u|^2\theta\\
=&\sum\limits_{i=1}^3VVIII_i,
\end{split}
\eeq  For $VVIII_2$ and $VVIII_3$, we have
\beq\label{VVIII_2}\begin{split}
VVIII_2+VVIII_3=&\int_{\mathbb{R}^3}\mu\left(\nabla u+(\nabla u)'\right):\nabla u\theta+\int_{\mathbb{R}^3}\lambda|\mathrm{div}u|^2\theta\\
=&-\int_{\mathbb{R}^3}\mu(\triangle u+\nabla\mathrm{div}u )\cdot u\theta-\int_{\mathbb{R}^3}\mu\left(\nabla u+(\nabla u)'\right):(\nabla\theta\otimes u)
\\&-\int_{\mathbb{R}^3}\lambda u\cdot\nabla\mathrm{div}u\theta-\int_{\mathbb{R}^3}\lambda\mathrm{div}u u\cdot\nabla \theta\\
=&-\int_{\mathbb{R}^3}(\rho u_t+\rho u\cdot \nabla u+\nabla P)\cdot
u\theta-\int_{\mathbb{R}^3}\mu\left(\nabla u+(\nabla
u)'\right):(\nabla\theta\otimes u)
\\&-\int_{\mathbb{R}^3}\lambda\mathrm{div}u u\cdot\nabla \theta\\
=&-\int_{\mathbb{R}^3}\rho u_t\cdot
u\theta-\int_{\mathbb{R}^3}\rho(u\cdot\nabla)u\cdot
u\theta+\int_{\mathbb{R}^3}\rho\theta^2\mathrm{div}u+\int_{\mathbb{R}^3}\rho\theta
u\cdot\nabla\theta
\\&-\int_{\mathbb{R}^3}\mu\left(\nabla u+(\nabla u)'\right):(\nabla\theta\otimes u)-\int_{\mathbb{R}^3}\lambda\mathrm{div}u u\cdot\nabla \theta,
\end{split}
\eeq where we have used integration by parts and (\ref{full N-S+1})$_2$.\\
Using H\"older inequality, Cauchy inequality, (\ref{blowup-2.1}) and
(\ref{u 2nabla u 2}), we have \beq\label{VVIII_2result}\begin{split}
&VVIII_2+VVIII_3\\ \leq&
\|\sqrt{\rho}u_t\|_{L^2}\|\sqrt[4]{\rho}u\|_{L^4}\|\sqrt[4]{\rho}\theta\|_{L^4}
+ \big\|u|\nabla u|\big\|_{L^2}\|\rho u\|_{L^3}\|\theta\|_{L^6}
+\int_{\mathbb{R}^3}\rho\theta^2\mathrm{div}u\\&+\|\nabla\theta\|_{L^2}\|\sqrt{\rho}u\|_{L^4}\|\sqrt{\rho}\theta\|_{L^4}
+C\big\|u|\nabla u|\big\|_{L^2}\|\nabla\theta\|_{L^2}\\
\leq&C\|\sqrt{\rho}u_t\|_{L^2}\|\sqrt[4]{\rho}\theta\|_{L^4}+\frac{1}{2}\|\nabla\theta\|_{L^2}^2+C\big\|u|\nabla
u|\big\|_{L^2}^2+\int_{\mathbb{R}^3}\rho\theta^2\mathrm{div}u+C\|\sqrt{\rho}\theta\|_{L^4}^2.
\end{split}
\eeq Substituting (\ref{VVIII_2result}) into (\ref{zyy}), we have
\beq\label{zyyy}\begin{split}
\int_{\mathbb{R}^3}|\nabla\theta|^2+\frac{d}{dt}\int_{\mathbb{R}^3}\rho|\theta|^2\leq
C\|\sqrt{\rho}u_t\|_{L^2}\|\sqrt[4]{\rho}\theta\|_{L^4}+C\big\|u|\nabla
u|\big\|_{L^2}^2+C\|\sqrt{\rho}\theta\|_{L^4}^2.
\end{split}
\eeq Integrating (\ref{zyyy}) over $[0,t]$ ($t<T^*$), and using
(\ref{u 2nabla u 2}), we have \beq\label{sy}\begin{split}
\int_0^t\int_{\mathbb{R}^3}|\nabla\theta|^2+\int_{\mathbb{R}^3}\rho|\theta|^2\leq
C\int_0^t\|\sqrt{\rho}u_t\|_{L^2}\|\sqrt[4]{\rho}\theta\|_{L^4}+C\int_0^t\|\sqrt{\rho}\theta\|_{L^4}^2+C.
\end{split}
\eeq Multiplying (\ref{sy}) by $2C$, and adding the resulting
inequality into (\ref{zy}), we have \bex\begin{split}
&\int_0^t\int_{\mathbb{R}^3}\rho|u_t|^2+\int_{\mathbb{R}^3}|\nabla u|^2+\int_0^t\int_{\mathbb{R}^3}|\nabla\theta|^2+\int_{\mathbb{R}^3}\rho|\theta|^2\\
\les&
\int_0^t\|\sqrt{\rho}u_t\|_{L^2}\|\sqrt[4]{\rho}\theta\|_{L^4}+\int_0^t\|\sqrt{\rho}\theta\|_{L^4}^2+1\\
\le&\frac{1}{2}\int_0^t\int_{\mathbb{R}^3}\rho|u_t|^2+C\int_0^t\|\sqrt{\rho}\theta\|_{L^2}^\frac{1}{2}\|\sqrt[6]{\rho}\theta\|_{L^6}^\frac{3}{2}+C\\
\le&\frac{1}{2}\int_0^t\int_{\mathbb{R}^3}\rho|u_t|^2+\frac{1}{2}\int_0^t\int_{\mathbb{R}^3}|\nabla\theta|^2+C\int_0^t\int_{\mathbb{R}^3}\rho|\theta|^2+C,
\end{split}
\eex where we have used Young inequality, H\"older inequality,
Sobolev inequality and (\ref{blowup-2.1}).
 This, together with Gronwall inequality, gives (\ref{H 1 of
u}).
\endpf

\begin{lemma}\label{blowup-le: int rho u t}Under the conditions of Theorem \ref{blowup-th:1.1} and (\ref{blowup-2.1}), it holds that for any $t\in(0,T^*)$
\be\label{H 2 of u}\int_{\mathbb{R}^3}(|\nabla
\theta|^2+\rho|\u|^2)+\int_0^t\int_{\mathbb{R}^3}(\rho
|\dot{\theta}|^2+|\nabla\u|^2)\le C.\ee
\end{lemma}
\pf By \cite{Wen-Zhu 3} (see (4.35) therein), we accurately
have\beq\label{dt rho u t-1}
\begin{split}
&\frac{1}{2}\frac{d}{dt}\int_{\mathbb{R}^3}\rho|\u|^2+\int_{\mathbb{R}^3}\left(\mu|\nabla\u|^2+(\mu+\lambda)|\di\u|^2\right)\\
=&\int_{\mathbb{R}^3}\left( P_t\di\u + u\otimes\nabla P:\nabla\u\right)+\mu\int_{\mathbb{R}^3}\Big(\di(\de u\otimes u)-\de(u\cdot\nabla u)\Big)\cdot\u \\
+&(\mu+\lambda)\int_{\mathbb{R}^3}\Big(\di(\nabla\di u\otimes
u)-\nabla\di(u\cdot\nabla u)\Big)\cdot\u =\sum\ls_{i=1}^3VVIV_i.
\end{split}
\eeq For $VVIV_1$, using (\ref{full N-S+1})$_1$, integration by
parts, (\ref{blowup-2.1}) and H\"older inequality, we have
\beq\label{VVIV 1}
\begin{split}
VVIV_1=&\int_{\mathbb{R}^3}\Big( (\rho\theta)_t\di\u-\rho\theta(\nabla u)^t:\nabla\u -\rho\theta u\cdot\nabla\di\u \Big)\\
=&\int_{\mathbb{R}^3}\Big(\rho\dot{\theta}\di\u-\rho\theta(\nabla u)^t:\nabla\u  \Big)\\
 \les&\|\sqrt{\rho}\dot{\theta}\|_{L^2}\|\nabla\u\|_{L^2}+\|\sqrt[4]{\rho}\theta\|_{L^4}\|\nabla u\|_{L^4}\|\nabla \u\|_{L^2}.
\end{split}
\eeq For $VVIV_2$ and $VVIV_3$, by \cite{Wen-Zhu 3} (see (4.37) and
(4.38) therein), we have we have \beq\label{VVIV 2, 3}
\begin{split}
VVIV_2+VVIV_3\les\|\nabla \u\|_{L^2}\|\nabla u\|_{L^4}^2.
\end{split}
\eeq Substituting (\ref{VVIV 1}) and (\ref{VVIV 2, 3}) into (\ref{dt
rho u t-1}), and using Cauchy inequality and (\ref{blowup-2.1}), we
have \bex
\begin{split}
\frac{1}{2}\frac{d}{dt}\int_{\mathbb{R}^3}\rho|\u|^2+\int_{\mathbb{R}^3}\left(\mu|\nabla\u|^2+(\mu+\lambda)|\di\u|^2\right)
\le\frac{\mu}{2}\|\nabla
\u\|_{L^2}^2+C\|\sqrt{\rho}\dot{\theta}\|_{L^2}^2+C\|\sqrt[4]{\rho}\theta\|_{L^4}^4+C\|\nabla
u\|_{L^4}^4.
\end{split}
\eex Integrating this inequality over $[0,t]$ for $t\in(0,T^*)$, and
using (\ref{blowup-2.1}), (\ref{non-energy inequality}), H\"older
inequality and Sobolev inequality, we have \beq\label{dt rho u t-2}
\begin{split}
\int_{\mathbb{R}^3}\rho|\u|^2+\int_0^t\int_{\mathbb{R}^3}|\nabla\u|^2\le
C\int_0^t\|\sqrt{\rho}\dot{\theta}\|_{L^2}^2+C\int_0^t\left(\|\nabla
u\|_{L^4}^4+\|\nabla \theta\|_{L^2}^4\right).
\end{split}
\eeq The next step is to get some estimates for $\theta$. We rewrite
(\ref{full N-S+1})$_3$ as follows: \be\label{rho theta+ =} \rho
\dot{\theta}+\rho\theta\mathrm{div}u=\frac{\mu}{2}\left|\nabla
u+(\nabla u)^\prime\right|^2+\lambda(\mathrm{div}u)^2+\Delta\theta.
\ee Multiplying (\ref{rho theta+ =}) by $\dot{\theta}$ and
integrating by parts over $\mathbb{R}^3$, we have \be\label{dt nabla
theta}\begin{split} &\int_{\mathbb{R}^3}\rho
|\dot{\theta}|^2+\frac{1}{2}\frac{d}{dt}\int_{\mathbb{R}^3}|\nabla
\theta|^2\\
=&-\int_{\mathbb{R}^3}\rho\theta\mathrm{div}u\dot{\theta}+\int_{\mathbb{R}^3}\left(\frac{\mu}{2}\left|\nabla
u+(\nabla
u)^\prime\right|^2+\lambda(\mathrm{div}u)^2\right)\theta_t\\&+\int_{\mathbb{R}^3}\left(\frac{\mu}{2}\left|\nabla
u+(\nabla
u)^\prime\right|^2+\lambda(\mathrm{div}u)^2\right)u\cdot\nabla\theta+\int_{\mathbb{R}^3}\Delta\theta
u\cdot\nabla\theta\\=& \sum\limits_{i=1}^4VVV_i.
\end{split}
\ee For $VVV_1$, using H\"older inequality, (\ref{blowup-2.1}) and
Cauchy inequality, we have \be\label{VVV 1}\begin{split} VVV_1\leq
C\|\sqrt{\rho}\dot{\theta}\|_{L^2}\|\nabla
u\|_{L^4}\|\rho^\frac{1}{4}\theta\|_{L^4}
\le\frac{1}{8}\|\sqrt{\rho}\dot{\theta}\|_{L^2}^2+C\|\nabla
u\|_{L^4}^4+C\|\nabla\theta\|_{L^2}^4.
\end{split}
\ee For $VVV_2$, we have \bex\begin{split}
VVV_2=&\frac{d}{dt}\int_{\mathbb{R}^3}\left(\frac{\mu}{2}\left|\nabla u+(\nabla u)^\prime\right|^2+\lambda(\mathrm{div}u)^2\right)\theta-\mu\int_{\mathbb{R}^3}\left(\nabla u+(\nabla u)^\prime\right):\left(\nabla u_t+(\nabla u_t)^\prime\right)\theta\\
&-2\lambda\int_{\mathbb{R}^3}\mathrm{div}u\mathrm{div}u_t\theta\\
=&\frac{d}{dt}\int_{\mathbb{R}^3}\left(\frac{\mu}{2}\left|\nabla
u+(\nabla u)^\prime\right|^2+\lambda(\mathrm{div}u)^2\right)\theta
-\mu\int_{\mathbb{R}^3}\left(\nabla u+(\nabla u)^\prime\right):\left(\nabla \u+(\nabla \u)^\prime\right)\theta\\
&+\mu\int_{\mathbb{R}^3}\left(\nabla u+(\nabla u)^\prime\right):\left(\nabla u\cdot\nabla u+(\nabla u\cdot\nabla u)^\prime\right)\theta\\
&+\mu\int_{\mathbb{R}^3}\left(\nabla u+(\nabla u)^\prime\right)\cdot
(u\cdot\nabla)\left(\nabla u+(\nabla u)^\prime\right)\theta
-2\lambda\int_{\mathbb{R}^3}\mathrm{div}u\mathrm{div}\dot{u}\theta\\
&+2\lambda\int_{\mathbb{R}^3}\mathrm{div}u(\nabla u)^\prime:\nabla
u\theta+2\lambda\int_{\mathbb{R}^3}u\cdot\nabla
\mathrm{div}u\mathrm{div}u\theta.
\end{split}
\eex Using integration by parts, we have \be\label{VVV
2-1}\begin{split}
VVV_2=&\frac{d}{dt}\int_{\mathbb{R}^3}\left(\frac{\mu}{2}\left|\nabla
u+(\nabla u)^\prime\right|^2+\lambda(\mathrm{div}u)^2\right)\theta
-\mu\int_{\mathbb{R}^3}\left(\nabla u+(\nabla u)^\prime\right):\left(\nabla \u+(\nabla \u)^\prime\right)\theta\\
&+\mu\int_{\mathbb{R}^3}\left(\nabla u+(\nabla
u)^\prime\right):\left(\nabla u\cdot\nabla u+(\nabla u\cdot\nabla
u)^\prime\right)\theta-\mu\int_{\mathbb{R}^3}\frac{|\nabla u+(\nabla
u)^\prime|^2}{2}\mathrm{div}u\theta\\
&-\mu\int_{\mathbb{R}^3}\frac{|\nabla u+(\nabla
u)^\prime|^2}{2}u\cdot\nabla\theta
-2\lambda\int_{\mathbb{R}^3}\mathrm{div}u\mathrm{div}\dot{u}\theta+2\lambda\int_{\mathbb{R}^3}\mathrm{div}u(\nabla
u)^\prime:\nabla
u\theta\\
&-\lambda\int_{\mathbb{R}^3} (\mathrm{div}u)^3\theta-\lambda\int_{\mathbb{R}^3}|\mathrm{div}u|^2u\cdot\nabla\theta\\
=&\sum\limits_{i=1}^9VVV_{2,i}.
\end{split}
\ee\\
For $VVV_{2,2}$ and $VVV_{2,6}$, using H\"older inequality, Sobolev
inequality, we have \be\label{VVV 2,2 and 2,6}\begin{split}
VVV_{2,2}+VVV_{2,6}\lesssim\|\nabla \u\|_{L^2}\|\nabla
u\|_{L^3}\|\theta\|_{L^6}\lesssim\|\nabla \u\|_{L^2}\|\nabla
u\|_{L^3}\|\nabla \theta\|_{L^2}.
\end{split}
\ee Since $\nabla
u=\nabla\Delta^{-1}\left(\nabla\mathrm{div}u-\nabla\times\mathrm{curl}u\right)$,
we apply Calderon-Zygmund inequality to
get\be\label{Calderon-Zygmund}\begin{split} \|\nabla
u\|_{L^3}\lesssim&\|\mathrm{curl}u\|_{L^3}+\|\mathrm{div}u\|_{L^3}.
\end{split}
\ee Taking $\mathrm{curl}$ on both sides of (\ref{full N-S+1})$_2$,
we have \be\label{equation of curlu}
\mu\Delta(\mathrm{curl}u)=\mathrm{curl}(\rho \dot{u}). \ee By
(\ref{equation of curlu}), the $L^2$-estimates of the elliptic
equations and (\ref{blowup-2.1}), we have \beq\label{H 1 of curl u}
\begin{split}\|\nabla \mathrm{curl}u\|_{L^{2}} \les \|\rho \dot{u}\|_{L^{2}}\les
 \|\sqrt{\rho} \dot{u}\|_{L^{2}}.
\end{split}
\eeq By (\ref{equation of G}), (\ref{Calderon-Zygmund}) and (\ref{H
1 of curl u}), together with Sobolev inequality, we have
 \be\label{tu3}\begin{split} \|\nabla
u\|_{L^3}\lesssim&\|\mathrm{curl}u\|_{L^3}+\|G\|_{L^3}+\|\rho\theta\|_{L^3}
\lesssim\|\mathrm{curl}u\|_{H^1}+\|G\|_{H^1}+\|\rho\|_{L^6}\|\theta\|_{L^6}\\
\lesssim&\|\mathrm{curl}u\|_{L^2}+\|\mathrm{div}u\|_{L^2}+\|\nabla(\mathrm{curl}u)\|_{L^2}+\|\nabla G\|_{L^2}+\|\nabla\theta\|_{L^2}\\
\lesssim&\|\nabla(\mathrm{curl}u)\|_{L^2}+\|\nabla G\|_{L^2}+\|\nabla\theta\|_{L^2}+1\\
\lesssim&\|\sqrt{\rho}\u\|_{L^2}+\|\nabla\theta\|_{L^2}+1,
\end{split}
\ee where we have used (\ref{blowup-2.1}), (\ref{non-energy
inequality}) and (\ref{H 1 of u}). Substituting (\ref{tu3}) into
(\ref{VVV 2,2 and 2,6}) and using Young inequality, we obtain
\be\label{VVV 2,2 and 2,6+1}\begin{split}
VVV_{2,2}+VVV_{2,6}\lesssim\|\nabla
\u\|_{L^2}\left(\|\sqrt{\rho}\u\|_{L^2}+\|\nabla\theta\|_{L^2}+1\right)\|\nabla
\theta\|_{L^2}.
\end{split}
\ee  For $VVV_{2,3}$, $VVV_{2,4}$, $VVV_{2,7}$ and $VVV_{2,8}$,
using H\"older inequality, Sobolev inequality and Calderon-Zygmund
inequality, we have \be\label{VVV 3,4,7,8}\begin{split}
&VVV_{2,3}+VVV_{2,4}+VVV_{2,7}+VVV_{2,8}\\
\lesssim&\int_{\mathbb{R}^3} |\nabla u|^3|\theta|\lesssim\|\nabla
u\|_{L^{\frac{18}{5}}}^3\|\theta\|_{L^6}\lesssim\|\nabla
u\|_{L^{\frac{18}{5}}}^3\|\nabla \theta\|_{L^2}\\
\lesssim&\|\mathrm{curl} u\|_{L^{\frac{18}{5}}}^3\|\nabla
\theta\|_{L^2}+\|\mathrm{div} u\|_{L^{\frac{18}{5}}}^3\|\nabla
\theta\|_{L^2}\\
\lesssim&\|\mathrm{curl}u\|_{L^{\frac{18}{5}}}^3\|\nabla
\theta\|_{L^2}+\|G\|_{L^{\frac{18}{5}}}^3\|\nabla
\theta\|_{L^2}+\|\rho\theta\|_{L^{\frac{18}{5}}}^3\|\nabla
\theta\|_{L^2}.
\end{split}
\ee Using H\"older inequality again, together with
(\ref{blowup-2.1}), (\ref{non-energy inequality}), Sobolev
inequality, Gagliardo-Nirenberg inequality, (\ref{H 1 of u}),
(\ref{equation of G}) and (\ref{H 1 of curl u}), we get
\be\label{VVV 3,4,7,8+1}\begin{split}
&VVV_{2,3}+VVV_{2,4}+VVV_{2,7}+VVV_{2,8}\\
\lesssim&\|\mathrm{curl}u\|_{L^2}\|\nabla\mathrm{curl}u\|_{L^2}^2\|\nabla
\theta\|_{L^2}+\|G\|_{L^2}\|\nabla G\|_{L^2}^2\|\nabla
\theta\|_{L^2}+\|\rho\|_{L^9}^3\|\theta\|_{L^6}^3\|\nabla
\theta\|_{L^2}\\
\lesssim& \|\sqrt{\rho} \dot{u}\|_{L^{2}}^2\|\nabla
\theta\|_{L^2}+\|\nabla \theta\|_{L^2}^4.
\end{split}
\ee  For $VVV_{2,5}$ and $VVV_{2,9}$, using H\"older inequality,
Cauchy inequality, Sobolev inequality, Gagliardo-Nirenberg
inequality and (\ref{H 1 of u}), we have \be\label{VVV
5,9}\begin{split}
VVV_{2,5}+VVV_{2,9}\lesssim&\int_{\mathbb{R}^3}|\nabla
u|^2|u||\nabla\theta|\les\|\nabla
u\|_{L^4}^2\|u\|_{L^6}\|\nabla\theta\|_{L^3}\\ \les&
\|\nabla u\|_{L^4}^4+\|\nabla u\|_{L^2}^2\|\nabla\theta\|_{L^2}\|\nabla^2\theta\|_{L^2}\\
\le& C\|\nabla
u\|_{L^4}^4+C\|\nabla\theta\|_{L^2}\|\nabla^2\theta\|_{L^2}.
\end{split}
\ee From the standard elliptic estimates and (\ref{rho theta+ =}),
we have \be\label{H 2 of theta}\begin{split}
\|\nabla^2\theta\|_{L^2}\les&\|\rho \dot{\theta}\|_{L^2}+\|\rho\theta\mathrm{div}u\|_{L^2}+\|\nabla u\|_{L^4}^2\\
\lesssim&\|\sqrt{\rho} \dot{\theta}\|_{L^2}+\|\rho\|_{L^{12}}\|\theta\|_{L^6}\|\nabla u\|_{L^4}+\|\nabla u\|_{L^4}^2\\
\le& C\|\sqrt{\rho}\dot{\theta}\|_{L^2}+C\|\nabla
u\|_{L^4}^2+C\|\nabla\theta\|_{L^2}^2,
\end{split}
\ee where we have used H\"older inequality, (\ref{blowup-2.1}),
(\ref{non-energy inequality}), Sobolev inequality and Cauchy
inequality. Substituting (\ref{H 2 of theta}) into (\ref{VVV 5,9}),
and using Cauchy inequality, we have \be\label{VVV
5,9+1}\begin{split} VVV_{2,5}+VVV_{2,9}\le
\frac{1}{8}\|\sqrt{\rho}\dot{\theta}\|_{L^2}^2+C\|\nabla
u\|_{L^4}^4+C\|\nabla\theta\|_{L^2}^4+C.
\end{split}
\ee Substituting (\ref{VVV 2,2 and 2,6+1}), (\ref{VVV 3,4,7,8+1})
and (\ref{VVV 5,9+1}) into (\ref{VVV 2-1}), and using Cauchy
inequality, we have \be\label{VVV 2}\begin{split} VVV_2\le&
\frac{d}{dt}\int_{\mathbb{R}^3}\left(\frac{\mu}{2}\left|\nabla
u+(\nabla u)^\prime\right|^2+\lambda(\mathrm{div}u)^2\right)\theta
\\ &+C\|\nabla
\u\|_{L^2}\left(\|\sqrt{\rho}\u\|_{L^2}+\|\nabla\theta\|_{L^2}+1\right)\|\nabla
\theta\|_{L^2}+C\|\sqrt{\rho} \dot{u}\|_{L^{2}}^2\|\nabla
\theta\|_{L^2}\\ &+C\|\nabla
\theta\|_{L^2}^4+\frac{1}{8}\|\sqrt{\rho}\dot{\theta}\|_{L^2}^2+C\|\nabla
u\|_{L^4}^4+C.
\end{split}
\ee For $VVV_3$, using (\ref{VVV 5,9}) and (\ref{VVV 5,9+1}), we
have \be\label{VVV 3}\begin{split}
VVV_3\les&\int_{\mathbb{R}^3}|\nabla
u|^2|u||\nabla\theta|\le\frac{1}{8}\|\sqrt{\rho}\dot{\theta}\|_{L^2}^2+C\|\nabla
u\|_{L^4}^4+C\|\nabla\theta\|_{L^2}^4+C.
\end{split}
\ee For $VVV_4$, using H\"older inequality, Sobolev inequality,
Gagliardo-Nirenberg inequality, (\ref{H 1 of u}), (\ref{H 2 of
theta}) and Young inequality, we have \be\label{VVV 4}\begin{split}
VVV_4\les&
\|\Delta\theta\|_{L^2}\|u\|_{L^6}\|\nabla\theta\|_{L^3}\les
\|\Delta\theta\|_{L^2}\|\nabla
u\|_{L^2}\|\nabla\theta\|_{L^2}^\frac{1}{2}\|\nabla^2\theta\|_{L^2}^\frac{1}{2}\\
\les&
\|\nabla\theta\|_{L^2}^\frac{1}{2}\|\nabla^2\theta\|_{L^2}^\frac{3}{2}\le
\frac{1}{8}\|\sqrt{\rho}\dot{\theta}\|_{L^2}^2+C\|\nabla\theta\|_{L^2}^4
+C\|\nabla u\|_{L^4}^4+C.
\end{split}
\ee Putting (\ref{VVV 1}), (\ref{VVV 2}), (\ref{VVV 3}) and
(\ref{VVV 4}) into (\ref{dt nabla theta}), integrating the resulting
inequality over $[0,t]$ for $t\in(0,T^*)$, and using Cauchy
inequality, (\ref{blowup-2.1}), (\ref{u 2nabla u 2}) and (\ref{H 1
of u}), we have \be\label{dt nabla theta+1a}\begin{split}
\int_0^t\int_{\mathbb{R}^3}\rho
|\dot{\theta}|^2+\int_{\mathbb{R}^3}|\nabla \theta|^2
\le&C\int_0^t\|\nabla u\|_{L^4}^4+C\int_{\mathbb{R}^3}|\nabla
u|^2|\theta|+\varepsilon\int_0^t\|\nabla
\u\|_{L^2}^2\\&+C_\varepsilon\int_0^t
\left(\|\sqrt{\rho}\u\|_{L^2}^2+\|\nabla
\theta\|_{L^2}^2+1\right)\|\nabla \theta\|_{L^2}^2+C.
\end{split}
\ee For the second term of the right hand side of (\ref{dt nabla
theta+1a}), we have \be\label{D}\begin{split}
C\int_{\mathbb{R}^3}|\nabla u|^2|\theta|\lesssim&\|\nabla u\|_{L^\frac{12}{5}}^2\|\theta\|_{L^6}\lesssim\|\mathrm{curl} u\|_{L^\frac{12}{5}}^2\|\nabla\theta\|_{L^2}+
\|\mathrm{div} u\|_{L^\frac{12}{5}}^2\|\nabla\theta\|_{L^2}\\
\lesssim&\left(\|\mathrm{curl}u\|_{L^2}^{\frac{3}{2}}\|\nabla\mathrm{curl}u\|_{L^2}^{\frac{1}{2}}+
\|G\|_{L^2}^{\frac{3}{2}}\|\nabla
G\|_{L^2}^{\frac{1}{2}}\right)\|\nabla\theta\|_{L^2}+\|\rho\theta\|_{L^\frac{12}{5}}^2\|\nabla\theta\|_{L^2}\\
\lesssim&
\|\sqrt{\rho}\u\|_{L^2}^{\frac{1}{2}}\|\nabla\theta\|_{L^2}+\|\rho\theta\|_{L^2}^{\frac{3}{2}}\|\rho\theta\|_{L^6}^{\frac{1}{2}}\|\nabla\theta\|_{L^2}
\\
\lesssim&\|\sqrt{\rho}\u\|_{L^2}^{\frac{1}{2}}\|\nabla\theta\|_{L^2}+\|\nabla\theta\|_{L^2}^{\frac{3}{2}}
\le\frac{1}{2}\|\nabla\theta\|_{L^2}^2+C\|\sqrt{\rho}\u\|_{L^2}+C,
\end{split}
\ee where we have used H\"older inequality, Calderon-Zygmund
inequality, Gagliardo-Nirenberg inequality, (\ref{blowup-2.1}),
(\ref{H 1 of u}), (\ref{equation of G}), (\ref{H 1 of curl u}),
Sobolev inequality and Young inequality. Substituting (\ref{D}) into
(\ref{dt nabla theta+1a}), we have \be\label{dt nabla
theta+1}\begin{split} \int_0^t\int_{\mathbb{R}^3}\rho
|\dot{\theta}|^2+\int_{\mathbb{R}^3}|\nabla \theta|^2
\le&C\int_0^t\|\nabla u\|_{L^4}^4+C
\|\sqrt{\rho}\u\|_{L^2}+\varepsilon\int_0^t\|\nabla \u\|_{L^2}^2
\\&+C_\varepsilon\int_0^t
\left(\|\sqrt{\rho}\u\|_{L^2}^2+\|\nabla
\theta\|_{L^2}^2+1\right)\|\nabla \theta\|_{L^2}^2+C.
\end{split}
\ee Multiplying (\ref{dt nabla theta+1}) by $2C$ and adding the
resulting inequality into (\ref{dt rho u t-2}), we have \bex
\begin{split}
&C\int_0^t\int_{\mathbb{R}^3}\rho |\dot{\theta}|^2+2C\int_{\mathbb{R}^3}|\nabla \theta|^2+\int_{\mathbb{R}^3}\rho|\u|^2+\int_0^t\int_{\mathbb{R}^3}|\nabla\u|^2\\
\leq&2C^2\int_0^t\|\nabla u\|_{L^4}^4+2C^2\|\sqrt{\rho}\u\|_{L^2}+2\varepsilon C\int_0^t\|\nabla \u\|_{L^2}^2\\
&+2CC_\varepsilon\int_0^t \left(\|\sqrt{\rho}\u\|_{L^2}^2+\|\nabla
\theta\|_{L^2}^2+1\right)\|\nabla \theta\|_{L^2}^2+2C^2.
\end{split}
\eex Taking $\varepsilon$ sufficiently small, together with Cauchy
inequality, we have \be\label{non-sum}
\begin{split}
&\int_{\mathbb{R}^3}(|\nabla
\theta|^2+\rho|\u|^2)+\int_0^t\int_{\mathbb{R}^3}(\rho
|\dot{\theta}|^2+|\nabla\u|^2) \\ \les&\int_0^t\|\nabla u\|_{L^4}^4+
\int_0^t \left(\|\sqrt{\rho}\u\|_{L^2}^2+\|\nabla
\theta\|_{L^2}^2+1\right)\|\nabla \theta\|_{L^2}^2+1.
\end{split}
\ee For the first term of the right hand side of (\ref{non-sum}),
similar to (\ref{tu3}), we have \be\label{non-sum+1}
\begin{split}
\int_0^t\|\nabla u\|_{L^4}^4
\les&\int_0^t\|\mathrm{curl}u\|_{L^4}^4+\int_0^t\|G\|_{L^4}^4+\int_0^t\|\nabla\theta\|_{L^2}^4\\
\les&\int_0^t\|\mathrm{curl}u\|_{L^2}\|\nabla\mathrm{curl}u\|_{L^2}^3+\int_0^t\|G\|_{L^2}\|\nabla
G\|_{L^2}^3+\int_0^t\|\nabla\theta\|_{L^2}^4\\
\les&\int_0^t\|\sqrt{\rho}\u\|_{L^2}^3+\int_0^t\|\nabla\theta\|_{L^2}^4.
\end{split}
\ee
 By (\ref{non-sum}), (\ref{non-sum+1}) and Cauchy inequality, we have \be\label{non-sum+2}
\begin{split}
&\int_{\mathbb{R}^3}(|\nabla
\theta|^2+\rho|\u|^2)+\int_0^t\int_{\mathbb{R}^3}(\rho
|\dot{\theta}|^2+|\nabla\u|^2)\\
\les&\int_0^t\|\sqrt{\rho}\u\|_{L^2}^3+ \int_0^t
\left(\|\sqrt{\rho}\u\|_{L^2}^2+\|\nabla
\theta\|_{L^2}^2+1\right)\|\nabla \theta\|_{L^2}^2+1\\ \les&\int_0^t
\left(\|\sqrt{\rho}\u\|_{L^2}^2+\|\nabla
\theta\|_{L^2}^2+1\right)\left(\|\nabla
\theta\|_{L^2}^2+\|\sqrt{\rho}\u\|_{L^2}^2\right)+1.
\end{split}
\ee By (\ref{blowup-2.1}), (\ref{u 2nabla u 2}), (\ref{H 1 of u}),
we have
$$
\int_0^t \left(\|\sqrt{\rho}\u\|_{L^2}^2+\|\nabla
\theta\|_{L^2}^2\right)\le C,
$$ for any $t\in(0,T^*)$. This, together with (\ref{non-sum+2}) and
Gronwall inequality, deduces (\ref{H 2 of u}).
\endpf

\begin{corollary}Under the conditions of Theorem \ref{blowup-th:1.1} and (\ref{blowup-2.1}), it holds that for any $T\in(0,T^*)$
\be\label{u infty}\sup\limits_{0\le t\le T}\left(\|\nabla
G\|_{L^2}+\|\nabla \mathrm{curl}u\|_{L^2}+\|\nabla
u\|_{L^6}+\|u\|_{L^\infty}\right)+\int_0^T\int_{\mathbb{R}^3}|\nabla^2\theta|^2\le
C.\ee
\end{corollary}
\pf The detailed proof of the lemma could be found in \cite{Wen-Zhu
3} (see Corollary 4.5 therein).
\endpf

\begin{lemma}\label{blowup-le: H 2 of theta}Under the conditions of Theorem \ref{blowup-th:1.1} and (\ref{blowup-2.1}), it holds that for any $t\in(0,T^*)$
\be\label{non-H 2 of theta}
\int_{\mathbb{R}^3}\rho|\theta_t|^2+\int_0^t\int_{\mathbb{R}^3}|\nabla\theta_t|^2\le
C. \ee
\end{lemma}
\pf By \cite{Wen-Zhu 3}, we have \be\label{dt H 2 of
theta}\begin{split}
&\frac{1}{2}\frac{d}{dt}\int_{\mathbb{R}^3}\rho|\theta_t|^2+\int_{\mathbb{R}^3}|\nabla\theta_t|^2\\
=&-
\int_{\mathbb{R}^3}\rho_t\left(\frac{\theta_t}{2}+u\cdot\nabla\theta+\theta\mathrm{div}u\right)\theta_t-\int_{\mathbb{R}^3}\rho(
u_t\cdot\nabla\theta+u\cdot\nabla\theta_t+\theta_t\mathrm{div}u)\theta_t\\&-\int_{\mathbb{R}^3}\rho\theta\mathrm{div}u_t\theta_t+\mu\int_{\mathbb{R}^3}\left(\nabla
u+(\nabla u)^\prime\right):\left(\nabla u_t+(\nabla
u_t)^\prime\right)\theta_t+2\lambda\int_{\mathbb{R}^3}\mathrm{div}u\mathrm{div}u_t\theta_t\\
=&\sum\limits_{i=1}^5VVVI_i,
\end{split}
\ee and \be\label{VVVI 1}\begin{split} VVVI_1=&
-\int_{\mathbb{R}^3}\rho
u\cdot\nabla\theta_t\left(\frac{\theta_t}{2}+u\cdot\nabla\theta+\theta\mathrm{div}u\right)-\int_{\mathbb{R}^3}\rho
u\cdot\frac{\nabla\theta_t}{2}\theta_t\\&-\int_{\mathbb{R}^3}\rho
u\cdot\left(\nabla
(u\cdot\nabla)\theta+u\cdot\nabla\nabla\theta\right)\theta_t-\int_{\mathbb{R}^3}\rho
u\cdot\left(\nabla\theta\mathrm{div}u+\theta\nabla\mathrm{div}u\right)\theta_t\\=&\sum\limits_{i=1}^4VVVI_{1,i}.
\end{split}
\ee For $VVVI_{1,1}$, we have \be\label{VVVI 1,1}\begin{split}
VVVI_{1,1} \le&
\frac{1}{24}\|\nabla \theta_t\|_{L^2}^2+C\|\sqrt{\rho}\theta_t\|_{L^2}^2+C\|\nabla\theta\|_{L^2}^2+C\|\sqrt[4]{\rho}\theta\|_{L^4}^2\|\nabla u\|_{L^4}^2\\
\le&\frac{1}{24}\int_{\mathbb{R}^3}|\nabla\theta_t|^2+C\int_{\mathbb{R}^3}\rho|\theta_t|^2
+C,
\end{split}
\ee where we have used Cauchy inequality, (\ref{blowup-2.1}),
(\ref{H 1 of u}), (\ref{H 2 of u}) and (\ref{u infty}).\\ For
$VVVI_{1,2}$, using Cauchy inequality, (\ref{blowup-2.1}) and
(\ref{u infty}) again, we have \be\label{VVVI 1,2}\begin{split}
VVVI_{1,2}\le&\frac{1}{24}\int_{\mathbb{R}^3}|\nabla\theta_t|^2+C\int_{\mathbb{R}^3}\rho
|\theta_t|^2.
\end{split}
\ee For $VVVI_{1,3}$, by \cite{Wen-Zhu 3} (see (4.65) therein), we
have \be\label{VVVI 1,3}\begin{split} VVVI_{1,3}
\les&\int_{\mathbb{R}^3}\rho
|\theta_t|^2+\int_{\mathbb{R}^3}|\nabla^2\theta|^2+1.
\end{split}
\ee For $VVVI_{1,4}$, integrating by parts, we have \be\label{VVVI
1,4a}\begin{split} VVVI_{1,4}
=&-\int_{\mathbb{R}^3}\rho u\cdot\nabla\theta\mathrm{div}u\theta_t-\displaystyle\frac{1}{2\mu+\lambda}\int_{\mathbb{R}^3}\rho\theta u\cdot\nabla G\theta_t
-\displaystyle\frac{1}{2\mu+\lambda}\int_{\mathbb{R}^3}\rho\theta u\cdot\nabla (\rho\theta)\theta_t\\
=&-\int_{\mathbb{R}^3}\rho
u\cdot\nabla\theta\mathrm{div}u\theta_t-\displaystyle\frac{1}{2\mu+\lambda}\int_{\mathbb{R}^3}\rho\theta
u\cdot\nabla G\theta_t
\\
&+\displaystyle\frac{1}{2(2\mu+\lambda)}\int_{\mathbb{R}^3}\rho^2\theta^2
u\cdot\nabla\theta_t+\displaystyle\frac{1}{2(2\mu+\lambda)}\int_{\mathbb{R}^3}\rho^2\theta^2\mathrm{div}u\theta_t
.
\end{split}
\ee Furthermore, we can get \be\label{VVVI 1,4}\begin{split}
VVVI_{1,4}\lesssim&\|\sqrt{\rho}\theta_t\|_{L^2}\|\nabla
u\|_{L^6}\|\nabla \theta\|_{L^3} +\|\nabla
G\|_{L^2}\|\theta_t\|_{L^6}\|\theta\|_{L^6}\|\rho\|_{L^{6}}\\&
+\|\theta\|_{L^6}^2\|\rho\|_{L^{12}}^2\|\nabla \theta_t\|_{L^2}+\|\nabla u\|_{L^3}\|\theta_t\|_{L^6}\|\theta\|_{L^6}^2\|\rho\|_{L^{12}}^2\\
\leq&\frac{1}{24}\|\nabla
\theta_t\|_{L^2}^2+C\|\nabla^2\theta\|_{L^2}^2+C\|\sqrt{\rho}\theta_t\|_{L^2}^2+C,
\end{split}
\ee where we have used H\"older inequality, Sobolev inequality,
(\ref{blowup-2.1}), (\ref{non-energy inequality}), (\ref{H 1 of u}),
(\ref{H 2 of u}), (\ref{u infty}) and Cauchy inequality.
Substituting (\ref{VVVI 1,1}), (\ref{VVVI 1,2}), (\ref{VVVI 1,3})
and (\ref{VVVI 1,4}) into (\ref{VVVI 1}), we have \be\label{VVVI
1+1}\begin{split}
VVVI_1\le&\frac{1}{8}\int_{\mathbb{R}^3}|\nabla\theta_t|^2
+C\int_{\mathbb{R}^3}\rho
|\theta_t|^2+C\int_{\mathbb{R}^3}|\nabla^2\theta|^2+C.
\end{split}
\ee For $VVVI_2$, \be\label{VVVI 2}\begin{split} VVVI_2
=&-\int_{\mathbb{R}^3}\rho
\u\cdot\nabla\theta\theta_t+\int_{\mathbb{R}^3}\rho
(u\cdot\nabla)u\cdot\nabla\theta\theta_t
-\int_{\mathbb{R}^3}\rho u\cdot\nabla\theta_t\theta_t-\int_{\mathbb{R}^3}\rho|\theta_t|^2\mathrm{div}u\\
\lesssim&\|\sqrt{\rho}\u\|_{L^2}\|\theta_t\|_{L^6}\|\nabla \theta\|_{L^3}+\|\sqrt{\rho}\theta_t\|_{L^2}\|\nabla u\|_{L^6}\|\nabla \theta\|_{L^3}+\|\sqrt{\rho}\theta_t\|_{L^2}\|\nabla \theta_t\|_{L^2}\\&+\|\sqrt{\rho}\theta_t\|_{L^2}\|\theta_t\|_{L^6}\|\nabla u\|_{L^3}\\
\leq&\frac{1}{8}\int_{\mathbb{R}^3}|\nabla\theta_t|^2
+C\int_{\mathbb{R}^3}\rho
|\theta_t|^2+C\int_{\mathbb{R}^3}|\nabla^2\theta|^2+C,
\end{split}
\ee
where we have used H\"older inequality, Sobolev inequality, (\ref{blowup-2.1}), (\ref{H 1 of u}), (\ref{H 2 of u}) and (\ref{u infty}).\\
For $VVVI_3$, integrating by parts, we have \be\label{VVVI
3}\begin{split}
VVVI_3=&-\int_{\mathbb{R}^3}\rho\theta\mathrm{div}\u\theta_t+\int_{\mathbb{R}^3}\rho\theta\mathrm{div}(u\cdot\nabla u)\theta_t\\
=&-\int_{\mathbb{R}^3}\rho\theta\mathrm{div}\u\theta_t+\int_{\mathbb{R}^3}\rho\theta\nabla u:(\nabla u)'\theta_t
+\int_{\mathbb{R}^3}\rho\theta u\cdot\nabla \mathrm{div}u\theta_t\\
=&-\int_{\mathbb{R}^3}\rho\theta\mathrm{div}\u\theta_t+\int_{\mathbb{R}^3}\rho\theta\nabla
u:(\nabla u)'\theta_t
+\displaystyle\frac{1}{2\mu+\lambda}\int_{\mathbb{R}^3}\rho\theta
u\cdot\nabla G\theta_t
\\
&-\displaystyle\frac{1}{2\mu+\lambda}\int_{\mathbb{R}^3}\frac{\rho^2}{2}\theta^2\mathrm{div}u\theta_t-\displaystyle\frac{1}{2\mu+\lambda}\int_{\mathbb{R}^3}\frac{\rho^2}{2}\theta^2u\cdot\nabla\theta_t.
\end{split}
\ee Furthermore, using H\"older inequality, Sobolev inequality,
(\ref{blowup-2.1}), (\ref{non-energy inequality}), (\ref{H 1 of u}),
(\ref{H 2 of u}), (\ref{u infty}) and Young inequality, we have
\be\label{VVVI 3+1}\begin{split} VVVI_3 \lesssim&\|\nabla
\u\|_{L^2}\|\theta_t\|_{L^6}\|\theta\|_{L^6}\|\rho\|_{L^{6}}+
\|\sqrt{\rho}\theta\|_{L^2}\|\nabla
u\|_{L^6}^2\|\theta_t\|_{L^6}\\&+ \|\nabla
G\|_{L^2}\|\theta_t\|_{L^6}\|\theta\|_{L^6}\|\rho\|_{L^{6}}
+\|\nabla u\|_{L^3}\|\theta_t\|_{L^6}\|\theta\|_{L^6}^2\|\rho\|_{L^{12}}^2\\&+\|\theta\|_{L^6}^2\|\rho\|_{L^{12}}^2\|\nabla \theta_t\|_{L^2}\\
\leq&\frac{1}{8}\|\nabla \theta_t\|_{L^2}^2+C\|\nabla
\u\|_{L^2}^2+C.
\end{split}
\ee Similar to $VVV_2$, for $VVVI_4$ and $VVVI_5$, we deduce
\be\label{VVVI 4 and VVVI 5}
\begin{split}
VVVI_4+VVVI_5\le&C\|\nabla\u\|_{L^2}\|\nabla
u\|_{L^3}\|\theta_t\|_{L^6}+C\int_{\mathbb{R}^3}|\nabla
u|^3|\theta_t|+C\int_{\mathbb{R}^3}|\nabla
u|^4+\frac{1}{16}\int_{\mathbb{R}^3}|\nabla\theta_t|^2\\
\le&C\|\nabla\u\|_{L^2}\|\nabla\theta_t\|_{L^2}+C\|\nabla
u\|_{L^\frac{18}{5}}^3\|\nabla\theta_t\|_{L^2}+\frac{1}{16}\int_{\mathbb{R}^3}|\nabla\theta_t|^2
+C\\
\le&\frac{1}{8}\int_{\mathbb{R}^3}|\nabla\theta_t|^2+C\int_{\mathbb{R}^3}|\nabla\u|^2+C,
\end{split}
\ee
where we have used H\"older inequality, integration by parts, Cauchy inequality, (\ref{H 1 of u}), (\ref{u infty}), the interpolation inequality and Sobolev inequality.\\
Putting (\ref{VVVI 1+1}), (\ref{VVVI 2}), (\ref{VVVI 3+1}) and
(\ref{VVVI 4 and VVVI 5}) into (\ref{dt H 2 of theta}), we have
\be\label{dt H 2 of theta+1}\begin{split}
\frac{d}{dt}\int_{\mathbb{R}^3}\rho|\theta_t|^2+\int_{\mathbb{R}^3}|\nabla\theta_t|^2\le&
C\int_{\mathbb{R}^3}\rho|\theta_t|^2
+C\int_{\mathbb{R}^3}(|\nabla\u|^2+|\nabla^2\theta|^2)+C.
\end{split}
\ee By (\ref{dt H 2 of theta+1}), (\ref{H 2 of u}), (\ref{u infty})
and Gronwall inequality, we complete the proof of Lemma
\ref{blowup-le: H 2 of theta}.
\endpf
\begin{corollary}Under the conditions of Theorem
\ref{blowup-th:1.1} and (\ref{blowup-2.1}), it holds that for any
$t\in( 0,T^*)$ \be\label{cor:H 2 of theta}
\int_{\mathbb{R}^3}|\nabla^2\theta|^2\le C. \ee
\end{corollary}
\pf It follows from (\ref{full N-S+1})$_3$, (\ref{blowup-2.1}),
(\ref{non-energy inequality}), (\ref{H 1 of u}), (\ref{H 2 of u}),
(\ref{u infty}), (\ref{non-H 2 of theta}) and the interpolation
inequality that \bex\begin{split} \|\nabla^2\theta\|_{L^2}\les&
\|\rho \theta_t\|_{L^2}+ \|\rho
u\cdot\nabla\theta\|_{L^2}+\|\rho\theta\mathrm{div}u\|_{L^2}+\|\nabla
u\|_{L^4}^2\\ \les&\|\sqrt{\rho} \theta_t\|_{L^2}+
\|\nabla\theta\|_{L^2}+\|\rho\|_{L^6}\|\theta\|_{L^6}\|\mathrm{div}u\|_{L^6}+1\\
\le&C.
\end{split}\eex
\endpf\\

By (\ref{H 2 of u}), (\ref{cor:H 2 of theta}) and Sobolev
inequality, we get the following corollary which is the desired one,
i.e., (\ref{non-uniform_est1}).
\begin{corollary}Under the conditions of Theorem \ref{blowup-th:1.1} and (\ref{blowup-2.1}), it holds that for any $t\in(0,T^*)$
\be\label{cor:H 2 of u} \|\theta\|_{L^\infty(0,t;L^\infty)}\le C.
\ee
\end{corollary}

\section*{Acknowledgements} The research was supported by the
National Natural Science Foundation of China $\#$10625105,
$\#$11071093, the PhD specialized grant of the Ministry of
Education of China $\#$20100144110001, and the Special Fund for
Basic Scientific Research  of Central Colleges $\#$CCNU10C01001.


\vskip 1cm

\addcontentsline{toc}{section}{\\References}

\begin{thebibliography}{99}

\bibitem{Bresch-Desjardins} D. Bresch, B. Desjardins, {\em On the existence of global weak solutions to the Navier-Stokes
equations for viscous compressible and heat conducting fluids,} J.
Math. Pures Appl., 87(2007), 57-90.

\bibitem{Cho-Choe-Kim} Y. Cho, H.J. Choe, H. Kim, {\em Unique solvability of the initial boundary value problems
for compressible viscous fluids,} J. Math. Pures Appl., 83(2004),
243-275.

\bibitem{Cho-Kim} Y. Cho, H. Kim, {\em On classical solutions of the compressible Navier-Stokes equations with
nonnegative initial densities,} Manuscripta Math., 120(2006),
91-129.

\bibitem{cho-Kim: perfect gas} Y. Cho, H. Kim, {\em Existence results for viscous polytropic fluids
with vacuum,} J. Differential Equations, 228(2006), 377-411.

\bibitem{Choe-Kim: Strong} H.J. Choe, H. Kim, {\em Strong solutions of the Navier-Stokes equations for isentropic
compressible fluids,} J. Differential Equations, 190(2003), 504-523.

\bibitem{Choe-Kim: Radially} H.J. Choe, H. Kim, {\em Global existence of the radially symmetric solutions of
the Navier-Stokes equations for the isentropic compressible fluids,}
Math. Methods Appl. Sci., 28(2005), 1-28.

\bibitem{Ding-Wen-Yao-Zhu} S.J. Ding, H.Y. Wen, L. Yao, C.J. Zhu, {\em Global classical spherically symmetric solution of compressible isentropic
Navier-Stokes equations with vacuum,} SIAM J. Math. Anal., 2011, to
appear.

\bibitem{Ding-Wen-Zhu} S.J. Ding, H.Y. Wen, C.J. Zhu, {\em Global classical large solutions to 1D compressible Navier-Stokes equations
with density-dependent viscosity and vacuum,} J. Differential
Equations, 251(2011), 1696-1725.

\bibitem{Evans} L.C. Evans, {\em Partial Differential Equations,}
Graduate Studies in Mathematics, Amer. Math. Soc., Providence, Rhode
Island, Vol. 19, 1998.

\bibitem{Fan-Jiang-Ni} J. Fan, S. Jiang, G. Ni, {\em Uniform boundedness
of the radially symmetric solutions of the Navier-Stokes equations
for isentropic compressible fluids,} Osaka J. Math., 46(2009),
863-876.

\bibitem{Fan-Jiang-Ou: Blow-up criterion} J. Fan, S. Jiang, Y. Ou,
{\em A blow-up criterion for compressible viscous heat-conductive
flows,} Ann. I. H. Poincar$\acute{\mathrm{e}}$-AN, 27(2010),
337-350.

\bibitem{Fang-Zi-Zhang} D.Y. Fang, R.Z. Zi, T. Zhang, {\em A blow-up criterion for two dimensional
compressible viscous heat-conductive flows}. arXiv:1107.4663v1
[math.AP] 23 Jul 2011.

\bibitem{Feireisl1} E. Feireisl, {\em On the motion of a viscous, compressible and heat conducting fluid,} Indiana Univ. Math. J., 53(2004),
1705-1738.

\bibitem{Feireisl-book} E. Feireisl, {\em Dynamics of Viscous Compressible Fluids,} Oxford Univ. Press, Oxford, 2004.

\bibitem{Feireisl2} E. Feireisl, A. Novotn$\acute{\mathrm{y}}$, H.
Petzeltov$\mathrm{\acute{a}}$, {\em On the existence of globally
defined weak solutions to the Navier-Stokes equations,} J. Math.
Fluid Mech., 3(2001), 358-392.

\bibitem{Frid} H. Frid, V. Shelukhin, {\em Vanishing shear viscosity in the equations
of compressible fluids for the flows with the cylinder symmetry,}
SIAM J. Math. Anal., 31(2000), 1144-1156.

\bibitem{Guo-Zhu} Z.H. Guo, C.J. Zhu, {\em Global weak solutions and asymptotic behavior
to 1D compressible Navier-Stokes equations with density-dependent
viscosity and vacuum,} J. Differential Equations, 248(2010),
2768-2799.

\bibitem{Hoff} D. Hoff, H. Jenssen, {\em Symmetric nonbarotropic flows
with large data and forces,} Arch. Rational Mech. Anal., 173
(2004), 297-343.

\bibitem{Huang-Li} X. Huang, J. Li, {\em Global classical and weak solutions to the
three-dimensional full compressible Navier-Stokes system with
vacuum and large oscillations,} arXiv:1107.4655v2 [math-ph] 26 Jul
2011.

\bibitem{Huang-Li-Xin} X. Huang, J. Li, Z. Xin, {\em
Global well-posedness of classical solutions with large
oscillations and vacuum to the three-dimensional isentropic
compressible Navier-Stokes equations,} arXiv:1004.4749v2 [math-ph]
8 Jul 2010.

\bibitem{Itaya} N. Itaya, {\em On the Cauchy problem for the system of fundamental equations describing the movement of compressible viscous
fluid,} Kodai Math. Sem. Rep., 23(1971), 60-120.

\bibitem{Jiang1} S. Jiang, {\em Global spherically symmetric solutions
to the equations of a viscous polytropic ideal gas in an exterior
domain,} Comm. Math. Phys., 178(1996), 339-374.

\bibitem{Jiang:Math Nachr} S. Jiang, {\em Global smooth solutions of the equations of a viscous,
heat-conducting one-dimensional gas with density-dependent
viscosity,} Math. Nachr., 190(1998), 169-183.

\bibitem{Jiang2} S. Jiang, {\em Large-time behavior of solutions to the equations
of a viscous polytropic ideal gas,} Annali di Matematica pura ed
applicata, (IV)(1998), 253-275.

\bibitem{Jiang3} S. Jiang, {\em Large-time behavior of solutions to the equations of
a one-dimensional viscous polytropic ideal gas in unbounded
domains,} Comm. Math. Phys., 200(1999), 181-193.

\bibitem{Jiang-Zhang} S. Jiang, J.W. Zhang, {\em Boundary layers for the Navier-Stokes equations
of compressible heat-conducting flows with cylindrical symmetry,}
SIAM J. Math. Anal., 41(2009), 237-268.

\bibitem{Jiang-Zhang1} S. Jiang, P. Zhang, {\em On spherically symmetric solutions
of the compressible isentropic Navier-Stokes equations,} Comm. Math.
Phys., 215(2001), 559-581.

\bibitem{Jiang-Zhang:weak solutions} S. Jiang, P. Zhang, {\em Global weak solutions to the Navier-Stokes equations for a 1D viscous polytropic ideal
gas,} Quart. Appl. Math., 61(2003), 435-449.

\bibitem{Kawohl} B. Kawohl, {\em Global existence of large solutions to initial boundary value
problems for a viscous, heat-conducting, one-dimensional real gas,}
J. Differential Equations, 58(1985), 76-103.

\bibitem{Kazhikhov-Shelukhi} A.V. Kazhikhov, V.V. Shelukhi, {\em Unique global solution with respect to time of
initial-boundary value problems for one-dimensional equations of a
viscous gas,} Prikl. Mat. Meh., 41(1977), 282-291.


\bibitem{Lions2} P.L. Lions, Mathematical Topics in Fluid Mechanics,
Vol. II, Compressible Models, Clarendon Press, Oxford, 1998.

\bibitem{Liu-Xin-Yang} T.P. Liu, Z. Xin, T. Yang, {\em Vacuum states of compressible flow,}
Discrete Contin. Dyn. Syst., 4(1998), 1-32.

\bibitem{Luo} T. Luo, Z.P. Xin, T. Yang, {\em Interface behavior of
compressible Navier-Stokes equations with vacuum,} SIAM J. Math.
Anal., 31(2000), 1175-1191.

\bibitem{Matsumura-Nishida: Kyoto Un} A. Matsumura, T. Nishida, {\em The initial value problem for the equations of motion of
viscous and heat-conductive gases,} J. Math. Kyoto Univ., 20(1980),
67-104.

\bibitem{Matsumura-Nishida: CMP} A. Matsumura, T. Nishida, {\em The initial boundary value problems for the equations of
motion of compressible and heat-conductive fluids,} Comm. Math.
Phys., 89(1983), 445-464.

\bibitem{Qin-Yao} X. Qin, Z. Yao,
{\em Global smooth solutions of the compressible Navier-Stokes
equations with density-dependent viscosity,} J. Differential
Equations, 244(2008), 2041-2061.

\bibitem{salvi}
R. Salvi, I. Stra$\breve{\mathrm{s}}$kraba, {\em Global existence
for viscous compressible fluids and their behavior as
$t\rightarrow\infty$, } J. Fac. Sci. Univ. Tokyo, Sect. IA Math.,
40(1993), 17-51.

\bibitem{Simon} J. Simon, {\em Nonhomogeneous viscous incompressible fluids: existence of vecocity, density and
pressure}, SIAM J. Math. Anal., 21(1990), 1093-1117.


\bibitem{Sun-Wang-Zhang 1} Y.Z. Sun, C. Wang, Z.F. Zhang, {\em A Beale-Kato-Majda criterion for three
dimensional compressible viscous heat-conductive Flows}. Arch.
Rational Mech. Anal., 201(2011), 727-742.

\bibitem{Tani} A. Tani, {\em On the first initial-boundary value problem of compressible viscous fluid
motion,} Publ. Res. Inst. Math. Sci. Kyoto Univ., 13(1977), 193-253.

\bibitem{Vong-Yang-Zhu} S.W. Vong, T. Yang,  C.J. Zhu, {\em Compressible Navier-Stokes equations
with degenerate viscosity coefficient and vacuum (II),} J.
Differential Equations, 192(2003), 475-501.

\bibitem{Wen-Zhu} H.Y. Wen, C.J. Zhu, {\em Global classical large solutions to Navier-Stokes
equations for viscous compressible and heat conducting fluids with
vacuum,} arXiv:1103.1421v1 [math.AP] 8 Mar 2011, 1-38.

\bibitem{Wen-Zhu 3} H.Y. Wen, C.J. Zhu. {\em Blow-up criterions of strong solutions to 3D compressible Navier-Stokes equations with
vacuum},
arXiv:1111.2657v1 [math.AP] 11 Nov 2011.

\bibitem{Xin} Z. Xin, {\em Blowup of smooth solutions to the compressible Navier-Stokes equation with
compact density,} Comm. Pure Appl. Math., 51(1998), 229-240.

\bibitem{Yang} T. Yang, {\em Singular behavior of vacuum states for compressible fluids,} Journal of Computational and Applied Mathematics, 190(2006), 211-231.

\bibitem{Yang-Yao-Zhu} T. Yang, Z.A. Yao, C.J. Zhu, {\em Compressible Navier-Stokes equations
with density-dependent viscosity and vacuum,} Comm. Partial
Differential Equations, 26(2001), 965-981.

\bibitem{Yang-Zhu} T. Yang, C.J. Zhu, {\em Compressible Navier-Stokes equations
with degenerate viscosity coefficient and vacuum,} Comm. Math.
Phys., 230(2002), 329-363.

\bibitem{Zhu} C.J. Zhu, {\em Asymptotic behavior of compressible Navier-Stokes
equations with density-dependent viscosity and vacuum,} Comm. Math.
Phys., 293(2010), 279-299.
\end{thebibliography}
\end{document}